\documentclass[11pt,a4thesis]{report}

\usepackage{amsmath}
\usepackage{amssymb}
\usepackage[usenames]{color}
\usepackage{theorem}
\usepackage [applemac] {inputenc} 
\usepackage{enumerate}

\newtheorem{theorem}{Theorem}[section]
\newtheorem{proposition}[theorem]{Proposition}
\newtheorem{lemma}[theorem]{Lemma}
\newtheorem{corollary}[theorem]{Corollary}
\newtheorem{definition}[theorem]{Definition}
\newtheorem{remark}[theorem]{Remark}

\newenvironment{proof}{\noindent\textbf{Proof}}
                      {$\Box$\vskip\theorempostskipamount}

\begin{document}

\title{\textbf{Algorithmic and combinatorial methods for enumerating the relators of a group presentation}}

\author{Carmelo Vaccaro}

\date{}

\maketitle

\begin{abstract}

 The main achievement of this thesis is an algorithm which given a finite group presentation and natural numbers $n$ and $k$, computes all the relators of length and area up to $n$ and $k$ respectively. The complexity of this algorithm is less by a factor of at least $a^{k^2} a^{kn}$ than those of the classical methods using van Kampen diagrams ($a$ is a constant).

\end{abstract}

\tableofcontents

\newpage

\chapter*{Introduction}

\addcontentsline{toc}{chapter}{Introduction}

The starting point of this thesis is the following: we have a finite group presentation and we have not additional information about the group presented. This is a very natural point of view since one can give a finite presentation by simply listing a finite set of words on a given alphabet; the words and the alphabet represent respectively the defining relators and the generators of the presentation.

A question that naturally arises is: what can be said about the group defined by this presentation? It is easy to convince oneself that the latter is a very hard question to answer in general due to the arbitrariness of the defining relators. Therefore it does not come as a surprise the fact that all the important problems concerning finite presentations have resulted to be unsolvable.

The situation is much better if the presentation verifies some strong hypothesis, like for instance a \textit{small cancellation} condition. When this is not the case one hope comes from the \textit{Knuth-Bendix algorithm} (see \cite{Sims} or \cite{Holt}). But as it is known this algorithm does not work for every presentation, in particular it can fail or run forever. 

\smallskip

Many decision problems for groups are based on the decidability of certain sets of relators. For instance the Word Problem is solvable if for every natural number $n$ one can decide the set of relators of length $n$; or the Conjugacy Problem is solvable if for every natural numbers $n_1$ and $n_2$ one can decide the set of relators of the form $u_1 \, a \, u_2 \, a^{-1}$ where $u_1$ and $u_2$ are words of length $n_1$ and $n_2$ respectively.


It must be observed that  for a finite presentation the set of relators is always recursively enumerable. 

If the presentation does not verify some special condition or if the Knuth-Bendix algorithm does not work, the only methods currently known in the literature to enumerate the relators are two methods based on \textit{van Kampen diagrams}. Let us fix natural numbers $k$ and $n$; the first method consists in enumerating all the van Kampen diagrams of area and length up to $k$ and $n$ respectively; the second uses the fact that a relator of area and length up to $k$ and $n$ is a product of conjugates of defining relators with the length of the conjugating elements bounded above by $mk+n$, where $m$ is the maximal length of a defining relator. 

Unfortunately both algorithms are quite bad and require a huge amount of calculations. 

 The main achievement of this thesis is an algorithm enumerating the relators of area and length up to $k$ and $n$ and whose complexity is less by a factor of at least $a^{k^2} a^{kn}$ than those of the classical methods using van Kampen diagrams ($a$ is a constant).

This is obtained in the following way. First we define by a double recursion a certain set $D$ of relators, whose elements can be enumerated by an algorithm whose complexity is the one shown above. Then we prove that every relator belongs in fact to this set, thus that algorithm enumerates in fact all the relators.

\smallskip

This result can be applied to every situation in which there is an enumeration of relators by means of the associated van Kampen diagrams, since the complexity of this method is better than the usual methods. In fact we show also how to associate with any element of $D$ a van Kampen diagram with certain additional properties, thus proving that it is not necessary to consider all the van Kampen diagrams but only those verifying a special property (see Section \ref{comple}).



 Another possible direction of research can be that of implementing the algorithm presented in this thesis in the form of a computer program for treating group presentations for which the Knuth-Bendix algorithm does not work.
 
It would be also interesting to generalize these results to monoid and semigroup presentations and finding an algorithm for these kinds of presentations.

\section*{Overview of the results} \label{overes}

\addcontentsline{toc}{section}{Overview of the results}

Let $\langle \, X \, | \, R \, \rangle$ be a group presentation where the set of defining relators $R$ is a set of cyclically reduced non-empty words such that $R^{-1}\subset R$. 
 
Let us denote $C_1$ the set of the cyclic conjugates (Definition \ref{cycon}) of the defining relators. For $k>1$ let us denote $C_k$ the set of the cyclic conjugates of the cyclically reduced products (Definition \ref{pi}) of words $u$ and $v$ such that the following three conditions are satisfied: \begin{enumerate}

  \item $u\in C_{k-1}$ and $v\in C_1$ or $u\in C_1$ and $v\in C_{k-1}$;
  
  \item $u\neq v^{-1}$;
  
  \item either $u^{-1}$ is a prefix of $v$, or $v^{-1}$ is a suffix of $u$ or the last letter of $u$ is the inverse of the first of $v$ but the first  letter of $u$ is not the inverse of the last of $v$.
  \end{enumerate}

  We call $k$-corollas the elements of $C_k$ (Definitions \ref{overlineR} and \ref{kg. corolla}, Remark \ref{Roverlinen}). If $R$ is finite then there is a finite algorithm which computes $C_k$ (Remark \ref{Roverlinen}). We denote $C$ the union of all the $C_k$.

  Let us consider for example the following presentation
   $$\langle \, a, b \, | \, aba^{-1} b^{-1}\, \rangle$$
which is the standard presentation for $\mathbb{Z}^2$, the free abelian group of rank 2. In fact, since we have assumed that the set of defining relators must contain the inverse of any of its elements, then we rewrite the above presentation in the following form
   $$\langle \, a, b \, | \, aba^{-1} b^{-1}, \, bab^{-1}a^{-1} \, \rangle,$$
 where the second defining relator is the inverse of the first. Let us compute the set of $k$-corollas for $k=1, 2$. Let us denote $\mathcal{A}$ the set of generators and inverses of them, that is 
      $$\mathcal{A}:=\{a, b, a^{-1},  b^{-1}\}.$$
$C_1$ is the set of the cyclic conjugates of the defining relators, that is
  \begin{equation} \label{vaf} C_1:=\{\textbf{xyx}^{-1}\textbf{y}^{-1}  \,\, : \,\, \textbf{x}, \textbf{y} \in \mathcal{A}, \,\, \textbf{x}^{\pm1}, \textbf{y}^{\pm1} \,\, \textrm{distinct}\}. \end{equation}
   We express the relators as monomials in the unknowns $\textbf{x}$ and $\textbf{y}$; to obtain the actual relators one has just to replace the unknowns with the letters of $\mathcal{A}$. 
  By the definition, the elements of $C_2$ are the cyclic conjugates of the cyclically reduced products of the pairs $(u, v)$ of elements of $C_1$ satisfying conditions 1-3. Since no element of $C_1$ is a proper prefix or suffix of another element of $C_1$, then the only product to consider is 
               $$(\textbf{xyx}^{-1}\textbf{y}^{-1}) (\textbf{yx}^{-1}\textbf{y}^{-1}\textbf{x})=\textbf{x}\textbf{yx}^{-2}\textbf{y}^{-1}\textbf{x};$$
 the second member of the above product is obtained by replacing in the general monomial of (\ref{vaf}) the letter $\textbf{x}$ with $\textbf{y}$ and $\textbf{y}$ with $\textbf{x}^{-1}$, that is if $m(\textbf{x}, \textbf{y})=\textbf{xyx}^{-1}\textbf{y}^{-1}$ then $\textbf{yx}^{-1}\textbf{y}^{-1}\textbf{x}=m(\textbf{y}, \textbf{x}^{-1})$.
 Now we have to compute the cyclic conjugates of this monomial to obtain $C_2$, that is 
 $$C_2:=\{\textbf{x}\textbf{yx}^{-2}\textbf{y}^{-1}\textbf{x}, \,\, \textbf{yx}^{-2}\textbf{y}^{-1}\textbf{x}^2, \,\, \textbf{x}^{-2}\textbf{y}^{-1}\textbf{x}^2\textbf{y},  \,\, \textbf{x}^{-1}\textbf{y}^{-1}\textbf{x}^2\textbf{y}\textbf{x}^{-1},$$
       $$\textbf{y}^{-1}\textbf{x}^2\textbf{y}\textbf{x}^{-2},\,\, \textbf{x}^2\textbf{y}\textbf{x}^{-2} \textbf{y}^{-1}, \,\,  \,\, : \,\, \textbf{x},\textbf{y} \in \mathcal{A}, \,\, \textbf{x}^{\pm1}, \textbf{y}^{\pm1} \textrm{distinct}\}.$$
   In fact the last three monomials are obtained from the first three by replacing each letter with its inverse, thus we can write $C_2$ in the following more compact form    
       $$C_2:=\{\textbf{x}^{-1}\textbf{y}^{-1}\textbf{x}^2\textbf{y}\textbf{x}^{-1},\,\, \textbf{y}^{-1}\textbf{x}^2\textbf{y}\textbf{x}^{-2},\,\, \textbf{x}^2\textbf{y}\textbf{x}^{-2} \textbf{y}^{-1}, \,\,  \,\, \textbf{x},\textbf{y} \in \mathcal{A}, \,\, \textbf{x}^{\pm1}, \textbf{y}^{\pm1} \textrm{distinct}\}.$$
       As seen in Remark \ref{bau}, with any $k$-corolla $c$ can be associated in a standard way a product of conjugates of defining relators whose reduced form is $c$. Since as we have seen, any element of $C_2$ can be obtained in two ways, then we have that all the 2-corollas of this presentation can be expressed in two different ways as product of conjugates of defining relators. 
       
       \smallskip

Let us return to the general case. Given words $w:=w_1 w_2$ and $u$, we say that the word $w_1 u w_2$ is \textit{the insertion of $u$ into $w$ at $w_1$} (Definition \ref{insert}).

We denote $D_1$ the set of words which are reduced and which are of the form $u w u^{-1}$ where $w \in C_1$. Suppose by induction to have defined $D_h$ for every $h<k$ and define $D_k$ as the set of reduced words which are either of the form $u w u^{-1}$ where $w \in C_k$ or are insertions of a word of $D_m$ into one of $D_n$ for $m+n=k$. We denote $D$ the union of all the $D_k$.

In a way analogous to that of Remarks \ref{bau} and \ref{miao} for $C_k$, we can associate with any element $w$ of $D_k$ a product of $k$ conjugates of defining relators whose reduced form is $w$ and a van Kampen diagram labeled by $w$ with $k$ faces. In particular that any element of $D$ is a relator.

The main result of the thesis is the following (Theorem \ref{maithe} and Section \ref{simpli}): \textbf{Let $w$ be a relator of $\langle \, X \, | \, R \, \rangle$; then $w$ belongs to $D$. Moreover $\mathrm{Area}(w)=\min\{k: w\in D_k\}$.}

To compute the relators of area and length up to $k$ and $n$ one can use the algorithm consisting in enumerating the elements of $D_1 \cup D_2 \cup \cdots \cup D_k$ of length up to $n$. This algorithm is presented in Section \ref{fuzumi}. There are currently two other methods in the literature for this kind of computation and our algorithm is much better than both of them.

 These two methods, which are illustrated in Section \ref{enume}, are the following. The first consists in enumerating all the van Kampen diagrams of area up to $k$ and boundary of up to $n$ edges. The second classical method uses the following fact (Remark \ref{smer}), which is a consequence of the properties of van Kampen diagrams: a relator of area and length up to $k$ and $n$ is a product of conjugates of defining relators with the length of the conjugating elements bounded above by $mk+n$, where $m$ is the maximal length of a defining relator. This implies that one can find all the relators of area and length up to $k$ and $n$ by computing a finite number of products of conjugates of defining relators. The complexity of this method is $q^k a^{kn} s^{k^2}$, where $q, a$ and $s$ are constants depending on the presentation (see Section \ref{enume} for details). We have the following result (Theorem \ref{furish}): \textbf{The complexity of the algorithm presented in this thesis is bounded above by $n^3 \, \beta^n \, \gamma^k$, where $\beta$ and $\gamma$ are constants, and is less than those of the classical methods described above by a factor of at least $a^{k(k-c_1)} a^{n(k-c_2)}$, with $c_1$ and $c_2$ constants.}

Let us make some observations about these complexities. The value $n^3 \, \beta^n \, \gamma^k$ given above is  far from being an optimal bound for the complexity of our method; it is only a value obtained with not too complicated calculations. With non-trivial arguments this result could be probably much improved. Instead the value $q^k a^{kn} s^{k^2}$ is the real complexity of the classical method, that is it is the exact number of words to be considered with this method to find all the relators of length and area up to $n$ up to $k$. The same is true also for the factor $2^k$ given as the improvement with respect to the first method and which could be probably much bigger.

\smallskip

Let us present other results; we say that a $k$-corolla $c$ is a \textit{proper $k$-corolla} if for $h<k$, $c$ does not belong to $D_h$. In particular if $c$ is a proper $k$-corolla then $\textrm{Area}(c)=k$. For every natural number $n$ we set 
$$\Delta_0(n):=\max\{\textrm{Area}(c) : \textrm{$c$ is a proper corolla} \,\,\textrm{and} \, \, |c|\leqslant n\}.$$
    Obviously $\Delta_0(n)\leqslant \Delta(n)$ where $\Delta(n)$ is the Dehn function (Definition \ref{Dehn}). We have (Theorem \ref{aska}): \textbf{Let $\langle \, X \, | \, R \, \rangle$ be a finite presentation; then the Word Problem is solvable if and only if $\Delta_0$ is bounded above by a computable function.}

This result improves a well known result (see Proposition \ref{gersho}) because $\Delta_0$ is bounded above by the Dehn function. An analogous result can be stated for non-necessarily finite presentations. In fact R. I. Grigorchuk and S. V. Ivanov in a recent paper \cite{GI} have introduced a function $f_1$ which plays for finitely generated decidable presentations the same rôle of the Dehn function. They proved that for such presentations, the Word Problem is solvable if and only if $f_1$ is bounded above by a computable function (see Proposition \ref{risGI} of this thesis).

As pointed out by Grigorchuk and Ivanov at the end of Sec. 2 of their paper \cite{GI}, the function $f_1$ is equivalent (in the usual sense, see Definition \ref{equiv}) to the function \textit{Work} introduced by J.-C. Birget in (\cite{Birg}, 1.4). We have proved in Section \ref{tcp} this equivalence; moreover we have introduced in Definition \ref{propcor} a function $\Omega_0$ such that if $\Omega$ is the function Work then $\Omega_0 \leqslant \Omega$ and the Word Problem is solvable for a finitely generated decidable presentation if and only if $\Omega_0$ is bounded by a computable function. This improves the result of Grigorchuk and Ivanov.

The main problem with the function $\Omega$ is that, unlike the Dehn function which is invariant (up to equivalence) for finite presentations of the same group, there is not a similar situation for finitely generated decidable presentations. Anyway we have found conditions under which the function $\Omega$ is invariant. Let us give some details.

What is special about finite presentations for the same group is that given any two of them, one can be obtained from the other by applications of \textit{elementary} Tietze transformations. Given two non-necessarily finite presentations for the same group this is not always the case. In fact we have generalized this situation by introducing the so-called \textit{Tietze transformations of bounded area} (see Section \ref{dpfg}) and we have proved that the Dehn function does not change if one applies to a (finite or infinite) presentation some Tietze transformations of bounded area (an elementary Tietze transformation is of bounded area).

In the same way we have introduced the so-called \textit{Tietze transformations of bounded work} (see Definition \ref{Ttobw}) and we have proved that $\Omega$ does not change if one applies to a presentation Tietze transformations of this kind. This improves a result of Grigorchuk and Ivanov (\cite{GI}, The. 1.6) which states that $\Omega$ is invariant under some transformations called \textit{T-transformations} and \textit{stabilizations} (they are special cases of Tietze transformations of bounded work).

We have also introduced a function, denoted $\Gamma$, which plays for the Conjugacy Problem the same rôle of the Dehn function for the Word Problem. We have proved (Proposition \ref{CP}) that the Conjugacy Problem is solvable for a finitely generated presentation if and only if the Word Problem is solvable and $\Gamma$ is bounded above by a computable function. It can be surprising that despite the fact that the Conjugacy Problem is harder to solve than the Word Problem, here the finite generation is sufficient to obtain the same results for the Word Problem. No hypothesis is needed for the set of defining relators, which can be infinite and also non-decidable (even non-enumerable). It is surprising too that the function $\Gamma$ is invariant for any two finitely generated presentations of the same group (Proposition \ref{new}).

Other results are the following: 

     Corollary \ref{cont}: \textbf{If $w$ is a relator of $\langle \, X \, | \, R \, \rangle$ then $w$ has a contiguous subword equal to a proper corolla whose area is at most equal to the area of $w$.}
 
       Corollary \ref{hyper2}: \textbf{Let $\alpha$ be a positive constant; if $\Delta_0(n)\leqslant \alpha n$ then $\Delta(n)\leqslant \alpha n$} (and the presentation is hyperbolic, Definition \ref{hyper}). This result is generalized in Propositions \ref{vega} and \ref{ange}: \textbf{If $\Delta_0$ is equivalent to a function $f$ which is polynomial or exponential, then $\Delta$ too is equivalent to $f$} (see Definition \ref{equiv} for equivalence of functions).

\section*{Structure of the thesis}

\addcontentsline{toc}{section}{Structure of the thesis}

The thesis is organized as follows: Chapter \ref{faipog} is about group presentations and decision problems. In this chapter we do not simply list all the known results, instead we generalize some of them and we present some new notions. In Chapter \ref{tsral} we introduce the main tool of the thesis, \textit{straight line algorithms}, and the main objects, the sets $C$, $L$ and the function $A$. In Chapter \ref{pr} we present results necessary for the proof of the Main Theorem, which occupies Chapter \ref{tmt}. In Chapter \ref{aac} we present some applications of the results proved and we compute the complexity of the algorithm presented in this thesis, which is shown to be better than the other methods currently known in the literature. Finally in the Appendix we prove some technical results.

\smallskip \smallskip

\chapter{Finite and infinite presentations of groups} \label{faipog}

This chapter presents result about group presentations and decision problems for groups. Sections \ref{intro} and \ref{dpfg} introduce the elementary facts about words, van Kampen diagrams, the area of a relator and decision problems for group presentations. Section \ref{tt} is about Tietze transformations and the invariance of the Dehn function under a change of group presentation. To address the last question we introduce the original notion of \textit{Tietze transformation of bounded area}. Sections \ref{ipog} and \ref{tcp} treat respectively the Word Problem for infinite presentations and the Conjugacy Problem. The approach is the same, we define functions analogous to the Dehn function and consider the classes of Tietze Transformations which preserve them. These two sections contain some original ideas and results. Finally Section \ref{enume} illustrates the two methods known in the literature for computing the relators of a group presentation. 

Generalizations of the results of this section to the case of monoid presentations can be found in \cite{VaccDecisProb}.

\section{Words and presentations} \label{intro}

Let $X$ be a set (finite or infinite), let $X^{-1}$ be a set disjoint from $X$ such that $|X| = |X^{-1}|$ and suppose given a bijection $X\rightarrow X^{-1}$. We denote $x^{-1}$ the image by this bijection of an element $x \in X$ and we call it \textit{the inverse of x}. If $y \in X^{-1}$ we denote $y^{-1}$ the element of $X$ such that $(y^{-1})^{-1}=y$. We call \textit{letters} the elements of $X\cup X^{-1}$.

Let $\mathcal{M}(X \cup X^{-1})$ be the free monoid on $X\cup X^{-1}$ and $\mathcal{F}(X)$ the free group on $X$. The elements of $\mathcal{M}(X \cup X^{-1})$ are called \textit{words}; the unity of $\mathcal{M}(X \cup X^{-1})$ is the word with zero letters, called the \textit{the empty word} and denoted 1. The elements of $\mathcal{F}(X)$ are called \textit{reduced words}. 

Let $w:=x_1 \cdots x_m$ be a word; the $x_i$ are called \textit{the letters of $w$}. Given $1 \leqslant i_1 < \cdots < i_n \leqslant m$, the word $\prod_{\alpha=1}^m x_{i_{\alpha}}$ is called a \textit{(non-necessarily contiguous) subword} of $w$. A \textit{contiguous subword} of $w$ is a subword in which the indices $i_1, \cdots, i_n$ are consecutive.

A \textit{null-word} is a word of the form $xx^{-1}$ where $x$ is a letter. A word is \textit{(freely) reduced} if it has no contiguous subword equal to a null word. Given a word $w$ of the form $w' xx^{-1} w''$, the word $w' w''$ is a \textit{one step (free) reduction} of $w$. A sequence of words $w=w_0, w_1, \cdots, w_n$ is a called a \textit{(free) reduction of $w$} if $w_i$ is a one step reduction of $w_{i-1}$ and $w_n$ is reduced. The word $w_n$ is called a \textit{reduced form of $w$}. By The. 1.2 of \cite{MKS}, two reduced forms of the same word are equal, then we can talk about \textit{the} reduced form of $w$.

\begin{definition} \rm \label{rho} We let $\rho: \mathcal{M}(X \cup X^{-1}) \rightarrow \mathcal{F}(X)$ be the function sending a word to its unique \textit{reduced form}.  \end{definition}

\begin{definition} \rm \label{presuf} Let $w_1$ and $w_2$ be words and let $w=w_1 w_2$. The words $w_1$ and $w_2$ are called respectively \textit{a prefix} and \textit{a suffix} of $w$. \end{definition}

\begin{definition} \label{reducprod}  \rm Given $u, v \in \mathcal{F}(X)$ there exist $u_1, v_1,a \in \mathcal{F}(X)$ such that $u=u_1 a$, $v=a^{-1} v_1$ and $\rho(uv)=u_1 v_1$. $\rho(uv)$ is called the \textit{reduced product} of $u$ by $v$ and is the product in $\mathcal{F}(X)$, whereas $uv$ denotes the product in $\mathcal{M}(X \cup X^{-1})$, which is the juxtaposition of words. Therefore $uv=u_1 a a^{-1} v_1$. The word $a a^{-1}$ is called \textit{the cancelled part in the reduced product of $u$ by $v$}. \end{definition}

Let $w:=x_1 \cdots x_n$ be a word; the word $x_n^{-1} \cdots  x_1^{-1}$ is \textit{the inverse of $w$} and is denoted $w^{-1}$. The \textit{length of w} is $|w| = n$. It is easy to see that $\rho(ww')=\rho\big(\rho(w) \rho(w')\big)$ and that $\rho(w^{-1})=\rho(w)^{-1}$.

\begin{definition} \label{cycon} \rm Let $w_1$ and $w_2$ be words, let $w:=w_1 w_2$ and let $n=|w_1|$. The word $w_2 w_1$ is called \textit{the n-th cyclic conjugate of $w$}. If $n > |w|$ we define the $n$-th cyclic conjugate of $w$ as $w$ itself.
     \end {definition}

\begin{definition} \label{cyrewo}  \rm A reduced word is \textit{cyclically reduced} if its last letter is not the inverse of the first one, that is if all its cyclic conjugates are reduced. We denote $\mathcal{F}(X)_c$ the set of cyclically reduced words.   \end {definition}

\begin{definition} \label{scope} \rm Given a word $u$, either $\rho(u)$ is cyclically reduced or there exist (a unique) $t \in \mathcal{F}(X)\setminus{\{1\}}$ and $w \in \mathcal{F}(X)_c$ such that $\rho(u)=twt^{-1}$. The word $w$ is called the \textit{cyclically reduced form of $u$.}  \end {definition}

If a word is cyclically reduced then it coincides with its own cyclically reduced form.

\begin{definition} \label{pi}  \rm We let $\pi: \mathcal{M}(X\cup X^{-1}) \times \mathcal{M}(X\cup X^{-1}) \rightarrow \mathcal{F}(X)_c$ be the function sending two words to the cyclically reduced form of their product. Given two words $u$ and $v$, the word $\pi(u, v)$ is called the \textit{cyclically reduced product} of $u$ by $v$.  \end {definition}

The restriction of $\pi$ to $\mathcal{F}(X)_c \times \mathcal{F}(X)_c$ is a non-associative product in $\mathcal{F}(X)_c$.

Let $\langle \, X \, | \, R \, \rangle$ be a group presentation where $X$ is the set of generators and $R$ that of defining relators. It is not restrictive to suppose that $R$ is a set of cyclically reduced non-empty words such that $R^{-1}\subset R$.

\medskip

\textbf{Van Kampen diagrams.} Let $w$ be a non-necessarily reduced relator, that is the reduced form of $w$ is a relator. This means that there exist defining relators $r_1, \cdots, r_k$ and words $a_1, \cdots, a_k$ such that the reduced form of $w$ is equal to the reduced form of $a_1 \, r_1 \, a_1^{-1} \cdots a_k \, r_k \, a_k^{-1}$. Then there exists a simply connected and planar 2-cell complex with the following properties: \begin{enumerate}

  \item the edges are labeled by elements of $X\cup X^{-1}$,
  
  \item there are at most $k$ faces and the set\footnote{$\{r_1, \cdots, r_k\}$ is in fact a multiset since some of the $r_i$ can be repeated} of the labels of their boundaries is a subset of $\{r_1, \cdots, r_k\}$,
  
  \item there is a cycle, called a \textit{boundary cycle}, labeled by $w$ and containing all the exterior edges.
  
\end{enumerate}

This complex is called \textit{a van Kampen diagram for $w$} and the word $w$ is called the \textit{boundary label} for the diagram. The way to construct it is the following (see V.1 of \cite{LS}). One starts with $k$ faces, $k$ simple cycles (which we denote $\gamma_1, \cdots, \gamma_k$) bounding each one a face and $k$ simple non-closed paths (which we denote $\pi_1, \cdots, \pi_k$). The initial vertices of all the $\pi_i$ coincide and for every $i$, the final vertex of $\pi_i$ coincides with the initial one of $\gamma_i$. All the other vertices are simple. $\pi_i$ is labeled by $a_i$ and $\gamma_i$ by $r_i$.

Let us consider the cycle which traverses once the $\gamma_i$ and which traverses twice the $\pi_i$ in the two opposite directions (one giving $a_i$ and the other $a_i^{-1}$). This cycle contains all the exterior edges of the complex (in fact there are no interior edges) and its label is $a_1 \, r_1 \, a_1^{-1} \cdots a_k \, r_k \, a_k^{-1}$. The word $w$ is obtained from $a_1 \, r_1 \, a_1^{-1} \cdots a_k \, r_k \, a_k^{-1}$ by canceling or inserting words of the form $xx^{-1}$. The idea is that the cancelation of $xx^{-1}$ “translates” into the folding of the two edges labeled by $x$ and $x^{-1}$, giving an interior edge; and the insertion of $xx^{-1}$ translates into an insertion of an edge labeled by $x$ (its opposite is labeled by $x^{-1}$). 

When all these “foldings” and insertions are carried out, the resulting complex is labeled by $w$. In fact the situation is a bit more complicated because in this process some 2-sphere can be discarded (see 2.1 of \cite{Carb}) and the final complex can have less than $k$ faces. 

Therefore a van Kampen diagram is determined by a product of conjugates of defining relators and a sequence of cancellations and insertions of null words.

\medskip

\textbf{Identities among the relations.} Let $r_1, \cdots, r_m$ be defining relators, let $a_1, \cdots, a_m$ be words and suppose that the reduced form of 
         $$a_1 r_1 a_1^{-1} \cdots a_m r_m a_m^{-1}$$
is equal to 1. Then we say that the expression $E:=\big((a_1, r_1), \cdots, (a_m, r_m)\big)$ \textit{determines an identity among the relations}. Identities among the relations exist in any presentation; for instance if $r$ and $s$ are defining relators then the reduced forms of $(r)(r^{-1})$ and $(rsr^{-1})(r)(s^{-1})(r^{-1})$ are equal to 1. 

The group of the “non-trivial” of such identities is the quotient of the group of all the identities among the relations over the subgroup generated by the so-called \textit{Peiffer relations} (see \cite{BroHue} or III.10 of \cite{LS}). This group is isomorphic to the second homotopy group of the presentation complex (see \cite{BroHue}; for the presentation complex see III.2 of \cite{LS}). We recall that the first homotopy group of the complex of a presentation of a group $\mathcal{G}$ is isomorphic to $\mathcal{G}$ itself (Prop. III.2.3 of \cite{LS}). 

Let $w$ be a relator and suppose that a single relator can be expressed in two ways as product of conjugates of defining relators, that is 
    $$w=\rho(a_1 r_1 a_1^{-1} \cdots a_m r_m a_m^{-1})=\rho(b_1 s_1 b_1^{-1} \cdots b_n s_n b_n^{-1}).$$
    It is obvious that this gives an identity among relations, namely
    \begin{equation} \label{smart} \rho(a_1 r_1 a_1^{-1} \cdots a_m r_m a_m^{-1}\, b_1 s_1^{-1} b_1^{-1} \cdots b_n s_n^{-1} b_n^{-1})=1. \end{equation} 
    Let us consider the van Kampen diagrams associated with some free reductions of $a_1 r_1 a_1^{-1} \cdots a_m r_m a_m^{-1}$ and $b_1 s_1 b_1^{-1} \cdots b_n s_n b_n^{-1}$. These two diagrams need not to be isomorphic. If they are, then the identity (\ref{smart}) is trivial (see \cite{CH}).

\medskip

\textbf{Area of a relator.} For a group presentation $\langle \, X \, | \, R \, \rangle$ we denote $\mathcal{N}$ the set of relators, that is the normal closure of $R$ in $\mathcal{F}(X)$.

\begin{definition} \label{defarea} \rm For every $w \in \mathcal{N}$ there exists a natural number $k$ and there exist $r_1, \cdots, r_k \in R$ and $a_1, \cdots, a_k \in \mathcal{F}(X)$ such that $w$ is the reduced form of $a_1 \, r_1 \, a_1^{-1} \cdots a_k \, r_k \, a_k^{-1}$. We call \textit{area of $w$} \big(denoted $\textrm{Area}(w)$\big) the least of such $k$, that is the least $k$ such that $w$ can be expressed in $\mathcal{F}(X)$ as product of $k$ conjugates of defining relators.
    \end{definition}

The area of a relator was introduced by K. Madlener and F. Otto in (\cite{MO}, Sec. 3) with the name of \textit{derivational complexity} (in the more general context of monoid presentations) and later by M. Gromov independently  in (\cite{Grom}, Sec. 2.3).

\begin{definition} \label{} \rm Let $w$ and $a$ be reduced words and let $r$ be a defining relator. Set $w':=\rho(w a r a^{-1})$; we say that the pair $(w, w')$ is a \textit{derivation of length one (by means of $r$) from $w$ to $w'$}. Since $R$ contains the inverse of any of its elements and since $w=\rho(w' a r^{-1} a^{-1})$ then also $(w', w)$ is a derivation of length one. A finite sequence $(w_1, \cdots, w_{k+1})$ of $k+1$ reduced words is a \textit{derivation of length $k$ from $w_1$ to $w_k$} if $(w_i, w_{i+1})$ is a derivation of length one for $i=1, \cdots, k-1$. More specifically, if $r_1, \cdots, r_k$ are defining relators such that $(w_i, w_{i+1})$ is a derivation by means of $r_i$ and if $S$ is a subset of $R$ containing $r_1, \cdots, r_k$, then we say that $(w_1, \cdots, w_{k+1})$ is an \textit{$S$-derivation of length $k$ (from $w_1$ to $w_k$)}.
       \end{definition}

\begin{remark} \label{cord} \rm Let $w$ and $w'$ be words; then there exists an $S$-derivation of length $k$ from $w$ to $w'$ if and only if there exist words $a_1, \cdots, a_k$ and defining relators $s_1, \cdots, s_k \in S$ such that $w=\rho(w' a_1 \, s_1 \, a_1^{-1} \cdots a_k \, s_k \, a_k^{-1})$.

In particular if $w$ is a relator then the area of $w$ can be defined as the minimal length of a derivation from $w$ to 1. We observe that the area of $w$ can be also defined as the minimal number of faces of van Kampen diagrams with boundary label equal to $w$.
 \end{remark}

We now prove some properties of the Area. For every $w \in \mathcal{N}$ and for every $t \in \mathcal{F}(X)$ we have    
       $$\textrm{Area}\big(\rho(t w t^{-1})\big) \leqslant \textrm{Area}(w)$$ 
    because if 
 $$w=\rho(a_1 \, r_1 \, a_1^{-1} \cdots a_k \, r_k \, a_k^{-1})$$
  with $k=\textrm{Area}(w)$, then 
 $$\rho(t w t^{-1})=\rho(b_1 \, r_1 \, b_1^{-1} \cdots b_k \, r_k \, b_k^{-1})$$
with $b_i =t a_i$. But conversely 
   $$\textrm{Area}(w)=\textrm{Area}\big(\rho(t^{-1}(t w t^{-1}) t)\big) \leqslant \textrm{Area}\big(\rho(t w t^{-1})\big),$$
     that is 

\begin{equation} \label{uwu-1} \textrm{Area}\big(\rho(t w t^{-1})\big)=\textrm{Area}(w). \end{equation}

 This implies that if $w \in \mathcal{N}$ and if $w'$ is a cyclic conjugate of $w$ (Definition \ref{cycon}), then $\textrm{Area}(w')=\textrm{Area}(w)$.

In the same way one proves that $\textrm{Area}(w^{-1})\leqslant \textrm{Area}(w)$ and $\textrm{Area}(w)\leqslant \textrm{Area}(w^{-1})$, therefore 
   \begin{equation} \label{Area w-1} \textrm{Area}(w^{-1})= \textrm{Area}(w). \end{equation}

Finally for every $v, w \in \mathcal{N}$ we have 
   \begin{equation} \label{Area(vw)} \textrm{Area}\big(\rho(v w)\big) \leqslant \textrm{Area}(v)+\textrm{Area}(w), \end{equation}
    because if $v$ and $w$ are respectively products of $h$ and $k$ conjugates of defining relators, then $vw$ is product of $k+h$ of them. This implies that    
        \begin{equation} \label{Area pi(v,w)} \textrm{Area}\big(\pi(v,w)\big) \leqslant \textrm{Area}(v)+\textrm{Area}(w) \end{equation}
      because $v w$ is conjugate to $\pi(v,w)$ and thus has the same area of the latter by (\ref{uwu-1}) ($\pi$ is the cyclically reduced product, see Definition \ref{pi}).

\section{Decision problems for groups} \label{dpfg}

For \textit{computable functions} and \textit{decidable sets} we follow \cite{Shen}. Let $S$ be a subset of $\mathbb{N}$ and let $f$ be a function from $S$ to $\mathbb{N}$. The function $f$ is said \textit{computable} if there exists an algorithm which taken as input an $n\in \mathbb{N}$ then \begin{itemize}

  \item it halts and gives $f(n)$ as output if $n\in S$;

  \item it does not halt if $n\notin S$.
 
\end{itemize}

    A computable function is also called \textit{recursive}. The set $S\subset \mathbb{N}$ is said \textit{decidable} if its characteristic function is computable, i.e., if there exists an algorithm which determines whether an arbitrary $n\in \mathbb{N}$ belongs to $S$. Any finite subset of $\mathbb{N}$ is decidable (see 1.2 of \cite{Shen}). The set $S$ is said \textit{enumerable} if there exists an algorithm which enumerates the elements of $S$. This means that there is an algorithm which for every $n$ outputs a certain $s_n$ and that $S=\{s_n : n \in \mathbb{N}\}$. We assume that for some $n$ the output of this algorithm can be empty, that is to be formal we must write $S=\{s_n : n \in \mathbb{N} \,\, \textrm{and} \, \, s_n \, \, \textrm{is not empty}\}$. We also assume that $s_m$ can be equal to $s_n$ for $m\neq n$, but that we are able to determine whether or not $s_m=s_n$.

    A decidable set is enumerable. It is obvious that $S$ is decidable if and only if $S$ and $\mathbb{N} \setminus{S}$ are enumerable. Decidable and enumerable sets are also called \textit{recursive} and \textit{recursively enumerable} respectively.

Let $U$ and $V$ be enumerable sets, that is $U=\{u_n : n \in \mathbb{N}\}$ and $V=\{v_n : n \in \mathbb{N}\}$, where there are two algorithms whose $n$-th inputs are $u_n$ and $v_n$ respectively. Let $V'$ be a subset of $V$ and let $f$ be a function from $V'$ to $U$. Let $\overline{f}$ be the following function: let $n$ be a natural such that $v_n \in V'$ and let $f(v_n)=u_m$; then set $\overline{f}(n)=m$. Let $E=\{n \in \mathbb{N} : v_n \in V'\}$; then $\overline{f}$ is a function from $E$ to $\mathbb{N}$. We say that $f$ is \textit{computable} if $\overline{f}$ is computable, that is if there exists an algorithm which taken $v_n$ as input, it halts and gives $f(v_n)$ as output if $v_n \in V'$, otherwise it does not halt. In an analogous way one can define \textit{decidable} and \textit{enumerable} subsets of an enumerable set.

Let $\mathcal{P}=\langle \, X \, | \, R \, \rangle$ be a group presentation. The normal subgroup of $\mathcal{F}(X)$ normally generated by $R$ is the set of relators of $\mathcal{P}$ which we denote $\mathcal{N}$.

If $X$ is enumerable then also $\mathcal{F}(X)$ is enumerable. Let $R$ be enumerable too; then the set of products of the form
               $$a_1 r_1 a_1^{-1} \cdots  a_k  r_k a_k^{-1}$$
is enumerable for $a_i \in \mathcal{F}(X)$, $r_i \in R$ and $k \in \mathbb{N}$, that is the set of relators is enumerable.

A group presentation $\mathcal{P}=\langle \, X \, | \, R \, \rangle$ is said \textit{finitely generated} if $X$ is finite; it is said \textit{finitely related} if $R$ is finite. A presentation is \textit{finite} if it is finitely generated and related. A presentation which is infinitely generated but finitely related is a free product of a finite presentation with an infinite rank free group presentation.

A finitely generated presentation $\mathcal{P}$ is said \textit{decidable} or \textit{enumerable} if $R$ is respectively decidable or enumerable in $\mathcal{F}(X)$. 

If a group admits an enumerable presentation then it admits also a decidable one. Indeed let $\mathcal{G}$ be a group and let $\langle \, X \, | \, R \, \rangle$ be an enumerable presentation of $\mathcal{G}$. This means that there is an algorithm which enumerates $R$, that is there is an algorithm whose $k$-th output is a word $w_k$ and $R=\{w_k : k \in \mathbb{N}\}$. Let $y$ be an element not belonging to $\mathcal{F}(X)$ and let $R'=\{ y^k w_k : k \in \mathbb{N}\}$. Then $\langle \, X\cup \{y\} \, | \, R'\cup \{y\} \, \rangle$ is a decidable presentation of $\mathcal{G}$.

A presentation \textit{has a solvable Word Problem} if given any $w\in \mathcal{F}(X)$ there exists an algorithm which determines whether $w$ is a relator; this is equivalent to say that the set of relators is decidable. A presentation \textit{has a solvable Word Search Problem} if it has a solvable Word Problem and for any relator $w$ there exists an algorithm which finds a product of conjugates of defining relators whose reduced form is $w$.

 If $X$ and $R$ are enumerable then the solvability of the Word Problem is equivalent to that of the Word Search Problem. Indeed if $w$ is a relator then we can execute the algorithm described above which computes the reduced form of the products of conjugates of defining relators; necessarily after a finite time this algorithm will give $w$ as output.

The complexity of the Word Search Problem can be much greater than that of the simple Word Problem. Madlener and Otto showed in (\cite{MO}, Cor. 6.14) group presentations with Word Problem solvable in polynomial time and with Word Search Problem\footnote{Madlener and Otto call \textit{pseudo-natural algorithm} an algorithm for solving the Word Search Problem} arbitrarily hard. Namely for every $m\geqslant 3$ they constructed a group $G(m)$ with these properties: the complexity of the Word Problem for $G(m)$ is at most polynomial, that of the Word Search Problem is bounded above by a function in the \textit{Grzegorczyk class} $E_m$ (see Ch. 12 of \cite{EC}) but by no function in the Grzegorczyk class $E_{m-1}$.

\begin{remark} \label{smer} \rm Let $\mathcal{P}:=\langle \, X \, | \, R \, \rangle$ be a finite presentation, let $m$ be the maximal length of elements of $R$, let $w$ be a relator of $\mathcal{P}$ and let $r_1, \cdots, r_k$ be defining relators and $a_1, \cdots, a_k$ be words such that $w$ is the reduced form of 
          $$a_1 r_1  a_1^{-1} \cdots a_k r_k u_k^{-1}.$$
Let $|w|=n$; by using the properties of the associated van Kampen diagram one can prove that there exist words $b_1, \cdots, b_k$ such that $|b_i|\leqslant mk+n$ and $w$ is the reduced form of 
              $$b_1 r_1  b_1^{-1} \cdots b_k r_k b_k^{-1},$$
that is we can choose the conjugating elements of length no more than $mk+n$. This implies that the sum of the lengths of the conjugating elements is bounded above by $mk^2+kn$. See Prop. 2.2 of \cite{Gerst} or The. 1.1 and 2.2 of \cite{Short} for a proof of this fact\footnote{In (\cite{LS}, Rem. after Lem. V.1.2) it is given $kn$ as bound which is worse than $mk+n$ since $m$ is a constant.} (see also Lem. 7.1 of \cite{Miller}).
    \end{remark}

\begin{definition} \label{Dehn} \rm Let $\mathcal{P}$ be a finitely generated presentation and let $n$ be a natural number; the \textit{Dehn function} at $n$ of $\mathcal{P}$ is
    $$\Delta(n):=\max\{\textrm{Area}(w) : w \in \mathcal{N} \, \textrm{and} \, |w|\leqslant n\}.$$
    \end{definition}
    
The Dehn function is a function from $\mathbb{N}$ to $\mathbb{N}$. We have the following

\begin{proposition} \label{istar} Let $\mathcal{P}$ be a presentation whose set of defining relators is enumerable and with solvable Word Problem and let $\Delta$ be the Dehn function of $\mathcal{P}$. If $\mathcal{P}$ is finitely generated then $\Delta$ is bounded above by a computable function; if $\mathcal{P}$ is finite then $\Delta$ is computable.
\end{proposition}

\begin{proof} Since the set of generators and that of the defining relators are enumerable, the solution of the Word Problem implies that of the Word Search Problem. Let $n$ be a natural number; since the presentation is finitely generated then the number of words of length bounded by $n$ is finite. For any of these words we solve the Word Search Problem; this gives for any relator $w$ of length less or equal to $n$ an upper bound $h_w$ for the area of $w$. We can compute in a finite time the maximum of these bounds and this is an upper bound for $\Delta(n)$.

Let $\mathcal{P}$ be finite and let $m$ be the maximal length of defining relators of $\mathcal{P}$. Let $n$ be a natural number and let $w$ be a relator such that $|w|=n'\leqslant n$. For any $k=1, \cdots, h_w$ let us compute the reduced forms of the products of conjugates of $k$ defining relators with the length of conjugating elements bounded above by $mk+n'$. By virtue of Remark \ref{smer}, at least one these reduced forms for some $k$ is equal to $w$ and a minimal $k$ for which this happens is equal to the area of $w$. Therefore we can compute $\Delta(n)$ in a finite time and thus $\Delta$ is computable.  \end{proof}

When $\mathcal{P}$ is finite the converse of Proposition \ref{istar} is true:

\begin{proposition} \label{gersho} A finite presentation has a solvable Word Problem if and only if its Dehn function is computable, if and only if its Dehn function is bounded above by a computable function.  \end{proposition}

\begin{proof} Let $\Delta$ be the Dehn function of $\mathcal{P}$. By Proposition \ref{istar}, if the Word Problem is solvable for $\mathcal{P}$ then $\Delta$ is computable which implies trivially that $\Delta$ is bounded above by a computable function. We have to prove that if $\Delta$ is bounded above by a computable function then the Word Problem is solvable.

Let $\{h_n\}_{n\in \mathbb{N}}$ be a computable function such that $\Delta(n)\leqslant h_n$ and let $w$ be a word of length less or equal to $n$. We have that $w$ is a relator for $\mathcal{P}$ if and only if it is the reduced form of a product of at most $h_n$ conjugates of defining relators with the length of the conjugating elements bounded above by $mh_n+n$ by Remark \ref{smer}. Since these products can be computed in a finite time then we have a solution of the Word Problem.   \end{proof}

Proposition \ref{gersho} was first proved (in a less general form) by Madlener and Otto in (\cite{MO}, Lem. 3.2).

Proposition \ref{gersho} is no longer true if the number of relators is infinite even when it is decidable: see Example 2.4 of \cite{GI}, where it is shown a decidable presentation with unsolvable Word Problem and Dehn function constantly equal to 2.

\medskip

Let $\mathcal{P}=\langle \, X \, | \, R \, \rangle$ be a presentation for a group $\mathcal{G}$; $\mathcal{P}$ has \textit{has a solvable Conjugacy Problem} if given $u, v\in \mathcal{F}(X)$ there exists an algorithm which determines whether the elements of $\mathcal{G}$ represented by $u$ and $v$ are conjugated, that is whether there exists a word $t\in \mathcal{F}(X)$ such that $u t v^{-1} t^{-1}$ is a relator. This is equivalent to say that the set
     \begin{equation} \label{sargi} \{(u, v, t) \,: \,\, u=t v t^{-1} \,\, \mathrm{in} \,\, \mathcal{G}\}\end{equation}
is decidable. The solvability of the Conjugacy Problem implies that of the Word Problem if we take $v=1$.

Let $\mathcal{N}$ be the set of relators of $\mathcal{P}$. If $\mathcal{F}(X)$ and $\mathcal{N}$ are enumerable, then the set (\ref{sargi}) is enumerable. This is because in this case the set $\mathcal{F}(X)^3 \times \mathcal{N}$ is enumerable and thus is equal to a set of the form $\{\alpha_n=(u_n, v_n, t_n, w_n) : n\in \mathcal{N}\}$, where there is an algorithm whose $n$-th output is $\alpha_n$. Consider the algorithm whose $n$-th output is $(u_n, v_n, t_n)$ if $u_n t_n v_n^{-1} t_n^{-1}=w_n$, otherwise the output is empty; this algorithm enumerates (\ref{sargi}).

We say that $\mathcal{P}$ has \textit{has a solvable Conjugacy Search Problem} if it has a solvable Conjugacy Problem and if there exists an algorithm which finds the conjugating word $t$. If $\mathcal{F}(X)$ is enumerable (in particular if $X$ is enumerable) then the solvability of the Conjugacy Problem is equivalent to that of the Conjugacy Search Problem. Indeed $\mathcal{F}(X) =\{t_n : n \in \mathbb{N}\}$, where there is an algorithm whose $n$-th output is $t_n$. Let $u$ and $v$ be conjugated in $\mathcal{G}$; for every $n$ let us apply the algorithm solving the Word Problem to the word $s_n=t_n \, u \, t_n^{-1}  \, v^{-1}$. For some $n$ the word $s_n$ will be equal to 1, thus $t_n$ will be the conjugating word.

We observe that if $u$ and $v$ are words and if $u'$ and $v'$ are the cyclically reduced forms of $u$ and $v$ respectively (see Definition \ref{scope}), then $u$ and $v$ are conjugated in $\mathcal{G}$ if and only this is true for $u'$ and $v'$. Therefore in studying the Conjugacy Problem it is not restrictive to consider only cyclically reduced words.

\medskip

Let $\mathcal{P}_1=\langle \, X_1\, | \, R_1 \, \rangle$ and $\mathcal{P}_2=\langle \, X_2 \, | \, R_2 \, \rangle$ be two group presentations, let $\mathcal{N}_i$ be the normal subgroup of $\mathcal{F}(X_i)$ normally generated by $R_i$ for $i=1,2$ and suppose that there exists a computable function $f$ from $\mathcal{F}(X_1)$ to $\mathcal{F}(X_2)$ which is not necessarily an homomorphism and such that for $u \in \mathcal{F}(X_1)$ we have that $u$ belongs to $\mathcal{N}_1$ if and only if  $f(u)$ belongs to $\mathcal{N}_2$. Then the Word Problem is solvable for $\mathcal{P}_1$ if and only if the same is true for $\mathcal{P}_2$. We observe that it is not necessary that the groups presented by $\mathcal{P}_1$ and $\mathcal{P}_2$ are isomorphic.

\section{Tietze transformations} \label{tt}

Given two finite presentations for the same group, it is well known (see for instance Prop. 1.3.3 of \cite{Brids}) that the Dehn functions relative to them are equivalent. This is a consequence of the fact that one of them can be obtained from the other by applications of elementary Tietze transformations. This result is no longer true if the presentations are not finite. In this section we define a kind of Tietze transformations, said \textit{of bounded area}, that when applied to a presentation, finite or infinite, give a presentation with the same Dehn function (up to equivalence).

\medskip

A group presentation is determined by two sets usually denoted $X$ and $R$, where $X$ is a set of letters and represents the generators of the group and $R$ is a set of (cyclically) reduced words in $X\cup X^{-1}$ and represents the defining relators of the group. We denote $\langle \langle R \rangle \rangle$ the normal subgroup of the free group $\mathcal{F}(X)$ which is normally generated by $R$.

\begin{definition} \label{Transf} \rm Let us consider the following transformations on a group presentation which consist in adding or deleting “superfluous” generators and defining relators. They are called \textit{Tietze transformations} (see Sec. 1.5 of \cite{MKS}):

\begin{enumerate}

  \item Let $S\subset \langle \langle R \rangle \rangle$; then replace $R$ with $R \cup S$, that is replace $\langle \, X \, | \, R \, \rangle$ with $\langle \, X \, | \, R \cup S \, \rangle$;
  
  \item Let $S\subset R$ be such that $\langle \langle R\setminus{S} \rangle \rangle=\langle \langle R \rangle \rangle$;  then replace $R$ with $R\setminus{S}$, that is replace $\langle \, X \, | \, R \, \rangle$ with $\langle \, X \, | \, R\setminus{S} \, \rangle$;
   
  \item Let $U \subset \mathcal{F}(X)$ and for every $u\in U$ let $y_u$ be an element not belonging to $\mathcal{F}(X)$ and such that if $u$ and $u'$ are distinct elements of $U$ we have $y_u\neq y_{u'}$. Set $Y:=\{y_u : u \in U\}$ and $T:=\{y_u u^{-1}: u \in U\}$; then replace $X$ with $X\cup Y$ and $R$ with $R \cup T$, that is replace $\langle \, X \, | \, R \, \rangle$ with $\langle \, X \cup Y \, | \, R \cup T \, \rangle$;
  
  \item Let $Y\subset X$ and suppose that for every $y \in Y$ there exists a word $u_y$ not containing neither $z$ nor $z^{-1}$ for every $z \in Y$ and such that $y u_y^{-1} \in R$. Let $\varphi$ be the homomorphism from $\mathcal{F}(X)$ to $\mathcal{F}(X \setminus{Y})$ which sends any $y \in Y$ to $u_y$ and any $x \in X\setminus{Y}$ to itself. Then set $T:=\{y u_y : y \in Y\}$, replace $X$ with $X\setminus{Y}$ and $R$ with $\varphi(R\setminus{T})$, that is replace $\langle \, X \, | \, R \, \rangle$ with $\langle \, X \setminus{Y} \, | \, \varphi(R\setminus{T}) \, \rangle$.\end{enumerate} \end{definition}

The transformations 1 and 2 add or delete respectively superfluous defining relators; transformation 3 and 4 add or delete superfluous generators. Transformations of types 1 and 2 are inverse one of the other; the inverse of a transformations of type 3 is a transformation of type 4. The inverse of a transformations of type 4 is a transformation of type 3, followed by one of type 1 and by one of type 2; this is because first we have to add $Y$ to the generators and $T$ to the defining relators obtaining $\langle \, X \, | \, T \cup \varphi(R\setminus{T}) \, \rangle$, then add $R$ and finally delete $\varphi(R\setminus{T})$ from the defining relators.

It is obvious that applying such transformations to a presentation, we obtain a presentation for the same group. The converse is also true, that is if $\mathcal{P}$ and $\mathcal{P}'$ are presentations for the same group then there exists a finite sequence of Tietze transformations which applied to $\mathcal{P}$ gives $\mathcal{P}'$ (see Th. 1.5 of \cite{MKS}). In fact by the proof of Th. 1.5 of \cite{MKS} we have the following stronger result:

\begin{proposition} \label{Tietze} Two presentations of the same group can be obtained one from the other by applications of at most four Tietze transformations, namely a transformation of type 3 followed by one of type 1, then by one of type 2 and finally by one of type 4 (some of these transformations can be empty). \end{proposition}

A Tietze transformation is said \textit{finite} if it adds or deletes only finitely many generators or defining relators, that is if the sets $S$, $T$ and $Y$ are finite\footnote{In the literature an \textit{elementary Tietze transformation} is defined as a Tietze transformation which adds or deletes exactly one generator or defining relator. An elementary Tietze transformation is obviously finite and any finite Tietze transformation is obtained by repeated applications of elementary Tietze transformations.}. In the same way a Tietze transformation is said \textit{decidable} if the sets $S$, $T$ and $Y$ are decidable.

By Cor. 1.5 of \cite{MKS} we have that two finite [respectively decidable] presentations define the same group if and only if one can be obtained from the other by repeated applications of finite [respectively decidable] Tietze transformations. Moreover Proposition \ref{Tietze} holds also in the special cases of finite or decidable presentations, that is two finite or decidable presentations for the same group can be obtained one from the other by applications of at most four Tietze transformations of the types specified in Proposition \ref{Tietze}.

  We say that a Tietze transformation of type 1 or 2 is of \textit{bounded area} if it adds or deletes a set of relators whose area is bounded, that is if $\sup\{\textrm{Area}(s) : s\in S\}$ is finite. We say that a Tietze transformation of type 3 or 4 is of \textit{bounded length} if it adds or deletes a set of generators whose length is bounded, that is if $\sup\{|u| : u\in U\}$ or $\sup\{|u_y| : y\in Y\}$ (respectively for type 3 or 4) are finite.

  A finite Tietze transformation is obviously of bounded area or length. If the presentation is finitely generated then any Tietze transformation of bounded length is necessarily finite. But if the presentation is finitely related then a Tietze transformation of bounded area can be finite or infinite.

   \begin{remark} \label{stef} \rm Let $\mathcal{P}$ be a presentation and let $\mathcal{P}'$ be the presentation obtained from $\mathcal{P}$ by application of a Tietze transformation of bounded area or length. We want to determine the relationship between the area of relators and the length of group elements between $\mathcal{P}$ and $\mathcal{P}'$.

  Suppose that $\mathcal{P}'$ is obtained from $\mathcal{P}$ by application of a Tietze transformation of type 1 of bounded area and let $c$ be the maximal area of defining relators added. If $\textrm{Area}$ and $\textrm{Area}'$ are the functions area relative to $\mathcal{P}$ and $\mathcal{P}'$ respectively, then we have for a relator $w$ that
              $$\textrm{Area}'(w) \leqslant \textrm{Area}(w) \leqslant c\textrm{Area}'(w).$$
  This implies that if $\Delta$ and $\Delta'$ are the Dehn functions relative to $\mathcal{P}$ and $\mathcal{P}'$ then
      \begin{equation} \label{T1} \Delta' \leqslant \Delta \leqslant c\Delta'.\end{equation}

 Let us now suppose that $\mathcal{P}'$ is obtained from $\mathcal{P}$ by application of a Tietze transformation of type 3 of bounded length. By using the notation of 3 of Definition \ref{Transf}, we have that there exists a natural number $c$ such that $|u|\leqslant c$ for every $u\in U$. Let $\varphi$ be the function from $\mathcal{F}(X\cup Y)$ to $\mathcal{F}(X)$ sending any element of $X$ to itself and any $y_u\in Y$ to $u$. We have that $\varphi$ is the identity on $\mathcal{F}(X)$, in particular $\varphi(r)=r$ for every $r \in R$, and that $\varphi(t)=1$ for every $t \in T$. Moreover for every $w\in \mathcal{F}(X\cup Y)$ we have that $|\varphi(w)|\leqslant c |w|$.

 If we denote $\mathcal{N}$ and $\mathcal{N}'$ the normal subgroups of $\mathcal{F}(X)$ and $\mathcal{F}(X\cup Y)$ normally generated by $R$ and $R \cup T$ respectively, then we have that $\mathcal{N} \subset \mathcal{N}'$ and $\varphi(\mathcal{N}')=\mathcal{N}$. 
 
Let Area and Area$'$ be the functions area relative to $\mathcal{P}$ and $\mathcal{P}'$ respectively. Let $w \in \mathcal{N}'$ and let 
     $$w=\rho(a_1 r_1 a_1^{-1} \cdots a_k r_k a_k^{-1})$$ 
where $r_i \in R\cup T$ and $k=\textrm{Area}'(w)$. Then 
   $$\rho \big(\varphi(w)\big)=\rho(b_1 s_1 b_1^{-1} \cdots b_k s_k b_k^{-1})$$
where $b_i=\varphi(a_i)\in \mathcal{F}(X)$ and $s_i$ is equal to $r_i$ if $r_i \in R$, $s_i$ is equal to 1 if $r_i \in T$. This implies that 
         $$\mathrm{Area}\big(\varphi(w)\big) \leqslant \textrm{Area}'(w).$$
     If $w \in \mathcal{N}$ then since $\varphi(w)=w$ and since obviously $\mathrm{Area}'(w) \leqslant \textrm{Area}(w)$, then $\mathrm{Area}(w) = \textrm{Area}'(w)$. In particular if $\Delta$ and $\Delta'$ are the Dehn functions relative to $\mathcal{P}$ and $\mathcal{P}'$ respectively, we have that 
                      $$\Delta \leqslant \Delta'.$$
Let $w \in \mathcal{N}'$; then there exists an $(R\cup T)$-derivation from $w$ to $\varphi(w)$ of length at most $|w|$. Since there exists an $R$-derivation from $\varphi(w)$ to 1 of length $\mathrm{Area}\big(\varphi(w)\big)$ then this gives an $(R\cup T)$-derivation from $w$ to 1 of length $\mathrm{Area}\big(\varphi(w)\big) + |w|$, that is 
             $$\mathrm{Area}'(w)\leqslant \mathrm{Area}\big(\varphi(w)\big) + |w|$$ and then 
                     $$\Delta'(n)\leqslant \Delta(cn)+ n,$$ 
that is 
          \begin{equation} \label{T3} \Delta(n) \leqslant \Delta'(n) \leqslant \Delta(cn)+n. \end{equation}

  Let us now determine the relationship between the length of group elements between $\mathcal{P}$ and $\mathcal{P}'$. If $\mathcal{G}$ is the group presented and if $g$ is any element of $\mathcal{G}$ then 
          \begin{equation} \label{hagen} |g|' \leqslant |g| \leqslant c|g|'\end{equation}  
where $|\cdot|$ and $|\cdot|'$ denote respectively the length with respect to the set of generators of $\mathcal{P}$ and $\mathcal{P}'$ respectively. \end{remark}

\medskip

 \begin{remark} \label{total} \rm Let $\mathcal{P}:=\langle \, X \, |  \, R \, \rangle$ be a group presentation such that $X$ and $R$ are enumerable and let $\mathcal{P}'$ be the presentation obtained from $\mathcal{P}$ by adding to $R$ all the relators of $\mathcal{P}$, that is if $\mathcal{N}$ is the set of all relators of $\mathcal{P}$ then $\mathcal{P}'=\langle \, X \, |  \, \mathcal{N} \, \rangle$.
 
The set of defining relatiors of $\mathcal{P}'$ is not finite, is enumerable and is decidable if and only if the Word Problem is solvable for $\mathcal{P}$. The Dehn function for $\mathcal{P}'$ is eventually equal to 1. \end{remark}

 \begin{definition} \label{equiv} \rm Let $f, g: \mathbb{R}_+ \rightarrow \mathbb{R}_+$ be two non-decreasing functions. We write $f \preceq g$ if there exists a positive constant $\alpha$ such that 
          $$f(n)\leqslant \alpha g(\alpha n) + \alpha n$$
   for every $n \in \mathbb{N}^*$. We say that $f$ and $g$ are \textit{equivalent} if $f \preceq g$ and $g \preceq f$ and in this case we write $f \simeq g$.
 
Let $f, g: \mathbb{R}^2_+ \rightarrow \mathbb{R}^2_+$ be two functions. We write $f \preceq_{\texttt{s}} g$ if there exists a positive constant $\alpha$ such that 
     $$f(n_1, n_2)\leqslant \alpha g(n_1, n_2)$$
    for every $n_1, n_2\in \mathbb{N}^*$. 
   We say that $f$ and $g$ are \textit{strongly equivalent} if $f \preceq_{\texttt{s}}g$ and $g \preceq_{\texttt{s}} f$ and in this case we write $f \simeq_{\texttt{s}} g$. \end{definition}

 If $f \leqslant g$ then $f \preceq_{\texttt{s}}g$ and $f \preceq g$. We consider $\Delta$ as a function defined on $\mathbb{R}_+$ by assigning the value $\Delta(n)$ to every $x\in ]n, n+1[$.

 \begin{remark} \label{lis} \rm Since $\alpha g(\alpha n) + \alpha n \leqslant \alpha' g(\alpha' n) + \alpha' n$ if $\alpha \leqslant \alpha'$, then $f \preceq g$ if and only if there exists $\alpha_0>0$ such that $f(n)\leqslant \alpha g(\alpha n) + \alpha n$ for every $\alpha\geqslant \alpha_0.$ \end{remark}

     \begin{proposition} \label{dasz} Let $\mathcal{P}$ and $\mathcal{P}'$ be group presentations, let $\Delta$ and $\Delta'$ be the Dehn functions for $\mathcal{P}$ and $\mathcal{P}'$ respectively. Let $\mathcal{P}'$ be obtained from $\mathcal{P}$ by applications of Tietze transformations of bounded area or length (in particular, let $\mathcal{P}$ and $\mathcal{P}'$ be finite presentations for the same group); then $\Delta$ and $\Delta'$ are equivalent. In particular $\Delta'$ is computable (or bounded above by a computable function) if and only if the same is true for $\Delta$.     \end{proposition}

\begin{proof} By proposition \ref{Tietze}, $\mathcal{P}'$ can be obtained from $\mathcal{P}$ by applications of finitely many Tietze transformations. By $(\ref{T1})$ and $(\ref{T3})$, we have that applying a Tietze transformation of type 1 or 3, the equivalence class of the Dehn function does not change. This is also true by applying a Tietze transformation of type 2 because its inverse is a transformation of type 1; and is true also for a Tietze transformation of type 4 since its inverse is a transformation of type 3 followed by one of type 1 and by one of type 2.    \end{proof}

The result of Proposition \ref{dasz} is stated in the literature only in the case where $\mathcal{P}$ and $\mathcal{P}'$ are finite, that is when $\mathcal{P}$ is finite and $\mathcal{P}'$ is obtained from $\mathcal{P}$ by applications of finite Tietze transformations. See for instance Prop. 1.3.3 of \cite{Brids}; the first proof of this result was given (in a less strong form) by Madlener and Otto in (\cite{MO}, Cor. 3.5).

\begin{remark} \label{chuck} \rm If $\mathcal{P}'$ is obtained from $\mathcal{P}$ by applying non-bounded Tietze transformations then $\Delta$ and $\Delta'$ are in general non-equivalent. For instance if $\mathcal{P}'$ and $\mathcal{P}$ are as in Remark \ref{total}, then $\Delta$ and $\Delta'$ are equivalent if and only if $\Delta$ is equivalent to a constant function. In particular if $\mathcal{P}$ is a finite presentation with unsolvable Word Problem then $\Delta$ is greater than any computable function while $\Delta'$ is constant. This shows that if in Proposition \ref{dasz} we remove the hypothesis of boundedness for the Tietze transformations applied to a presentation, then the Dehn functions of two presentations for the same group can be not related at all.

This shows also that a decidable Tietze transformation needs not be of bounded area. Indeed let $\mathcal{P}$ and $\mathcal{P}'$ as in Remark \ref{total}, let $\mathcal{P}$ be decidable, let $\Delta$ be not equivalent to a constant function and suppose that the Word Problem is solvable for $\mathcal{P}$. Then $\mathcal{P}'$ is decidable, thus one can obtain $\mathcal{P}'$ from $\mathcal{P}$ by applications of decidable Tietze transformations; but since $\Delta$ and $\Delta'$ are not equivalent then $\mathcal{P}'$ cannot be obtained from $\mathcal{P}$ by applications of Tietze transformations of bounded area.
  \end{remark}

 We have seen that two presentations define the same group if and only if they can be obtained one from the other by repeated applications of Tietze transformations. Let us consider a graph whose vertices are all the presentations of all the groups and with an edge joining two presentations if one can be obtained from the other by application of a Tietze transformation; any connected component of this graph corresponds to an isomorphism class of groups. Furthermore by Proposition \ref{Tietze}, two vertices are either non-connected or they are connected by a path of length at most four.

 Moreover let $\mathcal{G}$ be a group, let $\mathcal{P}$ and $\mathcal{P}'$ be two finite [respectively decidable] presentations of $\mathcal{G}$. Then there is a path (of length at most four) from $\mathcal{P}$ to $\mathcal{P}'$ whose edges are finite [respectively decidable] Tietze transformations.

 Suppose now that instead we join two presentations if one can be obtained from the other by application of a Tietze transformation of bounded area or length. All the presentations belonging to a same connected component define the same group and moreover by Proposition \ref{dasz} they have the same Dehn function up to equivalence. However in this case two different connected components can define the same group as we have seen in Remark \ref{chuck}.

  \bigskip 
    
 Finally the Dehn function is related also to the \textit{hyperbolicity} of a group.
 
   \begin{definition} \label{hyper} \rm Let $\langle \, X \, | \, R \, \rangle$ be a finite presentation of a group $\mathcal{G}$. The group $\mathcal{G}$ is said \textit{hyperbolic} if there exists a positive real constant $\alpha$ such that $\Delta(n) \leqslant \alpha n$ (see \cite{CDP}).     
         \end{definition}

By Proposition \ref{dasz}, if the inequality $\Delta(n) \leqslant \alpha n$ holds in a finite presentation of $\mathcal{G}$, then it holds in every finite presentations of $\mathcal{G}$.

\section{Infinite presentations of groups} \label{ipog}

In this section we consider a function (called \textit{work}) which plays for finitely generated decidable presentations the same rôle played by the Dehn function for finite presentations. This function has been introduced by J.-C. Birget in (\cite{Birg}, 1.4) in the more general context of semigroup presentations and was studied in a recent paper by R. I. Grigorchuk and S. V. Ivanov \cite{GI}. As done for the Dehn function, we define a kind of Tietze transformations, said \textit{of bounded work}, that when applied to a presentation give a presentation with the same work function. This improves a result of \cite{GI}.

\medskip

We recall that $\rho$ is the reduced form of a word and $\mathcal{F}(X)$ is the free group on $X$.

\begin{definition} \label{omega} \rm Let $\mathcal{P}:=\langle \, X\, | \, R \, \rangle$ be a group presentation and let $w$ be a relator. The \textit{work of $w$} is defined as
    $$\mathrm{Work}(w)=\min\{|r_1|+\cdots+|r_n| \, : \, \rho(a_1 r_1 a_1^{-1} \cdots a_n r_n a_n^{-1})=w, $$
    $$r_1, \cdots, r_n \in R, \, \, a_1, \cdots, a_n \in \mathcal{F}(X) \}.$$
  Let $\mathcal{P}$ be finitely generated; then for a natural number $n$ we define 
  $$\Omega(n)=\max\{\mathrm{Work}(w)\, : \, w \, \,\textrm{is a relator and} \, \, |w|\leqslant n\}.$$ \end{definition}

The functions $L_1$ and $f_1$ of Definition \ref{buda} have been introduced by Grigorchuk and Ivanov in \cite{GI}.

\begin{definition} \label{buda} \rm Let $\mathcal{P}:=\langle \, X\, | \, R \, \rangle$ be a group presentation and let $w$ be a relator. We set $L_1(w)$ as the minimal number of edges of van Kampen diagrams with boundary cycle equal to $w$. If $n$ is a natural number we set 
  $$f_1(n)=\max\{L_1(w) \, : \, w \, \,\textrm{is a relator and} \, \, |w|\leqslant n\}. $$ \end{definition}

We prove that the functions $\Omega$ and $f_1$ are equivalent. Let $a_1, \cdots, a_n \in \mathcal{F}(X)$ and $r_1, \cdots, r_n \in R$ and let $w$ be the reduced form of $a_1 r_1 a_1^{-1} \cdots a_n r_n a_n^{-1}$. We set $E:=\big( (a_1, r_1), \cdots, (a_n, r_n) \big)$ and let $\mathcal{V}(E)$ be the van Kampen diagram relative to a freely reduction of $a_1 r_1 a_1^{-1} \cdots a_n r_n a_n^{-1}$ to $w$. We set
    $$\textrm{Work}(E):=|r_1|+\cdots+|r_n|$$
and we define $L_1(E)$ as the number of edges of $\mathcal{V}(E)$. Obviously Work$(w)$ and $L_1(w)$ are equal to the the minimal Work$(E)$ and $L_1(E)$ respectively where $E$ is an expression for $w$.

Let $l_1(E), l_2(E)$ and $l'_2(E)$ be the number of edges of $\mathcal{V}(E)$ belonging respectively to no faces, to a single face and to two faces of $\mathcal{V}(E)$. Then 
     $$L_1(E)=l_1(E) + l_2(E) + l'_2(E).$$
Since in computing the function Work we are interested in finding an expression for which $|r_1|+\cdots+|r_n|$ is minimal, then we can suppose that in constructing $\mathcal{V}(E)$ no 2-sphere has been discarded, that is we can suppose that
     $$\textrm{Work}(E)=l_2(E) + 2 l'_2(E).$$
This implies that     
         $$\textrm{Work}(E)\leqslant 2 L_1(E)$$
and since $l_1(E) < 2|w|$ then 
                   $$L_1(E)\leqslant \textrm{Work}(E) + 2|w|.$$
This means that 
                 $$\textrm{Work}(w)\leqslant 2 L_1(w), \hspace{1cm} L_1(w)\leqslant \textrm{Work}(w) + 2|w|$$
    and then that
 $$\Omega(n)\leqslant 2 f_1(n), \hspace{1cm} f_1(n)\leqslant \Omega(n) + 2n,$$
 that is the the functions Work and $f_1$ are equivalent.

 \medskip

We have seen that Proposition \ref{gersho} is no longer true for a presentation whose set of defining relators is infinite. But if we replace the Dehn function by the function $\Omega$ then we have analogous results for decidable presentations.

\begin{proposition} \label{risGI} \begin{enumerate}

  \item \textit{A finitely generated decidable presentation has a solvable Word Problem if and only if its function $\Omega$ is computable} (Th. 1.2 of \cite{GI});
  
  \item \textit{For a finite presentation the function $\Omega$ is equivalent to the Dehn function} (Th. 1.7 of \cite{GI}).
  
\end{enumerate} \end{proposition}

We have seen in Proposition \ref{dasz} that the Dehn function is invariant under bounded Tietze transformations. The. 1.6 of \cite{GI} shows that $\Omega$ is invariant under \textit{T-transformations} and \textit{stabilizations}, which are special cases of finite Tietze transformations. We will show that indeed $\Omega$ is invariant under a more general class of transformations, the so-called Tietze transformations \textit{of bounded work}, which are not necessarily finite. For that we introduce

\begin{definition} \label{Ttobw} \rm A Tietze transformation of type 1 or 2 is said of \textit{bounded work} if it adds or deletes a set $S$ of defining relators whose work is bounded, that is if $\sup\{\textrm{Work}(s) \, : \, s\in S\}$ is finite. \end{definition}

A finite Tietze transformation of type 1 or 2 is of bounded work. Conversely if the presentation is finitely related then a Tietze transformation of bounded work is necessarily finite. If the presentation is finitely generated then a Tietze transformation of bounded work is also of bounded area; if the presentation is finitely related then a Tietze transformation of bounded area is also of bounded work. This implies in particular that for a finite presentation a Tietze transformation is of bounded work if and only if it is of bounded area.

 \begin{remark} \label{ania} \rm Let $\mathcal{P}$ be a presentation and let $\mathcal{P}'$ be the presentation obtained from $\mathcal{P}$ by application of a Tietze transformation of bounded work or length. We want to determine the relationship between the area of relators in $\mathcal{P}$ and $\mathcal{P}'$.

   Suppose that $\mathcal{P}'$ is obtained from $\mathcal{P}$ by application of a Tietze transformation of type 1 of bounded work and let $c$ be the maximal area of defining relators added. If $\textrm{Work}$ and $\textrm{Work}'$ are the functions work relative to $\mathcal{P}$ and $\mathcal{P}'$ respectively, then we have for a relator $w$ that
              $$\textrm{Work}'(w) \leqslant \textrm{Work}(w) \leqslant c\textrm{Work}'(w).$$
  This implies that if $\Omega$ and $\Omega'$ are the functions as defined in Definition \ref{omega} relative to $\mathcal{P}$ and $\mathcal{P}'$ then
      \begin{equation} \label{T'1} \Omega' \leqslant \Omega \leqslant c\Omega'.\end{equation}

 Let us now suppose that $\mathcal{P}'$ is obtained from $\mathcal{P}$ by application of a Tietze transformation of type 3 of bounded length. By using the notation of 3 of Definition \ref{Transf}, we have that there exists a natural number $c$ such that $|u|\leqslant c$ for every $u\in U$. Let $\varphi$ be the function from $\mathcal{F}(X\cup Y)$ to $\mathcal{F}(X)$ sending any element of $X$ to itself and any $y_u\in Y$ to $u$. We have that $\varphi$ is the identity on $\mathcal{F}(X)$, in particular $\varphi(r)=r$ for every $r \in R$, and that $\varphi(t)=1$ for every $t \in T$. Moreover for every $w \in \mathcal{F}(X\cup Y)$ we have that $|\varphi(w)| \leqslant c |w|$.

 If we denote $\mathcal{N}$ and $\mathcal{N}'$ the normal subgroups of $\mathcal{F}(X)$ and $\mathcal{F}(X\cup Y)$ normally generated by $R$ and $R \cup T$ respectively, then we have that $\mathcal{N} \subset \mathcal{N}'$ and $\varphi(\mathcal{N}')=\mathcal{N}$. 
 
Let Work and Work$'$ be the functions work relative to $\mathcal{P}$ and $\mathcal{P}'$ respectively. Let $w \in \mathcal{N}'$ and let 
     $$w=\rho(a_1 r_1 a_1^{-1} \cdots a_k r_k a_k^{-1})$$ 
where $r_i \in R\cup T$ and $k=\textrm{Work}'(w)$. Then 
   $$\rho \big(\varphi(w)\big)=\rho(b_1 s_1 b_1^{-1} \cdots b_k s_k b_k^{-1})$$
where $b_i=\varphi(a_i)\in \mathcal{F}(X)$ and $s_i$ is equal to $r_i$ if $r_i \in R$, $s_i$ is equal to 1 if $r_i \in T$. This implies that 
         $$\mathrm{Work}\big(\varphi(w)\big) \leqslant \textrm{Work}'(w).$$
     If $w \in \mathcal{N}$ then since $\varphi(w)=w$ and since obviously $\mathrm{Work}'(w) \leqslant \textrm{Work}(w)$, then $\mathrm{Work}(w) = \textrm{Work}'(w)$. In particular we have that 
                      $$\Omega \leqslant \Omega'.$$
Let $w \in \mathcal{N}'$; since the length of any $t \in T$ is at most $c+1$, then there exists an $(R\cup T)$-derivation from $w$ to $\varphi(w)$ of work at most $(c+1) |w|$. Since there is an $(R\cup T)$-derivation from $\varphi(w)$ to 1 of work $\mathrm{Work}\big(\varphi(w)\big)$, then we have
    $$\mathrm{Work}'(w)\leqslant \mathrm{Work}\big(\varphi(w)\big) + (c+1)|w|,$$
    thus                     
            \begin{equation} \label{T'3} \Omega(n) \leqslant \Omega'(n)\leqslant \Omega(cn)+(c+1)n.\end{equation}   
                     \end{remark}

\begin{proposition} \label{} Let $\mathcal{P}$ and $\mathcal{P}'$ be group presentations and let $\Omega$ and $\Omega'$ be the functions as defined in Definition \ref{omega} relative to $\mathcal{P}$ and $\mathcal{P}'$ respectively. Let $\mathcal{P}'$ be obtained from $\mathcal{P}$ by applications of Tietze transformations of bounded work or length; then $\Omega$ and $\Omega'$ are equivalent. \end{proposition}

\begin{proof} By proposition \ref{Tietze}, $\mathcal{P}'$ can be obtained from $\mathcal{P}$ by applications of finitely many Tietze transformations. By $(\ref{T'1})$ and $(\ref{T'3})$, we have that applying a Tietze transformation of type 1 or 3, the equivalence class of the function $\Omega$ does not change. This is also true by applying a Tietze transformation of type 2 because its inverse is a transformation of type 1; and is true also for a Tietze transformation of type 4 since its inverse is a transformation of type 3 followed by one of type 1 and by one of type 2.    \end{proof}

\section{The Conjugacy Problem} \label{tcp}

In this section we define a function, denoted $\Gamma$, which plays for the Conjugacy Problem the same rôle played by the Dehn function (or the function work) for the Word Problem. It is surprising that despite the fact that the Conjugacy Problem is harder to solve than the Word Problem, anyway to obtain results analogous to those valid for the Word Problem, the only hypotheses needed is the finite generation of the presentation. 

No hypothesis is necessary on the set of defining relators, which can be even non-decidable or non-enumerable. It is also surprising that the function $\Gamma$ is constant under a much larger class of Tietze transformations than for the Dehn (or work) function. The results of this section are apparently new. 

\medskip

We recall that $\mathcal{F}(X)_c$ denotes the set of cyclically reduced words on $X\cup X^{-1}$. 

Let $\mathcal{G}$ be a group, let $X$ be a generating set for $\mathcal{G}$ and let $w_1, w_2\in \mathcal{F}(X)_c$ be conjugated in $\mathcal{G}$. We set
       \begin{equation} \label{gamma} \gamma(w_1, w_2):=\min\{|t| \,\, : t \in \mathcal{F}(X) \,\,\mathrm{and} \,\,w_1=t w_2 t^{-1} \,\, \textrm{in} \,\, \mathcal{G} \}.\end{equation}
   Let $X$ be finite and let $n_1, n_2$ be natural numbers: we set
       \begin{equation} \label{Gamma} \Gamma(n_1, n_2):=\max\{ \gamma(w_1, w_2) \,\, : w_1 \,\, \mathrm{and} \,\, w_2 \,\, \mathrm{are \,\, conjugated \,\, in \,\, \mathcal{G}} \end{equation}
        $$ \mathrm{and} \,\, |w_1|\leqslant n_1, \,\, |w_2|\leqslant n_2\}.$$
        
 The functions $\gamma$ and $\Gamma$ depend on the generating set $X$. If $\mathcal{P}$ is a presentation, the function $\Gamma$ \textit{relative to $\mathcal{P}$} is the function $\Gamma$ for the group presented by $\mathcal{P}$ with respect to the set of generators for $\mathcal{P}$.

\smallskip

Let $\mathcal{G}$ be generated by a set $X$, let $u$ be a word in $X\cup X^{-1}$ such that the element of $\mathcal{G}$ represented by $u$ is central; if $v$ is another word then $u$ and $v$ are conjugated in $\mathcal{G}$ if and only if they determine the same element of $\mathcal{G}$; in this case $\gamma(u, v)=0$. Thus if $\mathcal{G}$ is an abelian group then $\Gamma$ is equal to the zero function. The converse is also true, that is a group is abelian if and only if $\Gamma(n_1, n_2)=0$ for every $n_1$ and $n_2$. This result holds with respect to every generating set of $\mathcal{G}$.

If $\mathcal{G}$ is a finite group then $\Gamma$ is bounded. Indeed if $X$ is a generating set for $\mathcal{G}$ and if $k$ is the longest length of elements of $\mathcal{G}$ then $\Gamma$ is eventually equal to a constant less or equal to $k$.

On the other side, let $\mathcal{G}$ be generated by a finite set $X$ and let $\Gamma$ be bounded; then this implies that there exists a finite subset $\{g_1, \cdots, g_n\}$ of $\mathcal{G}$ such that for every $h$ and $h'$ which are conjugated in $\mathcal{G}$ and which are equal to cyclically reduced products of elements of $X\cup X^{-1}$, then there exists $i=1, \cdots, n$ such that $h'=g_i h g_i^{-1}$.

\smallskip

\textbf{Problem.} \textit{Find a finitely generated infinite non-abelian group whose function $\Gamma$ (with respect to a finite generating set) is bounded; or prove that if $\Gamma$ is bounded then $\mathcal{G}$ is finite or abelian.}
        
\smallskip

The next result is the analogous for the Conjugacy Problem of Proposition \ref{gersho}. It has to be noticed that unlike for the Word Problem, here we do not need that the number of defining relators be finite. 

\begin{proposition} \label{CP} A finitely generated presentation has a solvable Conjugacy Problem if and only if it has a solvable Word Problem and $\Gamma$ is bounded above by a computable function. In this case $\Gamma$ is computable. \end{proposition}

\begin{proof} Let $\mathcal{G}$ be the group presented by $\mathcal{P}$. \begin{enumerate}

  \item Let $\mathcal{P}$ have a solvable Word Problem and let $\Gamma(n_1, n_2)\leqslant h(n_1, n_2)$ where $h(n_1, n_2)$ is computable. Since $\mathcal{P}$ is finitely generated, then for any natural numbers $n_1$ and $n_2$ the number of pairs of words $(w_1, w_2)$ such that $|w_1|\leqslant n_1$ and $|w_2|\leqslant n_2$ is finite. Let $w_1, w_2$ be words such that $|w_1|\leqslant n_1$ and $|w_2|\leqslant n_2$; $w_1$ and $w_2$ are conjugated in $\mathcal{G}$ if and only there exists a word $t$ of length no more than $h(n_1, n_2)$ such that $w_1 t w_2^{-1} t^{-1}$ is a relator. 

Let us solve the Word Problem for all the words of the form $w_1 t w_2^{-1} t^{-1}$ where $|t|\leqslant h(n_1, n_2)$; this solves the Conjugacy Problem for the pair $(w_1, w_2)$. Moreover if $w_1$ and $w_2$ are conjugate in $\mathcal{G}$, then the minimal length of a conjugating word is equal to $\gamma(w_1, w_2)$ and thus $\Gamma(n_1, n_2)$ is computable since it is the maximum of the $\gamma(w_1, w_2)$.
  
  \item Let $\mathcal{P}$ be have a solvable Conjugacy Problem. Let $n_1, n_2$ be natural numbers; since $\mathcal{P}$ is finitely generated, then the number of pairs of words $(w_1, w_2)$ such that $|w_1|\leqslant n_1$ and $|w_2|\leqslant n_2$ is finite. Let us solve the Conjugacy Problem for these pairs of words; this gives for any pair of conjugated words $(w_1, w_2)$ an upper bound for $\gamma(w_1, w_2)$. The maximum of these bounds is an upper bound for $\Gamma(n_1, n_2)$ which is then bounded above by a computable function.  
  \end{enumerate}
\end{proof}

It would be interesting to compute the function $\Gamma$ for presentations with a solvable Conjugacy Problem. For instance for a finite presentation $\mathcal{P}$ satisfying the small cancellation condition $C'(1/8)$ then $\Gamma(n_1, n_2)\leqslant 6 \max(n_1, n_2)$. If $\mathcal{P}$ satisfies $C'(1/6)$ then $\Gamma(n_1, n_2)\leqslant 12 \max(n_1, n_2)$ and if $\mathcal{P}$ satisfies $C'(1/4)$ and $T(4)$ then $\Gamma(n_1, n_2)\leqslant 16 \max(n_1, n_2)$ (see The. V.5.4 of \cite{LS}). More generally, if $\mathcal{P}$ is a finite presentation of an hyperbolic group, then there exists a constant $K$ such that $\Gamma(n_1, n_2)\leqslant K \max(n_1, n_2)$ (see III-$\Gamma$-2.11 and 2.12 of \cite{BriHae}).

\begin{proposition} \label{new} Let $\mathcal{P}$ and $\mathcal{P}'$ be group presentations and let $\Gamma$ and $\Gamma'$ be the functions as defined in (\ref{Gamma}) relative to $\mathcal{P}$ and $\mathcal{P}'$ respectively. Suppose that $\mathcal{P}'$ is obtained from $\mathcal{P}$ by applications of Tietze transformations of type 3 or 4 of bounded length and by any Tietze transformations of type 1 or 2 (in particular, let $\mathcal{P}$ and $\mathcal{P}'$ be finitely generated presentations for the same group). Then $\Gamma$ and $\Gamma'$ are strongly equivalent (Definition \ref{equiv}). \end{proposition}

\begin{proof} Let $\mathcal{G}$ be the group presented by $\mathcal{P}$ and $\mathcal{P}'$. Since $\Gamma$ is independent from the defining relators of a presentation, then by applying Tietze transformations of type 1 or 2 the function $\Gamma$ does not change. Let us apply a Tietze transformation of length bounded by a natural number $m$ and let $\gamma$ and $\gamma'$ be the functions as defined in (\ref{gamma}) relative to $\mathcal{P}$ and $\mathcal{P}'$ respectively. By (\ref{hagen}) we have that 
  $$\gamma'(w_1, w_2) \leqslant \gamma(w_1, w_2) \leqslant m\gamma'(w_1, w_2)$$
 for every $w_1, w_2\in \mathcal{F}(X)_c$ conjugated in $\mathcal{G}$,
 thus  
          $$\Gamma' \leqslant \Gamma \leqslant m \, \Gamma'$$
 and that proves the claim.  \end{proof}

 \section{Enumeration of relators} \label{enume}

Let $\mathcal{P}:=\langle \, X \, | \, R \, \rangle$ be a finite presentation and let $k$ and $n$ be natural numbers. Currently there are two known methods to compute all the relators of $\mathcal{P}$ of area and length at most $k$ and $n$: the first consists in constructing all the van Kampen diagrams with at most $k$ faces and with boundary of at most $n$ edges; the second uses the result of Remark \ref{smer} which states that if $m$ is the maximal length of elements of $R$ and if $w$ is a relator of area $k'\leqslant k$ and length $n' \leqslant n$, then $w$ is a product of $k'$ conjugates of defining relators, where the length of the conjugating elements is less or equal to $mk'+n'$. 

 \medskip

Let us consider the first method. The way of proceeding is the following: one first constructs all the van Kampen diagrams  with at most $k$ faces and with boundary of at most $n$ edges and which have no \textit{spines}, i.e., in which every edge belongs to a face. Then one considers all the van Kampen diagrams with at most $k$ faces and with boundary of at most $n$ edges and which are obtained by \textit{adjoining} (see \cite{adjunc}) one by one two or more of these \textit{spineless} diagrams and some (even zero) labeled edges in such a way that any adjunction is made on a single vertex.

Spineless van Kampen diagrams are constructed in this way: one starts with a face and then adds one by one another face by identifying a vertex of the face with a vertex of the boundary of the diagram already constructed; or by identifying one or more consecutive edges of the face with one or more consecutive edges of the boundary of the diagram already constructed.

We recall that any face of a van Kampen diagram is labeled by a defining relator and that since the identification of two edges corresponds to a cancellation of two letters between the labels of the boundaries, then the edges to be identified must be labeled by opposite letters.

It must be observed that with this method we have to construct also diagrams whose boundary label is not reduced. This is because there can exist diagrams whose boundary label is reduced but which cannot be adjunctions of two diagrams with reduced boundary labels.

Let us compute the complexity of this method expressed as the number of diagrams to be constructed. We will count the number of diagrams without spines and which are obtained by adjunctions only along vertices. The number of such diagrams will result to be more than factorial in $k$.

Let $q$ be the minimal length of a defining relator; suppose that we have a face $F_1$ and that we want to add a second face $F_2$. We have at least $q^2$ ways to do it, any corresponding to the identification of a given vertex of $F_1$ with one of $F_2$. These diagrams have boundary with length equal at least to $2q$ and there are at least $2q \cdot q= 2q^2$ ways to add a third face $F_3$, therefore there are at least $2q^2 \cdot q^2=2q^4$ diagrams with three faces. These diagrams have length $3q$ and thus the diagrams with four faces are at least $2q^4 \times 3q=3! \, q^5$. By continuing this way we see that for a natural number $j$, there are at least $(j-1)! \, q^{2(j-1)}$ spineless diagrams obtained by adjunctions along vertices.

Since we want to obtain diagrams with boundary of length at most $n$, then for $j=1, \cdots, k$ we are not interested in those diagrams with $k-j$ faces and whose boundary is longer than $n+jm$ (we recall that $m$ is the maximal length of a defining relator). If $j$ is such that $n+mj\leqslant (k-j)q$, anyway all the spineless diagrams obtained by adjunctions along vertices must be computed. This gives as condition that 
                    $$j \leqslant \frac{q}{m+q} k - \frac{n}{m+q}=c_1 k - c_2 n,$$
with $c_1$ and $c_2$ positive constants less than 1. We recall that $k$ is the number of faces of the diagrams, thus it is reasonable to take $k\geqslant \Delta(n)$, where $\Delta(n)$ is the Dehn function at $n$. If $\Delta(n)/n$ goes to infinity then we can assume that $c_1 k - c_2 n$ has the same asymptotic behavior of $k$ and the number of these diagrams is asymptotically equal at least to 
       \begin{equation} \label{schi} (k-1)! \, q^{2(k-1)}, \end{equation}
that is it is a product of a factorial and an exponential in $k$.

If $\Delta(n)/n$ does not go to infinity then the Dehn function is linear in $n$ and the group is hyperbolic (Definition \ref{hyper}). These groups are already well studied in the literature.

We observe that (\ref{schi}) does not take into account spineless diagrams obtained by adjunctions along vertices and with a number of faces different than $c_1 k - c_2 n$; and does not take into account the diagrams obtained by adjoining along edges and also those containing spines. Thus the complexity of this method must be much bigger than (\ref{schi}).

\medskip

Now let us consider the second method. We want to compute its complexity expressed as the number of words to be constructed to find the desired relators. Let $N$ be this number. To determine $N$ we have to compute the number of products of $h$ conjugates of defining relators for $h=1, \cdots, k$ where each conjugating element has length bounded above by $mh+n$.

Let $U$ be the set of reduced words of length bounded above by $mk+n$ and let $p:=|R|$. Then we have 
   $$N:=\sum_{h=1}^k (p|U|)^h=\frac{(p|U|)^{k+1}-p|U|}{p|U|-1}. $$
Let us compute $|U|$. Let $a=2|X|-1$; then we have that     
    $$|U|=1+(a+1)\sum_{i=0}^{mk+n-1}a^i.$$
If we set $M:=mk+n$ then 
     $$\sum_{i=0}^{mk+n-1}a^i=\frac{a^M-1}{a-1}$$ 
  and thus 
       $$|U|=1+(a+1)\frac{a^M-1}{a-1}.$$
    For simplicity we make the following approximation:
           $$|U|\approx \frac{a+1}{a-1} \, a^M.$$
We have that $N \approx (p|U|)^k$ and since $M=mk+n$ then
   \begin{equation} \label{frot} N \approx \Big(p\frac{a+1}{a-1}\Big)^k  a^{mk^2+kn}=q^k a^{kn} s^{k^2},\end{equation}
 where $q=p (a+1)/(a-1)$ and $s=a^m$. We recall that $p=|R|$, $a=2|X|-1$ and $m$ is the maximal length of the elements of $R$.

The complexity of this algorithm (expressed as the number of words among which looking for to find the desired relators) is the product of three exponentials: one in $k$, one in $kn$ and one in $k^2$. This complexity is less than (\ref{schi}), thus this algorithm is better than the first one.

   \smallskip
   
Given natural numbers $k$, $m$ and $n$, it can be quite difficult to find a finite presentation whose maximal length of defining relators is $m$ and with a relator $w$ of length $n$, of area $k$ and which is the reduced form of a product of $k$ conjugates of defining relators where the conjugating elements have length $mk+n$. Anyway the following method can go very near to such a situation.

Let $\mathcal{F}(X)$ be the free group on the set $X$, let $r_1, \cdots r_{k-1}$ be cyclically reduced words in $X$ (not necessarily distinct), let $m$ be a natural number such that $m\geqslant |r_i|$ for $i=1, \cdots, k-1$ and let $a_1, \cdots, a_{k-1}$ be reduced words of length $mk+n$. Let $w$ be the reduced form of $a_1 r_1 a_1^{-1} \cdots a_{k-1} r_{k-1} a_{k-1}^{-1}$, let $v$ be a word of length $n$, set $r_k:=w^{-1}v$ and set $R:=\{r_1, \cdots, r_k, r_1^{-1}, \cdots, r_k^{-1}\}$. Then in the presentation $\langle \, X \, | \, R \, \rangle$, we have that $v$ is a relator of length $n$ and is the reduced form of $a_1 r_1 a_1^{-1} \cdots a_k r_k a_k^{-1}$.

 \chapter{The sets $C$ and $L$ and the function $A$} \label{tsral}       
  
In this chapter we present the main tools and objects of the thesis. In Section \ref{recdef} we introduce a class of algorithms called \textit{straight line algorithms} by means of which we define and study in Sections \ref{second} and \ref{potc} some special relators called \textit{corollas} and in Section \ref{L} the set $L$ and the functions $A$ (which will be shown in the Main Theorem \ref{maithe} to coincide with the set of all relators and the function \textit{area} respectively). Finally in Section \ref{statres} we state the main results of the thesis.

  \section{Straight line algorithms} \label{recdef}

In this section we introduce the main tool of which we make use in this thesis: \textit{straight line algorithms}. They are algorithms which generate elements of a given set $U$ starting from a set $B\subset U$ and applying some given operations. As observed in Remark \ref{graphSLA}, a straight line algorithm is an algorithm such that there is one and only one path from a given step to the final one; in particular there are no cycles. In Section \ref{rdsaslp} it will be shown that they are a natural tool for describing recursively defined sets.

In the current literature the term “straight line algorithm” denotes a slightly different object (see for instance 3.1.3 of \cite{Holt})

\smallskip   \smallskip

\begin{definition} \label{SLA} \rm Let $U$ be a set and let $\Phi$ be a family of functions $\varphi : U^{n_\varphi} \rightarrow U$ (where $n_\varphi$ is a given non-zero natural number depending on $\varphi$) with codomain $U$ and with domain some Cartesian power of $U$. A \textit{straight line algorithm (or SLA) relative to $(U, B, \Phi)$} is a finite algorithm in which a step can be either an element $b$ of $B$ (in this case $b$ is the output of the step and there is no input) or the application of a function $\varphi \in \Phi$ to $n_\varphi$ outputs $t_1, \cdots, t_{n_\varphi}$ of preceding steps \big($t_1, \cdots, t_{n_\varphi}$ are the inputs and $\varphi(t_1, \cdots, t_{n_\varphi})$ is the output\big). We require also the conditions that $t_1, \cdots, t_{n_\varphi}$ be outputs of distinct steps and that the output of any step, except the last, be an input of one and only one of the successive steps. $U$ is called the \textit{universe set}, $B$ the \textit{base set} and $\Phi$ is called \textit{the set of operations}. The elements of $B$ are called \textit{base elements}, a step equal to an element of $B$ is called a \textit{a base step}. \end{definition}

This notion of straight line algorithm is used very often in mathematics. For instance let $U$ be a group, let $B$ be a subset of $U$ and let $\Phi$ be equal to the operations of product and inversion in $U$. Let $a, b, c$ be elements of $B$ (not necessarily distinct); the following is an SLA relative to $(U, B, \Phi)$:
$$\textrm{step 1}: a \,\,\,\ \,\,  \textrm{step 2}: b \,\,\,\,\,\ \textrm{step 3}: a b  \,\,\,\,\,\ \textrm{step 4}: c  \,\,\,\,\,\, \textrm{step 5}: c^{-1} \,\ \,\,\,\, \textrm{step 6}: a b c^{-1}.$$
The steps 1, 2 and 4 are equal to elements of $B$, the 3 applies the product (which is an operation of $\Phi$) to outputs of preceding steps (the 1 and 2), the 5 applies the inversion to the output of the step 4 and finally the 6 applies the product to the outputs 3 and 5. The output of any step (except the last) is input of one and only one of the successive steps. The result of any of these straight line algorithms is an element of the subgroup of $U$ generated by $B$. If we add to $\Phi$ also the operation of conjugation with an element of $U$ then the result is an element of the normal subgroup generated by $B$.

\begin{definition} \label{} \rm Let $\sigma$ be a straight line algorithm and let $s$ and $s'$ be steps such that the output of $s$ is one of the inputs of $s'$; we say that \textit{$s'$ depends directly on $s$} or that \textit{$s'$ uses directly $s$}. The latter is a relation in the set of steps of $\sigma$. We call \textit{relation of dependence} the transitive reflexive closure \cite{reclo} of this relation; that is, given steps $s$ and $s'$, we say that \textit{$s'$ depends on $s$} if $s=s'$ or if there exists a finite sequence of steps starting with $s$, ending with $s'$ and such that every step in the sequence depends directly on the preceding. If $s'$ depends on $s$ we can say also that \textit{$s'$ uses $s$}.  \end{definition}

The first step of a straight line algorithm is always a base step because it cannot use preceding steps. The final output of a straight line algorithm is its \textit{result}. By abuse of notation we will sometimes identify a step with its output, but this will not cause ambiguity.

\begin{proposition} \label{mudepends}  Let $\sigma$ be an SLA. Then: 

\begin{enumerate}

  \item the steps using a given step of $\sigma$ form a chain with respect to the relation of dependence;

  \item the final step of $\sigma$ uses every step;
  
 \item if a step $s''$ uses a step $s$, then there exists a step $s'$ used by $s''$ and using directly $s$;

  \item let $s$ and $s'$ be steps such that $s'$ uses $s$ and let $s'_1, \cdots, s'_n$ be all the steps used directly by $s'$; then one and only one of the $s'_i$ uses $s$.

\end{enumerate} \end{proposition}

\begin{proof}  \begin{enumerate}

  \item Take a step $s$ of $\sigma$ and let $s'$ be a step using $s$. There exists a chain of steps $s_0:=s, s_1, s_2, \cdots$ such that $s_i$ uses directly $s_{i-1}$; since $s_i$ is the unique step using directly $s_{i-1}$, this chain is unique. Since $s'$ uses $s$, then necessarily $s'$ is one of the $s_i$.

  \item Let $s$ be a step of $\sigma$. Since the chain of steps of $\sigma$ using $s$ is finite (being $\sigma$ a finite algorithm), this chain ends necessarily with the last step of $\sigma$, which therefore uses $s$.

\item By Part 1 there exists a chain of steps $s_1:=s, s_2, \cdots, s_m:=s''$ such that $s_i$ uses directly $s_{i-1}$. Thus $s':=s_2$ uses directly $s$ and is used by $s''$.

  \item There exist steps $s_0:=s, s_1, \cdots, s_m:=s'$ such that $s_i$ uses directly $s_{i-1}$. Since $s'$ uses directly $s_{m-1}$ then $s_{m-1}$ is one of the $s'_i$ and uses $s_0=s$.

Suppose that $s'_j$ and $s'_k$ are two different steps used directly by $s'$ and using $s$. By Part 1 the steps using $s'$ form a chain and therefore $s'_j$ uses $s'_k$ (or $s'_k$ uses $s'_j$). By Part 3 there exists a step used by $s'_j$ (therefore preceding it) and using directly $s'_k$. This step cannot be $s'$ because $s'$ uses $s'_j$ and therefore follows it; this is impossible because $s'$ is the only step using directly $s'_k$.
                \end{enumerate} \end{proof}

\begin{remark} \label{graphSLA} \rm There is a natural way to associate a directed graph with an algorithm: the vertices of this graph are the steps of the algorithm and there is an edge directed from a step $s_1$ to a step $s_2$ if the output of $s_1$ is one of the inputs of $s_2$. 

The graph associated with a straight line algorithm is such that for every vertex there is one and only one path beginning at that vertex and ending at the vertex corresponding to the final step\footnote{In particular this graph has no cycles and then is a \textit{tree}.}. Furthermore this property characterizes straight line algorithms in the class of finite algorithms. 
                                   \end{remark}

\begin{remark} \label{pSLsA} \rm Let $\sigma$ be an SLA and let $s$ be one of its steps. It is easy to see that the steps used by $s$ form a straight line algorithm. This SLA is called \textit{the proper straight line subalgorithm (pSLsA) determined by $s$}. Every base element of a pSLsA of $\sigma$ is also a base element of $\sigma$.             
\end{remark}

In Definition \ref{SLsA} we will generalize this notion of proper straight line subalgorithm; this explains the use here of the adjective “proper".

\begin{definition} \label{multiset} \rm  A \textit{multiset} is a set whose elements can be repeated; it is defined as a pair $(S, \lambda)$ where $S$ is a set and $\lambda$ is a function from $S$ to the natural numbers. The value of $\lambda$ on an element of $S$ is the \textit{multiplicity} of that element. \end{definition}

\begin{definition} \label{basele} \rm Let $\sigma$ be an SLA relative to $(U, B, \Phi)$, let $\lambda: B \rightarrow \mathbb{N}$ be the function such that for $b\in B$, $\lambda(b)$ is the number of steps of $\sigma$ equal to $b$. The multiset $(B, \lambda)$ is called \textit{the multiset of base elements of $\sigma$}. \end{definition}

\begin{definition} \label{unmult} \rm Given two multisets $M_1:=(S, \lambda_1)$ and $M_2:=(S, \lambda_2)$ we define their union as $M_1\cup M_2:=(S, \lambda_1 + \lambda_2)$, that is the multiplicy of an element of $S$ in $M_1\cup M_2$ is the sum of its multiplicities in $M_1$ and $M_2$.  \end{definition}

\begin{remark} \rm \label{lambda} Let $\sigma$ be an SLA and let $s$ be one of its steps. Let $s_1:=s$, $s_2, \cdots, s_m$ be the chain of steps of $\sigma$ depending on $s$ (see Part 1 of Proposition \ref{mudepends}); in particular $s_i$ depends directly on $s_{i-1}$ and $s_m$ is the last step of $\sigma$. We can reorder the steps of $\sigma$ in such a way that $s_1$ depends on every step preceding it and that for $i=2, \cdots, m$, $s_i$ depends on every step comprised between $s_{i-1}$ and $s_i$. This reorder of the steps of $\sigma$ does not change the relative order of $s_1, \cdots, s_m$ and obviously does not change the result and the multiset of base elements.
                  \end{remark}

\section{The set $C$ and the function $\eta$} \label{second}

In this section we introduce a very important object for this thesis, the so-called \textit{corollas}. Corollas together with \textit{stems}, i.e., words of the form $w w^{-1}$ are the “ingredients” to construct the set $L$ defined in the next section and which will be shown in the Main Theorem \ref{maithe} to coincide with the set of all relators. Thus corollas can be considered the “non-trivial” parts of the relators.

\smallskip

\begin{definition} \label{insert} \rm Let $w_1, w_2$ and $u$ be words and let $w:=w_1 w_2$. The word $w_1 u w_2$ is called \textit{the insertion of $u$ into $w$ at $w_1$}. If $n=|w_1|$, the word $w_1 u w_2$ is also called \textit{the n-th insertion of $u$ into $w$} or the \textit{insertion of $u$ into $w$ at the n-th letter}. If $n \geqslant |w|$ we define the $n$-th insertion of $u$ into $w$ as the product $wu$.
      \end {definition}

The $n$-insertions are binary operations in $\mathcal{M}(X \cup X^{-1})$.

\begin{definition} \label{stem} \rm Let $\mathcal{S}$ be the set of words of the form $w w ^{-1}$ where $w$ is reduced and $w\neq 1$. We call \textit{stems} the elements of $\mathcal{S}$. \end{definition}

By Definition \ref{reducprod}, the product of two words $u$ and $v$ is an insertion of cancelled part (which is a stem) into $\rho(uv)$.

\begin{remark} \label{shirv} \rm Let $u$ and $v$ be reduced words and let $w=\pi(u, v)$ be the cyclically reduced product of $u$ by $v$. Let us study the relation between $u$, $v$ and $\pi(u, v)$. As in Definition \ref{reducprod} let $u_1, v_1, a\in \mathcal{F}(X)$ be such that $u=u_1 a$, $v=a^{-1} v_1$ and $\rho(uv)=u_1 v_1$, that is $uv=u_1 a a^{-1} v_1$. If $w=1$ then $t=1$ because the word $u_1 v_1$ is reduced and in this case $u=v^{-1}$. Let $w\neq 1$; since $u_1 v_1=t w  t^{-1}$, three cases are possible: \begin{enumerate}

  \item $u_1$ is a prefix of $t$;
  
  \item $u_1$ is a prefix of $tw$ but not of $t$;
  
  \item $u_1$ is not a prefix of $tw$.

\end{enumerate}

Let us examine the three cases.

\begin{enumerate}

  \item there exists a word $t_1$ such that $t=u_1 t_1$. Therefore $v_1=t_1 w t^{-1}=t_1 w t_1^{-1} u_1^{-1}$ and $u=u_1 a$, $v=a^{-1} t_1 w t_1^{-1} u_1^{-1}$, thus $uv=u_1 a a^{-1} t_1 w t_1^{-1} u_1^{-1}$. Moreover $v^{-1}$ is a prefix of $v$. 
  
  \item since $u_1$ is not a prefix of $t$, $t$ is a prefix of $u_1$. Moreover $v_1$ is a suffix of $w  t^{-1}$. This means that there exist words $w_1$ and $w_2$ such that $w=w_1 w_2$ and $u_1=t w_1$ and $v_1=w_2 t^{-1}$. Therefore $u=t w_1 a$, $v=a^{-1} w_2 t^{-1}$ and thus $uv=t w_1 a a^{-1} w_2 t^{-1}$.
  
  \item since $u_1$ is not a prefix of $tw$, $tw$ is a prefix of $u_1$ and thus there exists a word $t_1$ such that $u_1=twt_1^{-1}$ and $t^{-1}=t_1^{-1} v_1$. Therefore $u_1=v_1^{-1} t_1wt_1^{-1}$ and $u=v_1^{-1} t_1 w  t_1^{-1} a$, $v=a^{-1} v_1$ and thus $uv=v_1^{-1} t_1 w  t_1^{-1} a a^{-1} v_1$. Moreover $v^{-1}$ is a suffix of $u$.
  
\end{enumerate} 

In Case 2 if $w_1=1$ then we obtain the Case 1, if $w_2=1$ the Case 3, therefore we can suppose that in Case 2, $w_1$ and $w_2$ are non-empty. In Case 1, $u$ is completely cancelled in the reduced product by $v$, in Case 3 it is $v$ to be completely cancelled. In Case 2 no one of them is completely cancelled. If $a$ is empty there is no cancellation in the reduced product of $u$ by $v$; if also $t$ is empty there is no cancellation also in their cyclically reduced product. 
\end{remark}

\begin{remark} \label{puyi} \rm Let $u$ and $v$ be reduced words such that $v^{-1}$ is not a suffix (Definition \ref{presuf}) of $u$ and $u^{-1}$ is not a prefix of $v$. Then there exist a cyclic conjugate $u'$ of $u$ and a cyclic conjugate $v'$ of $v$ such that $\rho(u' v')=\pi(u,v)$. To prove it let us consider the second case of Remark \ref{shirv}. We have that $u=t w_1 a$, $v=a^{-1} w_2 t^{-1}$ and $\pi(u, v)=w$. Set $u':=w_1 a t$ and $v':= t^{-1}a^{-1} w_2$; then $\rho(u' v')=\pi(u,v)$.  \end {remark}

\begin{proposition} \label{rpcrp} Let $u$ and $v$ be reduced words and let $w=\pi(u, v)$ be the cyclically reduced product of $u$ by $v$; then $\pi(v, u)$ is a cyclic conjugate of $\pi(u, v)$.  
     \end{proposition}

\begin{proof} Let us prove the claim part by showing that it holds in the three cases of Remark \ref{shirv}.

\begin{enumerate}

  \item we have that $v=a^{-1} t_1 w t_1^{-1} u_1^{-1}$ and $u=u_1 a$, therefore $\rho(vu)=a^{-1} t_1 w t_1^{-1} a$ and $\pi(v, u)=w=\pi(u, v)$.
  
  \item we have that $v=a^{-1} w_2 t^{-1}$ and $u=t w_1 a$, therefore $\rho(vu)=a^{-1} w_2 w_1 a$ and $\pi(v, u)=w_2 w_1$ which is a cyclic of $\pi(u, v)=w_1 w_2$.
    
  \item we have that $v=a^{-1} v_1$ and $u=v_1^{-1} t_1 w  t_1^{-1} a$, therefore $\rho(vu)=a^{-1} t_1 w  t_1^{-1} a$ and $\pi(v, u)=w=\pi(u, v)$.
  
\end{enumerate} 
\end{proof}

For every word $w$ and for every natural number $n$, we let $\psi_n(w)$ denote the reduced form of the $n$-th cyclic conjugate of $w$ and we set    
        $$\Psi:=\{\psi_n : n\in\mathbb{N}^*\} \cup \{\pi\},$$ 
where $\pi$ has been introduced in Definition \ref{pi}.

 Straight line algorithms have been defined in Section \ref{recdef}.

\begin{definition} \label{overlineR} \rm Let $R$ be a subset of $\mathcal{F}(X)_c\setminus\{1\}$ containing the inverse of any of its elements and let $\Sigma$ be the set of $(\mathcal{F}(X)_c, R, \Psi)$-straight line algorithms $\sigma$ such that if $\pi(u, v)$ is a step of $\sigma$ then $\pi(u, v)\neq 1$ and there is at least one cancellation in $\pi(u, v)$, that is $|\pi(u, v)|< |u|+|v|$. We denote $\overline{R}$ the set of results of $\Sigma$ and we call \textit{generalized corollas} or \textit{g. corollas} the elements of $\overline{R}$. 

Let $\Sigma_0$ the subset of $\Sigma$ consisting of straight line algorithms $\sigma$ such that if $\pi(u, v)$ is a step of $\sigma$ then either $u$ or $v$ or both are cyclic conjugates of elements of $R$. We denote $C$ the set of results of $\Sigma_0$ and we call \textit{corollas} the elements of $C$. 
\end{definition}

We will call \textit{straight line algorithms in $\overline{R}$} or \textit{$\overline{R}$-straight line algorithms} the elements of $\Sigma$. The elements of $\Sigma_0$ are called \textit{straight line algorithms in $C$} or \textit{$C$-straight line algorithms}

\begin{proposition}\label{rem} Let $\mathcal{N}$ be the normal closure of $R$ in $\mathcal{F}(X)$, i.e., the intersection of all normal subgroups of $\mathcal{F}(X)$ containing $R$. Then $\mathcal{N}\supset \overline{R}$. \end{proposition}

\begin{proof} Let $\sigma$ be an SLA in $\overline{R}$, let $c$ be the result of $\sigma$ and let $n$ be the number of steps of $\sigma$. We will prove the claim by induction on $n$, being trivial for $n=1$.

Let the last step of $\sigma$ be a cyclic conjugation of a word $c'$ and let $\sigma'$ be the proper straight line subalgorithm (Remark \ref{pSLsA}) computing $c'$. There exist words $d$ and $e$ such that $c'=de$ and $c=e d$. Then $c=\rho(d^{-1}c' d)$ and belongs to $\mathcal{N}$ since $c'\in \mathcal{N}$.

Let the last step of $\sigma$ be the cyclically reduced product of two words $c_1$ and $c_2$ and let $\sigma_1$ and $\sigma_1$ be the proper straight line subalgorithms computing $c_1$ and $c_2$. Then $c$ is the reduced form of a conjugate of the product of $c_1$ by $c_2$ and thus belongs to $\mathcal{N}$.
                                                       \end{proof}

We recall that a base step (Definition \ref{SLA}) of a straight line algorithm is a step equal to a base element. A base step for a straight line algorithm in $\overline{R}$ is a step equal to a defining relators (the other two kinds of steps are the cyclic conjugations and the cyclically reduced products)

\begin{definition} \label{eta} \rm Let $\tau$ be a straight line algorithm in $\overline{R}$. We set $\eta(\tau)$ as the number of steps of $\tau$ equal to base steps.  \end{definition}

$\eta$ is a function going from the set of SLA's in $\overline{R}$ to the natural numbers. If $\tau$ is an SLA with only one step then this step is necessarily a base step and $\eta(\tau)=1$. Suppose that $\tau$ has more than one step. Let the final step of $\tau$ be the cyclic conjugation of a preceding output $c$. If $\tau'$ is the pSLsA of $\tau$ (Remark \ref{pSLsA}) computing $c$ then $\tau'$ has the same number of base steps as $\tau$ and thus $\eta(\tau)=\eta(\tau')$. Let the final step of $\tau$ be the cyclically reduced product of preceding outputs $c_1$ and $c_2$ and let $\tau_1$ and $\tau_2$ be the pSLsA's of $\tau$ computing $c_1$ and $c_2$; then the base steps of $\tau$ are those of $\tau_1$ and of $\tau_2$ and thus $\eta(\tau)=\eta(\tau_1)+\eta(\tau_2)$.

 \begin{definition} \label{eta(c)} \rm Let $c \in \overline{R}$; we set $\eta(c):=\min\{ \eta(\tau) :$ $\tau \, \,$ \textrm{is an SLA in $\overline{R}$ computing} $\, \, c\}$.  \end{definition}

\begin{proposition} \label{minoreuguale} Let $\tau$ be an SLA in $\overline{R}$ and let $c$ be its result. Then $\textrm{Area}(c) \leqslant \eta(\tau)$.
      \end{proposition}

\begin{proof} We prove the claim by induction on the number of steps of $\tau$. If $\tau$ has only one step then $\eta(\tau)=1$ and $c$ is a base element; this means that $c$ belongs to $R$ and therefore $\textrm{Area}(c)=1$. 

Let $\tau$ have more than one step and let the claim be true for every SLA with less steps than $\tau$. Let the final step of $\tau$ be the cyclic conjugation of a preceding output $c'$ and let $c$ be the final output of $\tau$; if $\tau'$ is the pSLsA of $\tau$ computing $c'$ we have $\eta(\tau)=\eta(\tau')$ by the construction of $\eta$ and $\textrm{Area}(c)=\textrm{Area}(c')$ by (\ref{uwu-1}). By induction hypothesis we have $\textrm{Area}(c') \leqslant \eta(\tau')$, therefore $\textrm{Area}(c) \leqslant \eta(\tau)$.

Let the final step of $\tau$ be the cyclically reduced product of preceding outputs $c_1$ and $c_2$, that is $c=\pi(c_1,c_2)$. If $\tau_1$ and $\tau_2$ are the pSLsA's of $\tau$ computing $c_1$ and $c_2$, we have $\eta(\tau)=\eta(\tau_1)+\eta(\tau_2)$ by the construction of $\eta$ and $\textrm{Area}(c) \leqslant \textrm{Area}(c_1)+\textrm{Area}(c_2)$ by (\ref{Area pi(v,w)}). By induction hypothesis we have $\textrm{Area}(c_1) \leqslant \eta(\tau_1)$ and $\textrm{Area}(c_2) \leqslant \eta(\tau_2)$, therefore $\textrm{Area}(c) \leqslant \eta(\tau)$.
 \end{proof}

\begin{corollary} \label{areaeta} If $c\in \overline{R}$ then $\textrm{Area}(c) \leqslant \eta(c)$.
      \end{corollary}

 \begin{proof} Follows from Definition \ref{eta(c)} and Proposition \ref{minoreuguale}.
 \end{proof}

\section{Properties of the corollas} \label{potc}

In this section we prove some properties of the corollas. In Remarks \ref{bau} and \ref{miao} we show how to associate a product of conjugates of defining relators and a van Kampen diagram with a straight line algorithm computing a corolla. Remark \ref{mj} says that if the presentation is finite and $m$ is the maximal length of a defining relator then a corolla of length $n$ and area $k$ can be expressed as product of conjugates of defining relators with the sum of the lengths of the conjugating elements bounded above by $(mk-1)/2+n-1$, while for a general relator this sum is bounded above by $mk^2+kn$, see Remark \ref{smer}. Remark \ref{Roverlinen} gives an algorithm for computing the corollas when the presentation is finite.

\medskip

 \begin{remark} \label{eta(psi)} \rm Let $c$ be a g. corolla and let $c'$ be a cyclic conjugate of $c$; then $\eta(c')=\eta(c)$. Indeed let $\tau$ be an SLA in $\overline{R}$ computing $c$ and such that $\eta(\tau)=\eta(c)$. Let $\tau'$ be the SLA in $\overline{R}$ obtained by adding to $\tau$ a cyclic conjugation from $c$ to $c'$; $\tau'$ computes $c'$ then $\eta(c')\leqslant \eta(\tau')$ and $\eta(\tau')=\eta(\tau)$. This implies that $\eta(c')\leqslant \eta(c)$. Since $c$ is a cyclic conjugate of $c'$ then in the same way we prove that $\eta(c)\leqslant \eta(c')$ and thus that $\eta(c')=\eta(c)$.
     \end{remark}

 \begin{remark} \label{eta(pi(cc'))} \rm Let $c$ and $c'$ be g. corollas such that $|\pi(c, c')|<|c|+|c'|$, let $\tau$ and $\tau'$ be SLA's in $\overline{R}$ computing $c$ and $c'$ respectively and such that $\eta(c)=\eta(\tau)$ and $\eta(c')=\eta(\tau')$. Let $\tau''$ be the SLA obtained by adding to $\tau$ the steps of $\tau'$ and finally a step equal to $\pi(c,c')$; $\tau''$ is an SLA in $\overline{R}$ (see Definition \ref{overlineR}). We have that $\eta(\tau'')=\eta(\tau)+\eta(\tau')$ and thus that $\eta(\tau'')=\eta(c)+\eta(c')$. Since $\eta\big(\pi(c,c')\big)\leqslant \eta(\tau'')$ then we have that $\eta\big(\pi(c,c')\big)\leqslant \eta(c)+\eta(c')$.
   \end{remark}

\begin{remark} \label{bau} \rm Let $\sigma$ be a straight line algorithm in $\overline{R}$; we show how to find a product of conjugates of $\eta(\sigma)$ defining relators whose reduced form is the result of $\sigma$.

Let $\sigma$ have one step; then this step is a defining relator $r$ and we associate $r$ with $\sigma$. Let $\sigma$ have more than one step and suppose to have proved the claim for every SLA with less steps than $\sigma$. Let the last step of $\sigma$ be the cyclic conjugation of a word $c'$ into a word $c$; then there exist words $d$ and $e$ such that $c'=de$ and $c=e d$. Let $\sigma'$ be the proper straight line subalgorithm (Remark \ref{pSLsA}) computing $c'$; by induction hypothesis 
    $$c'=\rho(a_1 r_1 a_1^{-1} \cdots a_m r_m a_m^{-1})$$
where $r_1, \cdots, r_m$ are defining relators, $a_1, \cdots, a_m$ are words and $m=\eta(\sigma')$. Then $m=\eta(\sigma)$ and
     $$c=\rho(b_1 r_1 b_1^{-1} \cdots b_m r_m b_m^{-1})$$
where $b_i=d^{-1}a_i$.

Let the last step of $\sigma$ be the cyclically reduced product of two words $c_1$ and $c_2$, that is the result of $\sigma$ is $c=\pi(c_1, c_2)$ (Definition \ref{pi}). Let $\sigma_1$ and $\sigma_2$ be the proper straight line subalgorithms computing $c_1$ and $c_2$; then by induction hypothesis
   $$c_1=\rho(a_1 r_1 a_1^{-1} \cdots a_k r_k a_k^{-1}), \,\,\,\,\,\,\,\, c_2=\rho(a_{k+1} r_{k+1} a_{k+1}^{-1} \cdots a_m r_m a_m^{-1})$$
   where $k=\eta(\sigma_1)$ and $m-k=\eta(\sigma_2)$.
 
  By Definition \ref{scope} there exist a word $t$ such that 
$$tct^{-1}=\rho(a_1 r_1 a_1^{-1} \cdots a_m r_m a_m^{-1}).$$
   and thus 
 $$c=\rho(b_1 r_1 b_1^{-1} \cdots b_m r_m b_m^{-1})$$
where $b_i=t^{-1}a_i$ and $m=\eta(\sigma)$.    \end {remark}

 \begin{remark} \label{miao} \rm Let $\sigma$ be a straight line algorithm in $\overline{R}$, let $k=\eta(\sigma)$ and let $c$ be its result. By Remark \ref{bau}, $c$ is the reduced form of a product of $k$ conjugates of defining relators. Thus we can associate with $c$ a van Kampen diagram with $k$ faces as seen in Section \ref{intro} but we show a more direct way. The procedure is analogous to that of Remark \ref{bau}. We will prove that the diagram associated with $\sigma$ is planar, contractible and homeomorphic to a 2-disc. In particular, unlike the general situation described in Section \ref{intro}, in the construction of this diagram no 2-sphere is discarded and the number of faces of the diagram is equal to the number of base steps of $\sigma$.

 Let $\sigma$ have one step; then this step is a defining relator $r=x_1\cdots x_m$ and we associate with $\sigma$ a complex with a single face whose boundary has $m$ edges labeled consecutively by $x_1,\cdots, x_m$. This diagram verifies trivially the desired properties.
 
 Let $\sigma$ have more than one step and suppose to have proved the claim for every SLA with less steps than $\sigma$. Let the last step of $\sigma$ be the cyclic conjugation of a word $c'$ into a word $c$; then there exist words $d$ and $e$ such that $c'=de$ and $c=e d$. Let $\sigma'$ be the proper straight line subalgorithm computing $c'$; by induction hypothesis there exists a diagram verifying the three properties above. The boundary cycle of this diagram is of the form $\delta \epsilon$ where $\delta$ is labeled by $d$ and $\epsilon$ by $e$. Then we associate with $\sigma$ the same van Kampen diagram with boundary label $\epsilon \delta$.

 Let the last step of $\sigma$ be the cyclically reduced product of two words $c_1$ and $c_2$, that is the result of $\sigma$ is $c=\pi(c_1, c_2)$. Let $\sigma_1$ and $\sigma_2$ be the proper straight line subalgorithms computing $c_1$ and $c_2$; then by induction hypothesis we can suppose that we have associated van Kampen diagrams $\mathcal{V}_1$ and $\mathcal{V}_2$ with $\sigma_1$ and $\sigma_2$. By Definition \ref{scope} there exists a word $t$ such that $\rho(c_1 c_2)=t c t^{-1}$ where $\rho(c_1 c_2)$ is the reduced product of $c_1$ by $c_2$. By Definition \ref{reducprod}, there exist words $a, c_1', c_2'$ such that $c_1=c_1' a$, $c_2=a^{-1} c_2'$ and $\rho(c_1 c_2)=c_1' c_2'$. First we associate a complex with $\rho(c_1 c_2)$ by folding the subpath of $\mathcal{V}_1$ labeled by $a$ onto the subpath of $\mathcal{V}_2$ labeled by $a^{-1}$; this complex is the \textit{adjunction} (see \cite{adjunc}) of $\mathcal{V}_1$ to $\mathcal{V}_2$ along the subgraph of $\mathcal{V}_1$ labeled by $a$ into the subgraph of $\mathcal{V}_2$ labeled by $a^{-1}$. Then we consider the quotient of the obtained complex given by identifying the edges in the path labeled by $t$ with the opposite of the edges of the path labeled by $t^{-1}$. The complex obtained is planar and homeomorphic to a 2-disc because it is the adjunction of two complexes with these properties. It is contractible by Ex. 23 in Chapter 0 of \cite{H} because it is the union of two contractible 2-complexes whose intersection is the path labeled by $t$ which is contractible.    \end {remark}

\begin{remark} \label{squit} \rm Let $\sigma$ be a straight line algorithm in $\overline{R}$; in Remark \ref{bau} we have shown how to associate with $\sigma$ a product of $\eta(\sigma)$ conjugates of defining relators whose reduced form is the result of $\sigma$. We show now another way to obtain the same product.

If $\sigma$ has one step, then this step is a defining relator $r$ and we associate $r$ with $\sigma$. Let $\sigma$ have more than one step and suppose that the last step of $\sigma$ is a cyclic conjugation of a word $c'$ into a word $c$; then there exist words $d$ and $e$ such that $c'=de$ and $c=e d$ and thus $c$ is the reduced form of $d^{-1} c d$.

Let the last step of $\sigma$ be the cyclically reduced product of two words $c_1$ and $c_2$, that is the result of $\sigma$ is $c=\pi(c_1, c_2)$. By Definition \ref{scope} there exist a word $t$ such that $tct^{-1}=\rho(c_1 c_2)$, that is $c$ is the reduced form of $t^{-1}c_1 t \, t^{-1} c_2 t$.

In both cases $c$ is the reduced form of a product of conjugates of g. corollas of $\sigma$. Any of these g. corollas is a cyclic conjugate of another g. corolla of $\sigma$ or the cyclically reduced product of two g. corollas of $\sigma$, then we can repeat the same procedure to these g. corollas. Continuing this way until all the g. corollas obtained are the base elements of $\sigma$, we obtain an expression of $c$ as product of conjugates of defining relators. This product obviously coincides with that defined in Remark \ref{bau}.      \end {remark}

\begin{remark} \label{mj} \rm Let $\sigma$ be a straight line algorithm in $\overline{R}$, let $c$ be the result of $\sigma$ and let  
    $$c=\rho(a_1 r_1 a_1^{-1} \cdots a_m r_m a_m^{-1})$$
where this expression is that obtained as in Remark \ref{bau}. We show that there exists a sequence of cancellations reducing $a_1 r_1 a_1^{-1} \cdots a_m r_m a_m^{-1}$ to $c$ such that for $i=1, \cdots, m$:\begin{enumerate}

  \item if $l$ is a letter of $a_i$ which cancels, then $l$ cancels either with a letter of some of the $r_j$ or with a letter of $a_{i-1}^{-1}$;
   
  \item if $l$ is a letter of $a_i^{-1}$ which cancels, then $l$ cancels either with a letter of some of the $r_j$ or with a letter of $a_{i+1}$.
  \end{enumerate}

In particular this implies that a letter of $a_1$ or one of $a_m^{-1}$ can cancel only with letters of some of the $r_j$.

We use Remark \ref{squit}. If $\sigma$ has only one step then the claim is trivial. Let $\sigma$ have more than one step; then by induction hypothesis $c$ is the reduced form of some product $b_1 c_1 b_1^{-1} \cdots b_n c_n b_n^{-1}$ where the $c_i$ are g. corollas of $\sigma$ and the $b_i$ verify the two properties above. Let $c_i$ for some $i$ be the output of a step of $\sigma$ consisting in a cyclic conjugation of a g. corolla $c'$. Then there exist words $d$ and $e$ such that $c'=de$ and $c=e d$, that is $c$ is the reduced form of $a_1 c'_1 a_1^{-1} \cdots a_n c'_n a_n^{-1}$, where $c'_i=c$, $a_i=b_i d^{-1}$ and $a_j=b_j$, $c'_j=c_j$ for $j\neq i$. The $a_i$ verify the two properties above because the suffix $d^{-1}$ of $a_i$ cancels with $c'_i$. 

Let $c_i$ for some $i$ be the output of a step of $\sigma$ consisting in a cyclically reduced product of two g. corollas $d_1$ and $d_2$, that is $c_i=\pi(d_1, d_2)$. By Definition \ref{scope} there exists a word $t$ such that $t c_i t^{-1}=\rho(d_1 d_2)$, that is $c_i$ is the reduced form of $t^{-1} d_1 t \, t^{-1} d_2 t$. In this case $c$ is the reduced form of $a_1 c'_1 a_1^{-1} \cdots a_{n+1} c'_{n+1} a_{n+1}^{-1}$, where $a_j=b_j$ and $c'_j=c_j$ for $j< i$, $a_i=a_{i+1}=b_i t^{-1}$, $c'_i=d_1$, $c'_{i+1}=d_2$ and $a_j=b_{j-1}$ and $c'_j=c_{j-1}$ for $j>i$. The $a_i$ verify the two properties above because the suffix $t^{-1}$ of $a_i$ cancels with $c'_i$, $a_i^{-1}$ cancels with $a_{i+1}$ and the prefix $t$ of $a_{i+1}^{-1}$ cancels with $c'_{i+1}$. 
\end {remark}      

\begin{remark} \label{annie} \rm Let $\sigma$ be a straight line algorithm in $\overline{R}$, let $c$ be the result of $\sigma$ and let  
    $$c=\rho(a_1 r_1 a_1^{-1} \cdots a_k r_k a_k^{-1})$$
where this expression is that obtained as in Remark \ref{bau}. Let 
     $$\alpha=|a_1|+\cdots+|a_k|, \hspace{1cm} \beta=|r_1|+\cdots+|r_k|.$$
  We show that either $\alpha<\beta$ or 
        $$\alpha\leqslant (\beta-1)/2 + |c|-1.$$ 
    In particular let $m$ be the maximal length of an element of $R$ and let $|c|\leqslant n$; then if $\alpha\geqslant \beta$ this implies that 
        $$\alpha\leqslant (mk-1)/2 + n-1.$$

Let $\Lambda$ be the set of letters of $a_1 r_1 a_1^{-1} \cdots a_k r_k a_k^{-1}$ which are cancelled in the free reduction to $c$ and define a function 
       $$\varphi : \Lambda \rightarrow \Lambda $$
   such that for every $l \in \Lambda$, the letter $l$ cancels with $\varphi(l)$ and the cancellations verify the two properties of Remark \ref{mj}. The function $\varphi$ is a bijection and $\varphi^2=id$.     

Consider the set of sequences $(l_1, l_2, \cdots, l_{2h})$ of elements of $\Lambda$ such that:\begin{enumerate}

  \item $l_{2i}=\varphi(l_{2i-1})$ for $i=1, \cdots, h$;
  
  \item $l_i$ belongs to some of the $a_j^{\pm 1}$ for $i=2, \cdots, 2h-1$;
  
  \item $l_{2i+1}=l_{2i}^{-1}$ for $i=1, \cdots, h-1$;
  
   \item either $l_1$ is a letter of some of the $r_j$ or it is a letter of some of the $a_j^{\pm 1}$ and $l_1^{-1}$ does not belong to $\Lambda$ (i.e., it is not cancelled).
   
    \item either $l_{2h}$ is a letter of some of the $r_j$ or it is a letter of some of the $a_j^{\pm 1}$ and $l_{2h}^{-1}$ does not belong to $\Lambda$ (i.e., it is not cancelled).
   
\end{enumerate} 

By condition 2, the $l_{2i}$ of condition 3 belong to some of the $a_j^{\pm 1}$, therefore condition 3 says that there exists $j:=1, \cdots, k$ such that $a_j=bxc$ for words $b$ and $c$ and a letter $x$ and $\{l_{2i}, l_{2i+1}\}=\{x, x^{-1}\}$.

If a sequence $(l_1, l_2, \cdots, l_{2h})$ verifies the four properties above then also does the sequence $(l_{2h}, l_{2h-1}, \cdots,l_1)$.

We observe that any sequence contains at most one letter of a given $a_j$ or $a_j^{-1}$: indeed let $l_i$ be a letter of $a_j$. If $i$ is odd, then $l_{i+1}$ is the letter cancelled with $l_i$. By Remark \ref{mj}, either $l_{i+1}$ belongs to some of the $r_j$ and thus $l_{i+1}$ is the last element of the sequence or $l_{i+1}$ belongs to $a_{j-1}^{-1}$. This means that the indices of the $a_j$ which have letters in the sequence are strictly decreasing and in particular there is at most one letter for any of the $a_j$. The same is true if $l_i$ is a letter of $a_j^{-1}$, in which case the indices are strictly increasing. If $i$ is even then $l_{i+1}$ is a letter of $a_j^{-1}$, thus $i+1$ is odd what said above applies.

Moreover if two sequences have an element in common then one is obtained from the other by reversing the order of its elements. Any letter cancelled uniquely determines the first and the last letters of the sequence to which it belongs. If a sequence contains the letter $x$ of some $a_j=bxc$ and the letter $x^{-1}$ of $a_j^{-1}=c^{-1}x^{-1}b^{-1}$ does not belong to the sequence, then $x^{-1}$ is not cancelled. Otherwise, if $x$ and $x^{-1}$ are cancelled then they belong to the same sequence and that sequence is uniquely determined by its first and last element. Therefore corresponding to $x$ and $x^{-1}$ there are two letters which are a cancelled letter of some of the $r_j$ and a non-cancelled one of some of the $a_j^{\pm 1}$ or there are two letters which are both of one of these two types. 

This means that the sum of cancelled letters of the $a_j^{\pm 1}$ is no more than the sum of non-cancelled letters of the $a_j^{\pm 1}$ plus the cancelled letters of the $r_j$.

Let $\alpha=|a_1|+\cdots+|a_k|$ and $\beta=|r_1|+\cdots+|r_k|$. Then if $w=a_1 r_1 a_1^{-1} \cdots a_k r_k a_k^{-1}$, then $|w|=2\alpha+\beta$. Let $\alpha_1$ [respectively $\beta_1$] be the number of letters of the $a_j^{\pm 1}$ [respectively of the $r_j$] which are cancelled and $\alpha_2$ [respectively $\beta_2$] that of the non-cancelled. Then by what seen above we have
 \begin{equation} \label{b&w} \alpha_1\leqslant \alpha_2+\beta_1. \end{equation}
Let $\alpha\geqslant \beta$, that is $\alpha=\beta+\gamma$ where $\gamma\geqslant 0$. Since $2\alpha=\alpha_1+\alpha_2$ then $\alpha_1+\alpha_2=2\beta_1+2\beta_2+2\gamma$ and $\alpha_1=-\alpha_2+2\beta_1+2\beta_2+2\gamma$, that is by (\ref{b&w}) we have
     $$-\alpha_2+2\beta_1+2\beta_2+2\gamma \leqslant \alpha_2+\beta_1$$
and thus
       $$\alpha_2\geqslant \gamma+\beta_2 + \frac{\beta_1}{2}=\frac{\beta+\beta_2}{2}+\gamma$$
     since $\beta=\beta_1+\beta_2$.
Then 
      $$|c|=\alpha_2+\beta_2\geqslant \frac{\beta+3\beta_2}{2}+\gamma.$$
  Let $n$ be a natural number and let $|c|\leqslant n$; then 
         $$\frac{\beta+3\beta_2}{2}+\gamma \leqslant n$$
      that is
          $$\gamma \leqslant n-\frac{\beta+3\beta_2}{2}.$$
  Since $\beta_2 \geqslant 1$ then     
          $$\gamma \leqslant n-\frac{\beta+3}{2}.$$
  Then 
             $$\alpha=\beta+\gamma \leqslant \beta+n-\frac{\beta+3}{2}=\frac{\beta-3}{2}+n.$$   
      Let $m$ be the maximal length of elements of $R$; then $\beta\leqslant mk$, that is    
          $$\alpha=\beta+\gamma \leqslant \frac{mk-3}{2}+n=\frac{mk-1}{2}+n-1.$$
  \end {remark}

      \begin{lemma} \label{unair} Let $\tau$ be an SLA in $\overline{R}$, let $c$ be the result of $\tau$ and suppose that $\eta(\tau)=\eta(c)$. Then if $\tau'$ is a pSLsA of $\tau$ and if $c'$ is the result of $\tau'$ then $\eta(\tau')=\eta(c')$.
    \end{lemma}

 \begin{proof} Let $n$ be the number of steps of $\tau$; we prove the claim by induction of $n$ being evident for $n=1$. Let $n>2$ and the claim be true for all the SLA's in $\overline{R}$ with less steps than $\tau$. Let the last step of $\tau$ be the cyclic conjugation of a word $c''$ and let $\tau''$ be the pSLsA of $\tau$ computing $c''$. If $\tau'$ is a pSLsA of $\tau$ then either $\tau'=\tau$ and the claim is obvious or $\tau'$ is a pSLsA of $\tau''$, in which case the claim is true by induction hypothesis. Let the last step of $\tau$ be the cyclic reduced product of two words $c_1$ and $c_2$ and let $\tau_1$ and $\tau_2$ be the pSLsA's of $\tau$ computing $c_1$ and $c_2$ respectively. If $\tau'=\tau$ the claim is obvious, otherwise $\tau'$ is a pSLsA of either $\tau_1$ or $\tau_2$ and the claim is true by induction hypothesis.
 \end{proof}

\begin{definition} \label{kg. corolla} \rm Let $\sigma$ be an SLA in $\overline{R}$ such that $\eta(\sigma)=k$ and let $c$ be the result of $\sigma$. Then we say that $c$ is a \textit{$k$-g. corolla}. We denote $\overline{R}_k$ the set of $k$-g. corollas. If $\sigma$ is an SLA in $C$ then we call $c$ a \textit{$k$-corolla}. We denote $C_k$ the set of $k$-corollas. \end{definition}       
    
Let $c$ be a $k$-g. corolla, that is there exists $\sigma$ whose result is $c$ and such that $\eta(\sigma)=k$. Let $c'$ be a cyclic conjugate of $c$ and let us call $\sigma'$ the SLA obtained by adding to $\sigma$ the cyclic conjugation from $c$ to $c'$. Then $\eta(\sigma')=k$ and the result of $\sigma'$ is $c'$, that is $c'$ is a $k$-g. corolla and thus $\overline{R}_k$ contains the cyclic conjugates of any of its elements.

\medskip

We have the following result, which we prove in the appendix.

\medskip

\textbf{Proposition} \ref{amex} \textit{Let $\sigma$ be an $\overline{R}$-SLA; then there exists a $C$-SLA $\sigma'$ with the same set of base elements and result of $\sigma$ and such that $\eta(\sigma')=\eta(\sigma)$. In particular $\overline{R}_k=C_k$.}

\medskip

\begin{remark} \label{bart} \rm We prove that $\overline{R}_k$ and $C_k$ contain the inverse of any of its elements, that is if $\sigma$ computes $c$ and $\eta(\sigma)=k$ then there exists $\tau$ which computes $c^{-1}$ and $\eta(\tau)=k$. We give the proof for $\overline{R}_k$, that for $C_k$ being analogous.

Let $\sigma$ be a $\overline{R}$-straight line algorithm computing $c$ and such that $\eta(\sigma)=k$ and let $n$ be the number of steps of $\sigma$. If $n=1$ then $c\in R$, thus $c^{-1}\in R$ and therefore there is an SLA with one step equal to $c^{-1}$. 

Let $n>1$ and the claim be true for every SLA with less steps than $\sigma$. Let the last step of $\sigma$ be a cyclic conjugation of a word $c_0$ and let $\sigma_0$ be the straight line subalgorithm computing $c_0$; then $\eta(\sigma_0)=\eta(\sigma)$=k. By induction hypothesis $c_0^{-1} \in \overline{R}_k$ and this means that there exists an SLA $\tau_0$ computing $c_0^{-1}$ and such that $\eta(\tau_0)=k$. Since $c^{-1}$ is a cyclic conjugate of $c_0^{-1}$ then if we add to $\tau_0$ the cyclic conjugation from $c_0^{-1}$ to $c^{-1}$ then we obtain an SLA which we denote $\tau$ that computes $c^{-1}$ and such that $\eta(\tau)=k$.

    Let the last step of $\sigma$ be a cyclically reduced product of two words $c_1$ and $c_2$ and let $\sigma_1$ and $\sigma_2$ be the proper straight line subalgorithms computing $c_1$ and $c_2$. We have that $k=\eta(\sigma)=\eta(\sigma_1)+\eta(\sigma_2)$. By induction hypothesis we have that there exist straight line algorithms $\tau_1$ and $\tau_2$ computing $c_1^{-1}$ and $c_2^{-1}$ and such that $\eta(\tau_1)+\eta(\tau_2)=k$. Let $\tau$ be the straight line algorithm obtained by adding to the steps of $\tau_1$ and $\tau_2$ the cyclically reduced product of $c_2^{-1}$ by $c_1^{-1}$. Then $\tau$ computes $c^{-1}$ and $\eta(\tau)=k$.
   \end{remark}

 \begin{lemma} \label{olgo} $\overline{R}_1$ is the set of cyclic conjugates of elements of $R$ and is equal to $C_1$. Let $k>1$; then $\overline{R}_k$ is the set of the cyclic conjugates of the cyclically reduced products of an element of $\overline{R}_m$ by one of $\overline{R}_n$, for $m$ and $n$ non-zero natural numbers such that $m+n=k$ and $m\leqslant n$; $C_k$ is the set of the cyclic conjugates of the cyclically reduced products of an element of $\overline{R}_{k-1}$ by one of $C_1$.
    \end{lemma}

  \begin{proof} For $k=1$ the claim is trivial. Let $k>1$, let $\sigma$ be an SLA in $\overline{R}$ such that $\eta(\tau)=k$ and let $c$ be the result of $\tau$. Let $s$ be the last step of $\tau$ which is a not a cyclic conjugation and let $c'$ be the output of $s$. Let $\tau'$ be the pSLsA of $\tau$ computing $c'$. There exist preceding outputs $c_1$ and $c_2$ such that $c'=\pi(c_1, c_2)$. Let $\tau_1$ and $\tau_2$ be the pSLsA's of $\tau$ computing $c_1$ and $c_2$ respectively. Let $m=\eta(\tau_1)$ and $n=\eta(\tau_2)$, that is $c_1 \in \overline{R}_m$ and $c_2 \in \overline{R}_n$. This means that $c$ is a cyclic conjugate of the cyclically reduced product of $c_1$ by $c_2$ and $m+n=k$. Finally by Proposition \ref{rpcrp} we have that $\pi(c_2, c_1)$ is a cyclic conjugate of $\pi(c_1, c_2)$, therefore we can assume that $m \leqslant n$.  
  
 Analogously is for $C_k$.   \end{proof}

\begin{theorem} \label{algo} If $R$ is finite (in particular $\langle \, X \, | \, R \, \rangle$ is a finite presentation) then $\overline{R}_k$ is finite for every natural number $k$. 
    \end{theorem}

    \begin{proof} We prove the claim by induction on $k$. Let $k=1$; by Lemma \ref{olgo}, $\overline{R}_1$ is equal to the set of cyclic conjugates of the elements of $R$ and is finite since $R$ is finite. Let $k>1$ and $\overline{R}_{k'}$ be finite for every $k'<k$; by Lemma \ref{olgo}, $\overline{R}_k$ is contained in the set of cyclic conjugates of the cyclically reduced products of an element of $\overline{R}_m$ by one of $\overline{R}_n$, for $m$ and $n$ non-zero natural numbers such that $m+n=k$. The claim follows from the fact that by induction hypothesis $\overline{R}_m$ and $\overline{R}_n$ are finite and that the number of cyclic conjugates of a given word is finite.
       \end{proof}

  \begin{remark} \label{Roverlinen} \rm Let $R$ be finite, let $k$ be a natural number and let us use Lemma \ref{olgo} and Remark \ref{puyi} to find an algorithm computing the elements of $C_k$. If $k=1$ then we compute the cyclic conjugates of the elements of $R$. Let $k>1$ and suppose to have computed $C_1, \cdots, C_{k-1}$; if $c \in C_k$ then there exist $c_1, c_2$ such that $c_1\neq c_2^{-1}$ and $c_1  \in C_{k-1}$, $c_2 \in C_1$ or $c_1  \in C_1$, $c_2 \in C_{k-1}$ and $c$ is a cyclic conjugate of $\pi(c_1, c_2)$. By Remark \ref{puyi} we have that either $c_1^{-1}$ is a prefix of $c_2$ or $c_2^{-1}$ is a suffix of $c_1$ or there exist cyclic conjugates $w_1$ of $c_1$ and $w_2$ of $c_2$ such that $\pi(c_1, c_2)=\rho(w_1 w_2)$. Since any $C_h$ is closed under cyclic conjugations, then to compute the elements of $C_k$ it is not necessary to compute all the $\pi(c_1, c_2)$ with $c_1  \in C_{k-1}$ and $c_2 \in C_1$ or $c_1  \in C_1$ and $c_2 \in C_{k-1}$, but just those for which there is cancellation between a suffix of $c_1$ and a prefix of $c_2$ but not between a prefix of $c_1$ and a suffix of $c_2$. 
 \end{remark}

\section{The set $L$ and the function $A$} \label{L}

We recall that we have called \textit{stem} (Definition \ref{stem}) a word of the form $w w ^{-1}$ where $w$ is reduced and $w\neq 1$, and that we have denoted $\mathcal{S}$ the set of all stems. The set $C$ has been introduced in Definition \ref{overlineR} (we recall that the elements of $C$ are called \textit{ corollas}) and insertions of words in Definition \ref{insert}.

\begin{definition} \label{maiobj} \rm We denote $L$ the set of results of straight line algorithms (Definition \ref{SLA}) whose universe set is $\mathcal{M}(X \cup X^{-1})$, whose base set is $B=\mathcal{S}\cup C$ and whose operations are the insertions of words. We will call \textit{$L$-straight line algorithms} or \textit{straight line algorithms in $L$} these straight line algorithms. 

We denote $L_g$ the set of results of straight line algorithms whose universe set is $\mathcal{M}(X \cup X^{-1})$, whose base set is $B_g=\mathcal{S}\cup \overline{R}$ and whose operations are the insertions of words. We will call \textit{$L_g$-straight line algorithms} or \textit{straight line algorithms in $L_g$} these straight line algorithms.  \end {definition}

We observe that $\overline{R} \cap \mathcal{S}=\emptyset$ because every element of $\overline{R}$ is reduced and every one of $\mathcal{S}$ is not. Obviously $L_g\supset L$ since $\overline{R}\supset C$. We recall from that an element of $B_g$ is called a \textit{base element}.

Conjugating a word $v$ with a reduced word $w$ is equivalent to inserting $v$ into the stem $w w ^{-1}$ at $w$. By Definition \ref{reducprod}, the product of two reduced words is equal to an insertion of the cancelled part (which is a stem) into their reduced product.

$L_g$ is closed under product (which is a special case of insertion) and under conjugation with a reduced word, because if $l \in L_g$ and if $w$ is reduced then $w w^{-1}$ is a stem, therefore belongs to $L_g$ and $w l w ^{-1}$ is an insertion of $l$ into $w\, w ^{-1}$. This means that $L_g$ contains any “non-cancelled" product of conjugates of elements of $C$. This implies

 \begin{proposition} \label{inv} Let $\mathcal{N}$ be the normal closure of $R$ in $\mathcal{F}(X)$, i.e., the set of relators of the presentation $\langle \, X \, | \, R \, \rangle$. If $w\in \mathcal{N}$ then there exists an element of $L_g$ whose reduced form is $w$. \end{proposition}  
 
 \begin{proof} Follows from what said above because every element of $\mathcal{N}$ is the reduced form of a product of conjugates of elements of $R\subset C$.     \end{proof}       
   
The following result is a converse of Proposition \ref{inv}.

 \begin{proposition} \label{veta} The reduced form of any element of $L_g$ is a relator, that is it belongs to $\mathcal{N}$. \end{proposition}

\begin{proof} We have to prove that $\rho(L) \subset \mathcal{N}$ (recall that $\rho$ denotes the reduced form, Definition \ref{rho}). Let $\sigma$ be an SLA in $L_g$, let $w$ be the result of $\sigma$ and let $n$ be the number of steps of $\sigma$; we prove the claim by induction on $n$.

Let $n=1$; then $w$ belongs either to $C$ (in which case $\rho(w)=w \in C\subset \mathcal{N}$) or $w \in \mathcal{S}$ (in which case $\rho(w)=1 \in \mathcal{N}$). 

Let $n>1$ and the claim be true for every SLA with less steps than $\sigma$. Then $w$ is of the form $u w' v$ and the last step is the insertion of $w'$ into $uv$ at $u$. By induction hypothesis, $\rho(uv)$ and $\rho(w')$ belong to $\mathcal{N}$; since $\rho(w)=\rho(uw'u^{-1} uv)$ then $\rho(w)$ belongs to $\mathcal{N}$.                                        \end{proof}

We recall that the \textit{g. corollas} are the elements of $C$. The function $\eta(c)$ for a g. corolla $c$ has been introduced in Definition \ref{eta(c)}.

 \begin{definition} \label{A} \rm  Let $\sigma$ be an SLA in $L_g$; we set $A(\sigma):=\sum \eta(c)$ where $c$ varies in the set of g. corollas of $\sigma$. Let $w \in L_g$; we set $A(w):=\min\{ A(\sigma) : \sigma \textrm{ is an } SLA \textrm{ in } L \textrm{ computing } w\}$.
         \end{definition}

\begin{remark} \label{etaA} \rm Let $c$ be a g. corolla and let $\sigma$ be the SLA in $L_g$ with a single step equal to $c$. Then $\sigma$ computes $c$ and $A(\sigma)=\eta(c)$, therefore $A(c)\leqslant \eta(c)$.
     \end{remark}

 For the area and the work of a relator see Definitions \ref{defarea} and \ref{omega} respectively.

\begin{definition} \label{propcor} \rm We call \textit{proper (g.) corolla} a (g.) corolla $c$ such that $A(c)=\eta(c)$. For every natural number $n$ set 
    $$\Delta_0(n):=\max\{\textrm{Area}(w) : \textrm{$w$ is a proper corolla} \,\,\textrm{and} \, \, |w|\leqslant n\},$$
    $$\Omega_0(n):=\max\{\textrm{Work}(w) : \textrm{$w$ is a proper corolla} \,\,\textrm{and} \, \, |w|\leqslant n\}.$$    \end{definition}   
    
 Obviously $\Delta_0\leqslant \Delta$ and  $\Omega_0\leqslant \Omega$ where $\Delta$ is the Dehn function (Definition \ref{Dehn}) and $\Omega$ has been introduced in Definition \ref{omega}.

\section{Statement of the results} \label{statres}

\begin{definition} \label{CMDR} \rm Let $\sigma$ be a straight line algorithm in $L_g$, let $c_1, \cdots, c_m$ be the corollas of $\sigma$ (counted with their multiplicity) and let $\tau_1, \cdots, \tau_m$ be straight line algorithms in $\overline{R}$ computing them and such that $\eta(\tau_i)=\eta(c_i)$ for $i=1, \cdots, m$. The union\footnote{for union of multisets see Definition \ref{unmult}} of the multisets of base elements (Definition \ref{basele}) of $\tau_1, \cdots, \tau_m$ is called a \textit{complete multiset of defining relators} (abbreviated CMDR) \textit{for $\sigma$}. 
      \end{definition}

\begin{remark} \label{union} \rm Let $\sigma$ be a straight line algorithm in $L_g$ and let $M$ be a CMDR for $\sigma$. Then it is obvious that $A(\sigma)=|M|$.

Let $\sigma, \sigma_1$ and $\sigma_2$ be straight line algorithms in $L_g$ such that the multiset of corollas of $\sigma$ is the union of those of $\sigma_1$ and of $\sigma_2$ and let $M_1$ and $M_2$ be CMDR for $\sigma_1$ and $\sigma_2$ respectively. Then $A(\sigma)=A(\sigma_1)+A(\sigma_2)$ and $M_1\cup M_2$ is a CMDR for $\sigma$. 

  In particular this is the case when the last step of an SLA $\sigma$ is the insertion of a word $w_2$ into a word $w_1$ and $\sigma_1$ and $\sigma_2$ are the pSLsA's (Remark \ref{pSLsA}) of $\sigma$ computing $w_1$ and $w_2$.
  \end{remark}

  \begin{proposition} \label{impo} Let $\sigma$ be an SLA in $L_g$ and let $w$ be its result. Then $Area\big(\rho(w)\big)\leqslant A(\sigma)$; in particular $\textrm{Area}\big(\rho(w)\big) \leqslant A(w)$.
      \end{proposition}
         
\begin{proof} We prove the claim by induction on the number of steps of $\sigma$. If $\sigma$ has only one step then $w$ is a stem or a corolla. If $w$ is a stem then $\rho(w)=1$, $\textrm{Area}\big(\rho(w)\big)=0$ and $A(\sigma)=0$. If $w$ is a corolla (in particular it is reduced, i.e., $\rho(w)=w$) then $A(\sigma)=\eta(w)$ and the claim follows from Corollary \ref{areaeta}. 

Let $\sigma$ have more than one step and the claim be true for every SLA with less steps than $\sigma$. The last step of $\sigma$ is the insertion of a word $w_2$ into a word $w_1$, that is there exist words $u, v$ such that $w_1=uv$ and $w=u w_2 v$. Let $\sigma_1$ and $\sigma_2$ be the pSLsA's computing $w_1$ and $w_2$; by induction hypothesis $\textrm{Area}\big(\rho(w_i)\big)\leqslant A(\sigma_i)$ for $i=1,2$. Since $A(\sigma)=A(\sigma_1)+A(\sigma_2)$ by Remark \ref{union}, then 
    \begin{equation} \label{eq1} \textrm{Area}\big(\rho(w_1)\big) + \textrm{Area}\big(\rho(w_2)\big)  \leqslant A(\sigma).\end{equation}
 Since $\rho(w)=\rho(u w_2 u^{-1} u v)$ and $w_1=uv$, then 
    $$\rho(w)=\rho(u w_2 u^{-1}) \rho(w_1);$$
by (\ref{Area(vw)}) we have that 
    \begin{equation} \label{eq2} \textrm{Area}\big(\rho(w)\big) \leqslant \textrm{Area}\big(\rho(u w_2 u^{-1})\big)+\textrm{Area}\big(\rho(w_1)\big).\end{equation}
  Finally the claim follows from (\ref{eq2}) by virtue of (\ref{eq1}) and of the fact that $\textrm{Area}\big(\rho(u w_2 u^{-1})\big)=\textrm{Area}\big(\rho(w_2)\big)$ by (\ref{uwu-1}) of Section \ref{intro}.
         \end{proof}

\begin{lemma} \label{quellemma} Let $w\in L_g$ and let $\sigma$ be an SLA computing $w$ and such that $\textrm{Area}\big(\rho(w)\big)=A(\sigma)$. Then any corolla of $\sigma$ is a proper corolla (Definition \ref{propcor}).
      \end{lemma}

\begin{proof} We prove the claim by induction on the number of steps of $\sigma$. If $\sigma$ has one step then $w$ is the only corolla of $\sigma$. Moreover it is reduced, i.e., $\rho(w)=w$ and $A(\sigma)=\eta(w)$. By Definition \ref{A}, $A(w)\leqslant A(\sigma)$ and $A(\sigma)=\textrm{Area}\big(\rho(w)\big)$ by hypothesis. Since $\textrm{Area}\big(\rho(w)\big) \leqslant A(w)$ by Proposition \ref{impo}, then 
  $$A(w)\leqslant A(\sigma)=\eta(w)=\textrm{Area}\big(\rho(w)\big)  \leqslant A(w),$$
  therefore $A(w)=\eta(w)$ and $w$ is a proper corolla.

 Let $\sigma$ have more than one step and the claim be true for every SLA with less steps than $\sigma$. The last step of $\sigma$ is the insertion of a word $w_2$ into a word $w_1$, that is there exist words $u, v$ such that $w_1=uv$ and $w=u w_2 v$. As in the proof of Proposition \ref{impo}, $\rho(w)=\rho(u w_2 u^{-1} w_1)$ and  
   \begin{equation} \label{rho(w)} \textrm{Area}\big(\rho(w)\big) \leqslant \textrm{Area}\big(\rho(w_2)\big)+\textrm{Area}\big(\rho(w_1)\big).\end{equation}
 Let $\sigma_1$ and $\sigma_2$ be the pSLsA's computing $w_1$ and $w_2$; we have that $A(\sigma)=A(\sigma_1)+A(\sigma_2)$ by Remark \ref{union}. If $\textrm{Area}\big(\rho(w_1)\big)<A(\sigma_1)$ or $\textrm{Area}\big(\rho(w_2)\big)<A(\sigma_2)$, then by (\ref{rho(w)})
          $$\textrm{Area}\big(\rho(w)\big) < A(\sigma_1)+A(\sigma_2)=A(\sigma)$$
which is contrary to the hypothesis. Thus $A(\sigma_1)=\textrm{Area}\big(\rho(w_1)\big)$ and $A(\sigma_2)=\textrm{Area}\big(\rho(w_2)\big)$, therefore by induction hypothesis for every corolla $c$ of $\sigma_1$ or of $\sigma_2$ we have that $A(c)=\eta(c)$. The claim follows from the fact that every corolla of $\sigma$ is a corolla of $\sigma_1$ or of $\sigma_2$.
         \end{proof}

We can now state the Main Theorem of this thesis; we recall that $\overline{R}$ (the set of corollas) has been introduced in Definition \ref{overlineR}, $L$ in Definition \ref{maiobj}, the function $A$ in Definition \ref{A} and complete multiset of defining relators in Definition \ref{CMDR}.

\medskip

\textbf{Main Theorem \ref{maithe}} \textit{Let $\langle \, X \, | \, R \, \rangle$ be a group presentation and let $\mathcal{N}$ be the set of reduced relators. Then $\mathcal{N}$ coincides with the subset of $L$ consisting of reduced words. Let $w$ be the reduced form of $f_1 r_1 f_1^{-1} \cdots f_n r_n f_n^{-1}$, where $r_i\in R$; then there exist a submultiset $M$ of $\{r_1, \cdots, r_n\}$ and a straight line algorithm $\sigma$ computing $w$ which has $M$ as a complete multiset of defining relators (CMDR) and such that $A(\sigma)\leqslant n$. If $n=\textrm{Area}(w)$ then $A(\sigma)=n$, every corolla of $\sigma$ is a proper corolla and $\{r_1, \cdots, r_n\}$ is a CMDR for $\sigma$. Finally $A(w)=\textrm{Area}(w)$.}

\medskip

The Main Theorem implies the following 

\medskip

\textbf{Corollary \ref{maicor}} \textit{Let $w \in \mathcal{N}$ and let $\sigma$ be an SLA computing $w$ and such that $A(\sigma)=A(w)$. Then the area of $w$ is equal to the sum of the areas of the corollas of $\sigma$, that is if $c_1, \cdots, c_m$ are the corollas of $\sigma$ then $\textrm{Area}(w)=\displaystyle \sum_{i=1}^m \textrm{Area}(c_i)$ and $|w|\geqslant \sum_{i=1}^m |c_i|$.} 

\medskip

Let $\langle \, X \, | \, R \, \rangle$ be a group presentation and suppose that there exists a positive real constant $\alpha$ such that $\textrm{Area}(c) \leqslant \alpha |c|$ for every proper corolla $c$. Let $w$ be a relator; by virtue of Corollary \ref{maicor} there exist proper corollas $c_1, \cdots, c_m$ such that $\textrm{Area}(w)=\displaystyle \sum_{i=1}^m \textrm{Area}(c_i)$ and $|w|\geqslant \sum_{i=1}^m |c_i|$. Therefore $$\textrm{Area}(w)=\displaystyle \sum_{i=1}^m \textrm{Area}(c_i)\leqslant \sum_{i=1}^m \alpha |c_i|\leqslant \alpha |w|$$
   and $\langle \, X \, | \, R \, \rangle$ is hyperbolic (Definition \ref{hyper}). Since the converse is obvious, this proves the following

\medskip

\textbf{Corollary \ref{hyper2}}\textit{ The presentation $\langle \, X \, | \, R \, \rangle$ is hyperbolic if and only if there exists a positive real constant $\alpha$ such that $\textrm{Area}(c) \leqslant \alpha |c|$ for every proper corolla $c$.}

\medskip

Proposition \ref{vega} is stronger than Corollary \ref{hyper2}, anyway the proof of Corollary \ref{hyper2} is direct and does not use Proposition \ref{vega}.

For every natural number $n$ we have set $\Delta_0(n)$ as the maximal area of proper corollas of length at most $n$ (Definition \ref{propcor}); obviously $\Delta_0(n) \leqslant \Delta(n)$. By the preceding calculation we have

\medskip

 \textbf{Corollary \ref{hyper3}}\textit{ Let $\alpha$ be a positive real number; then the Dehn function $\Delta$ is bounded by the linear function $\alpha n$ if and only if the $\Delta_0$ is. }
 
 \medskip

To prove the Main Theorem it is sufficient to prove

\medskip

\textbf{Lemma \ref{maile}} \textit{ Let $l:=l_1 z z^{-1} l_2$ (where $l_1$ and $l_2$ are words and $z$ a letter) be an element of $L_g$; then $l_1 l_2\in L_g$. In particular if $\sigma$ is a straight line algorithm computing $l$ and if $M$ is a CMDR for $\sigma$ (Definition \ref{CMDR}), then there exists a straight line algorithm $\sigma'$ computing $l_1 l_2$, such that $A(\sigma')\leqslant A(\sigma)$ and such that if $M'$ is a CMDR for $\sigma'$ then $M'\subset M$. Moreover if $A(\sigma)=\textrm{Area}\big(\rho(l)\big)$ then $A(\sigma')=A(\sigma)$ and $M'=M$.}

\medskip

Let Lemma \ref{maile} be true; take $w \in \mathcal{N}$ and let 
        $$w=\rho(f_1 r_1 f_1^{-1} \cdots f_n r_n f_n^{-1})$$
where $f_i \in \mathcal{F}(X)$ and $r_i \in R$ for $i=1, \cdots, n$. Since $R\subset \overline{R}$, then the $r_i$ are corollas. Let $\sigma$ be the SLA consisting in the insertions of $r_i$ into the stem $f_i f_i^{-1}$ at $f_i$ which give $g_i:=f_i r_i f_i^{-1}$ and then in the products $g_1 g_2$, $g_1 g_2 g_3$, $\cdots$, $g_1 g_2\cdots g_n$. Its result is $f_1 r_1 f_1^{-1} \cdots f_n r_n f_n^{-1}$, $A(\sigma)=n$ and $\{r_1, \cdots, r_n\}$ is a complete multiset of defining relators for $\sigma$. $w$ is the reduced form of $f_1 r_1 f_1^{-1} \cdots f_n r_n f_n^{-1}$ and is obtained from the latter by performing all the possible cancellations. By applying repeatedly Lemma \ref{maile}, we obtain an $L_g$-SLA $\sigma'$ whose result is $w$, such that $A(\sigma')\leqslant A(\sigma)$ and such that if $M'$ is a complete multiset of defining relators for $\sigma'$ then $M'$ is contained in $\{r_1, \cdots, r_n\}$. By Proposition \ref{amex} it follows that $L=L_g$.

Let $n=\textrm{Area}(w)$; then $n\leqslant A(\sigma')$ by Proposition \ref{impo} and since $A(\sigma')\leqslant A(\sigma)=n$ then $A(\sigma')=n$. Thus we have that $A(w)\leqslant A(\sigma')=\textrm{Area}(w)$ and since $\textrm{Area}(w)\leqslant A(w)$ by Proposition \ref{impo}, this implies that $A(w)=\textrm{Area}(w)$. By Remark \ref{union} we have that $|M'|=n$ therefore $M'=\{r_1, \cdots, r_n\}$ and $\{r_1, \cdots, r_n\}$ is a CMDR for $\sigma'$. Finally Lemma \ref{quellemma} implies that every corolla of $\sigma'$ is a proper corolla.

Lemma \ref{maile} also implies  

\medskip

\textbf{Lemma \ref{maile2}} \textit{1. Let $l \in L$ and let $w \in \mathcal{M}(X\cup X^{-1})$ be such that $\rho(w)=\rho(l)$. Then $w \in L$.}

\medskip
  
  \textit{2. $L$ contains the cyclic conjugate of any of its elements.}
  
\medskip
  
\begin{proof} 

\smallskip  

1. $w$ is obtained from $l$ by insertions and deletions of words of the form $zz^{-1}$ where $z\in X\cup X^{-1}$. Since $zz^{-1}$ is a stem, if we insert $zz^{-1}$ into an element of $L$ we obtain an element of $L$; Lemma \ref{maile} says that if we delete a subword of the form $zz^{-1}$ from an element of $L$ we still have an element of $L$.

\medskip

2. Let $l \in L$ and let $l'$ be a cyclic conjugate of $l$; then there exist words $l_1$ and $l_2$ such that $l=l_1 l_2$ and $l'=l_2 l_1$. The word $l_1^{-1} l_1 l_2 l_1$ belongs to $L$ because it is an insertion of $l$ into the stem $l_1^{-1} l_1$. By applying repeatedly Lemma \ref{maile} to $l_1^{-1} l_1 l_2 l_1$ we have that $l_2 l_1$ belongs to $L$.

            \end{proof}

\medskip

Lemma \ref{maile2} implies the following interesting result,

\medskip

\textbf{Theorem \ref{maithe2}} \textit{ $L$ is the subset of $\mathcal{M}(X\cup X^{-1})$ of words whose reduced form belongs to $\mathcal{N}$, i.e., if $\rho:\mathcal{M}(X\cup X^{-1}) \rightarrow \mathcal{F}(X)$ is the function “reduced form" (Definition \ref{rho}) then $L=\rho^{-1}(\mathcal{N})$.} 

\medskip

\begin{proof}
    By Proposition \ref{veta} we have that $L\subset \rho^{-1}(\mathcal{N})$; to prove Theorem \ref{maithe2} it is sufficient to prove the reverse inclusion. Let $w\in \mathcal{M}(X\cup X^{-1})$ be such that $\rho(w) \in \mathcal{N}$; by Proposition \ref{inv} there exists $l \in L$ such that $\rho(l)=\rho(w)$ and by Part 1 of Lemma \ref{maile2}, $w \in L$.
 \end{proof}

\chapter{Preliminary results}  \label{pr}

In Section \ref{ggior} we show how to associate a 2-cell complex with a straight line algorithm in $L$. Sections \ref{prel1} and \ref{prel2} introduce useful tools like \textit{stem elements}, \textit{flowers}, \textit{flower elements}, \textit{straight line subalgorithms}, \textit{ramifications} and \textit{surrounds}. These notions, although being technical, have a clear and intuitive interpretation as their names suggest.

\section{2-cell complexes for straight line algorithms}  \label{ggior}

In this section we show how to associate a 2-cell complex with any $L_g$-straight line algorithm. For definitions and properties of 2-cell complexes see III of \cite{LS}. In this thesis we will consider only 2-complexes equipped with a cycle $\gamma$ such that for every edge $e$ of the complex either $e$ or $e^{-1}$ are contained in $\gamma$. We call $\gamma$ and its initial vertex respectively \textit{the boundary cycle} and \textit{the initial vertex} of the 2-complex. We will suppose also that these complexes are connected and that every edge is labeled by an element of $X\cup X^{-1}$.

\begin{definition} \rm \label{orientation} We say that the \textit{orientation} of a 2-complex is \textit{compatible} with the \textit{orientation} of a face it contains if the boundary cycle of the face is a (non-necessarily contiguous) subpath of the boundary cycle of the 2-complex. \end{definition} 

This means in particular that the edges of the boundary of the face are in the same order in the boundary of the 2-complex.

\begin{definition} \label{extremal} \rm Given two faces $F$ and $F'$ of a complex, we say that $F$ is \textit{comprised} in $F'$ if all the edges of the boundary cycle of $F$ are comprised between two consecutive edges of the boundary cycle of $F'$. The latter is a transitive and antisymmetric relation in the set of faces of a given complex. Consider the reflexive closure \cite{reclo} of this relation; since it is a finite partial order, by Zorn's lemma \cite{zorn} there are minimal elements. We call such minimal elements \textit{extremal faces}. \end{definition}

From now on, unless otherwise specified, with the term \textit{straight line algorithm} (or SLA) without other specifications we mean \textit{straight line algorithm in $L_g$}. 

Let $\sigma$ be an SLA and let $s$ be one of its steps; three cases are possible: 

\textbf{1.} $s$ is a corolla; \textbf{2.} $s$ is a stem; \textbf{3.} there exist two steps $s_1$ and $s_2$ preceding $s$ such that $s$ is an insertion of $s_2$ into $s_1$ (we recall that we identify a step with its output).

\smallskip

\textbf{First case}: We say that \textit{$s$ is a corolla of $\sigma$}. If $s$ is the empty word then we associate with $s$ a graph consisting of a single vertex. If $s=x_1 \cdots x_m$ then we associate with $s$ the following contractible 2-cell complex

 
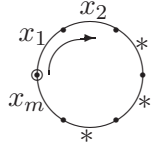
\begin{picture}(180,100)(12,-5) 

\put(177,40){\circle{40}}

\put(157,40){\circle*{2}}
\put(167,57){\circle*{2}}
\put(187,57){\circle*{2}}
\put(197,40){\circle*{2}}
\put(187,23){\circle*{2}}
\put(167,23){\circle*{2}}

\put(175,54){\vector(1,0){4}}
\put(175,40){\oval(27,27)[tl]}

\put(157,40){\circle{4}}

\put(150,53){$x_1$}
\put(173,63){$x_2$}
\put(194,46){*}
\put(195,24){*}
\put(173,11){*}
\put(146,28){$x_m$}

\end{picture}

whose edges are labeled consecutively by $x_1, x_2, \cdots, x_m$. The vertex with a circle surrounding it, is the initial vertex. The boundary, whose orientation is determined by the arrow inside it, is a simple cycle.

\smallskip

\textbf{Second case}: We say that \textit{$s$ is a stem of $\sigma$}. If $s=x_1 \cdots x_m \,  x_m^{-1} \cdots  x_1^{-1}$ then we associate with $s$ the following contractible 2-cell complex

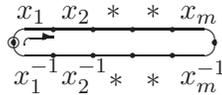
\begin{picture}(180,85)(6,-3)

\put(173,40){\oval(76,10)}

\put(147,39){\oval(16,5)[tl]}
\put(146,42){\vector(1,0){4}}

\put(135,40){\circle{4}}

\put(135,40){\circle*{2}}
\put(150,45){\circle*{2}}
\put(165,45){\circle*{2}}
\put(180,45){\circle*{2}}
\put(195,45){\circle*{2}}
\put(211,40){\circle*{2}}

\put(195,35){\circle*{2}}
\put(180,35){\circle*{2}}
\put(165,35){\circle*{2}}
\put(150,35){\circle*{2}}

\put(136,49){$x_1$}
\put(153,49){$x_2$}
\put(170,46){*}
\put(185,46){*}
\put(198,49){$x_m$}
\put(198,24){$x_m^{-1}$}

\put(185,20){*}
\put(171,20){*}
\put(153,24){$x_2^{-1}$}
\put(135,24){$x_1^{-1}$}

\end{picture}

whose edges are labeled consecutively by $x_1$, $\cdots$, $x_m$, $x_m^{-1}$, $\cdots$, $x_1^{-1}$. The boundary is a simple cycle. The second vertex of the edge labeled by $x_m$ (which coincides with the first vertex of that labeled by $x_m^{-1}$) is called \textit{vertex in the middle} or \textit{mid-vertex} of the stem.

\smallskip

We use the terms \textit{corolla} and \textit{stem} also for the associated complexes, that is we call \textit{corolla} (or \textit{stem}) a 2-cell whose boundary is a simple cycle labeled by a word of $\overline{R}$ (or of $\mathcal{S}$). This will not cause ambiguity.

\smallskip

\textbf{Third case}: Let $s_1:=x_1 \cdots x_m$, $s_2:=y_1 \cdots y_p$. Then there exists $n: 1 \leqslant n \leqslant m$ such that 
    $$s=x_1 \cdots x_n \, y_1 \cdots y_p \, x_{n+1} \cdots x_m.$$

 We say that \textit{$s$ is the insertion of $s_2$ into $s_1$ at the $n$-th letter or at $x_n$ or at the subword $x_1 \cdots x_n$}. Let $i: 1 \leqslant i \leqslant n$; we say that \textit{the $i$-th letter of $s$ comes directly from (the $i$-th letter of) $s_1$}. Let $i: n+1 \leqslant i \leqslant n+p$; we say that \textit{the i-th letter of $s$ comes directly from (the $(i-n)$-th letter of) $s_2$}. Let $i: p+n+1 \leqslant i \leqslant p+m$; we say that \textit{the i-th letter of $s$ comes directly from (the $(i-p)$-th letter of) $s_1$}.

Let $\mathcal{C}_1$ and $\mathcal{C}_2$ be the complexes associated with $s_1$ and $s_2$ respectively and let $v_2$ be the initial vertex of $\mathcal{C}_2$. Let $v_1$ be the final vertex of the edge of $\mathcal{C}_1$ labeled by $x_n$; we associate with $s$ the complex obtained by grafting $\mathcal{C}_2$ into $\mathcal{C}_1$ at $v_2$ and $v_1$ (see \cite{adjunc}), that is the complex $\mathcal{C}_1\cup_f \mathcal{C}_2$ where $f$ is the function from $\{v_2\}$ to $\mathcal{C}_1$ such that $f(v_2)=v_1$. This complex is obtained by joining $\mathcal{C}_1$ to $\mathcal{C}_2$ in such a way that $v_2$ coincides with $v_1$. The intersection of $\mathcal{C}_1$ and $\mathcal{C}_2$ is a single vertex and $\mathcal{C}_1\cup_f \mathcal{C}_2$ is their union. The initial vertex of $\mathcal{C}_1\cup_f \mathcal{C}_2$ is the initial vertex of $\mathcal{C}_1$, its boundary cycle is the cycle obtained by inserting the boundary cycle of $\mathcal{C}_2$ into that of $\mathcal{C}_1$ between the edges labeled by $x_n$ and $x_{n+1}$.

Suppose for instance that we have associated the following labeled complex with $s_1$ 

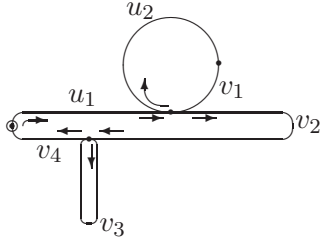
\begin{picture}(180,147)(18,-29)

\put(129,41){\oval(10,2)[tl]}    

\put(129,42){\vector(1,0){4}}

 \put(168,43){\vector(1,0){9}}

\put(179,54){\oval(18,13)[bl]}   

\put(170,54){\vector(0,1){4}}

 \put(188,43){\vector(1,0){9}}

 \put(146,38){\vector(-1,0){9}}

 \put(162,38){\vector(-1,0){9}}

 \put(150,33){\vector(0,-1){9}}

\put(141,49){$u_1$}
 
\put(163,82){$u_2$}

\put(198,52){$v_1$}

\put(227,39){$v_2$}

\put(153,0){$v_3$}

\put(129,28){$v_4$}

\put(120,40){\circle{4}}

\put(120,40){\circle*{2}}    

\put(180,45){\circle*{2}}    

\put(198,64){\circle*{2}}   

\put(149,35){\circle*{2}}

\put(173,40){\oval(106,10)}   

\put(180,63){\circle{35}}        

\put(149,19){\oval(6,32)}

\end{picture}

where $u_1, u_2, v_1, \cdots, v_4$ are words such that $x_1\cdots x_n=u_1 u_2$ and $x_{n+1}\cdots x_m= v_1 v_2 v_3 v_4$; and the following one to $s_2$

 
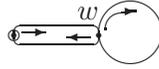
\begin{picture}(180,80)(40,5)

 \put(174,41){\vector(1,0){9}}

 \put(200,39){\vector(-1,0){9}}

\put(213,42){\oval(15,13)[tl]}    

\put(213,49){\vector(1,0){4}}

\put(187,40){\oval(32,7)}

\put(171,40){\circle{4}}

\put(171,40){\circle*{2}}

\put(203,40){\circle*{2}}

\put(215,41){\circle{25}}

\put(195,46){$w$}

\end{picture}

where $w=y_1 \cdots y_p$; then we associate with $s$ the following complex

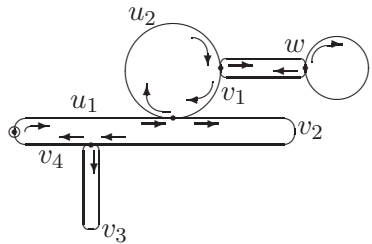
\begin{picture}(180,140)(24,-25)

\put(129,41){\oval(10,2)[tl]}    

\put(129,42){\vector(1,0){4}}

 \put(168,43){\vector(1,0){9}}

\put(179,54){\oval(18,13)[bl]}    

\put(170,54){\vector(0,1){4}}

\put(187,69){\oval(12,13)[tr]}   

\put(193,69){\vector(0,-1){4}}

\put(189,59){\oval(12,13)[br]}   

\put(189,52){\vector(-1,0){4}}

 \put(188,43){\vector(1,0){9}}

 \put(162,38){\vector(-1,0){9}}

 \put(146,38){\vector(-1,0){9}}

 \put(150,33){\vector(0,-1){9}}

\put(141,49){$u_1$}
 
\put(163,82){$u_2$}

\put(198,52){$v_1$}

\put(227,39){$v_2$}

\put(153,0){$v_3$}

\put(129,28){$v_4$}

\put(120,40){\circle{4}}

\put(120,40){\circle*{2}}    

\put(180,45){\circle*{2}}    

\put(198,64){\circle*{2}}    

\put(149,35){\circle*{2}}

\put(173,40){\oval(106,10)}    

\put(180,63){\circle{35}}        

\put(149,19){\oval(6,32)}

 \put(201,65){\vector(1,0){9}}

 \put(227,63){\vector(-1,0){9}}

\put(240,66){\oval(15,13)[tl]}    

\put(240,73){\vector(1,0){4}}

\put(214,64){\oval(32,7)}

\put(230,64){\circle*{2}}

\put(242,64){\circle{25}}

\put(222,70){$w$}

\end{picture}

whose label is $u_1 u_2 \, w \, v_1 v_2 v_3 v_4=x_1 \cdots x_n \, s_2 \, x_{n+1} \cdots x_m$. We say that the complex of $s_2$ has been \textit{grafted into the complex of $s_1$ at the $n$-th vertex}. 

We associate with a straight line algorithm $\sigma$ the complex that has been associated with its last step. 

\smallskip   \smallskip

If a step $s$ depends directly on a step $s'$, then the output of $s'$ is a subword of the output of $s$. More generally this is true also if $s$ depends on $s'$.

\begin{definition} \label{coming from} \rm We have defined in the Third case when a letter of a step comes directly from one of another step. We now define in the set of the letters of (the outputs of) the steps of $\sigma$, the reflexive transitive closure \cite{reclo} of the relation “coming directly from" and we call it \textit{relation of coming from}. If a letter $x$ is in that relation with another one $x'$, then we say that \textit{$x$ comes from $x'$}. Then $x$ comes from $x'$ if either they are the same letter of the same step or if there is a finite sequence of letters starting with $x$, ending with $x'$ and such that every letter of the sequence comes directly from the previous one.

If a letter $x$ of a step $s$ comes from a letter $x'$ of $x'$, then we can say improperly that \textit{$s$ contains $c'$}. 
         \end{definition}

\begin{remark} \label{s'} \rm 
If a step $s$ contains a letter of a step $s'$, then this means that $s$ depends on $s'$; therefore $s'$ is a subword of $s$, that is $s$ contains all the letters of $s'$.
          \end{remark}

 \begin{definition} \label{mathcalL} \rm We call $\mathcal{L}$ the set of 2-complexes associated with straight line algorithms in the way shown in the three cases above.
 \end{definition}
 
  Let $\mathcal{M}$ be the set of connected 2-cell complexes whose edges are labeled by elements of $X\cup X^{-1}$ and let $\mathcal{B}$ be the set of stems and corollas. $\mathcal{L}$ is the set of results of SLA's whose universe set is $\mathcal{M}$, whose base set is $\mathcal{B}$ and whose operation is the \textit{grafting} of complexes as seen in the Third Case. We call SLA's in $\mathcal{L}$ these SLA's.

  We have the following result.

 \begin{theorem} \label{th1} $\mathcal{L}$ is the set of labeled connected 2-cell complexes

\begin{itemize}

  \item whose faces are stems or corollas,
  
  \item whose orientation is compatible with the orientation of its faces (see Definition \ref{orientation}),

  \item in which given an edge $e$, either $e$ or $e^{-1}$ belong to the boundary of some face.
  
\end{itemize}
 
            \end{theorem}

\begin{proof} Let $\sigma$ be an SLA in $\mathcal{L}$, let $\mathcal{C}$ be the result of $\sigma$ and let $n$ be the number of steps of $\sigma$. We prove by induction on $n$ that $\mathcal{C}$ has the three properties stated in the thesis of the proposition, being this evident for $n=1$.

   Let $n>1$ and the claim be true for all the SLA's in $\mathcal{L}$ with less steps than $\sigma$. The last step of $\sigma$ is a grafting of a complex $\mathcal{C}_1$ into a complex $\mathcal{C}_2$ at vertices $v_1$ and $v_2$ of $\mathcal{C}_1$ and $\mathcal{C}_2$ respectively. By induction hypothesis  $\mathcal{C}_1$ and $\mathcal{C}_2$ verify the properties of the claim. $\mathcal{C}$ is connected because it is the non-disjoint union of two connected complexes; its faces are the faces of $\mathcal{C}_1$ and those of $\mathcal{C}_2$, therefore they are stems and corollas. Let $\gamma'_1$ and $\gamma''_1$ be paths such that $\gamma_1=\gamma'_1 \gamma''_1$ and the final vertex of $\gamma'_1$ is $v_1$ and let $\gamma=\gamma'_1 \gamma_2 \gamma''_1$; $\gamma_1$ and $\gamma_2$  are subpaths of $\gamma$. A face of $\mathcal{C}$ is a face of $\mathcal{C}_1$ or of $\mathcal{C}_2$, thus its boundary cycle is a subpath of $\gamma_1$ or $\gamma_2$ and therefore of $\gamma$, thus the orientation of $\mathcal{C}$ is compatible with that of its faces. Finally every edge $e$ of $\mathcal{C}$ is an edge of $\mathcal{C}_1$ or of $\mathcal{C}_2$, therefore either $e$ or $e^{-1}$ are contained in the boundary of a face of $\mathcal{C}_1$ or $\mathcal{C}_2$ and thus in the boundary of a face of $\mathcal{C}$.

We now prove by induction on the number of faces that any complex $\mathcal{C}$ verifying the three properties of the theorem belongs to $\mathcal{L}$. If $\mathcal{C}$ has only one face then it is a stem or a corolla and thus belongs to $\mathcal{L}$. Suppose to have proved the claim for any complex with less faces than $\mathcal{C}$. Remove from $\mathcal{C}$ an extremal face (see Definition \ref{extremal}) and any edge and vertex belonging to its boundary and not belonging to the boundary of another face. We obtain a 2-cell complex with less faces than $\mathcal{C}$, which is still connected, whose faces are also faces of $\mathcal{C}$ (therefore they are stems or corollas) and whose edges are also edges of $\mathcal{C}$ (therefore given an edge $e$ either $e$ or $e^{-1}$ are contained in the boundary of some face). Moreover the boundary cycle of this complex is a subpath of that of $\mathcal{C}$, therefore its orientation is compatible with the orientation of its faces. Since $\mathcal{C}'$ has less faces than $\mathcal{C}$, by induction hypothesis $\mathcal{C}'$ belongs to $\mathcal{L}$. $\mathcal{C}$ is obtained from $\mathcal{C}'$ by adding to it the face, the edges and vertices removed; $\mathcal{C}$ is thus equal to the grafting of an element of $\mathcal{L}$ (with only one face) into $\mathcal{C}'$ and therefore belongs to $\mathcal{L}$. 
                        \end{proof}

\begin{proposition} \label{th2} Let $\mathcal{C}$ be an element of $\mathcal{L}$. Then \begin{enumerate}
  
  \item $\mathcal{C}$ is planar and contractible;
  
  \item the boundary cycle of $\mathcal{C}$ has no repeated edges and its only non-simple vertices are the initial vertices of its stems and corollas (except at most the initial vertex of $\mathcal{C}$). 
  \end{enumerate}
\end{proposition}

 \begin{proof} Let $\sigma$ be an SLA in $\mathcal{L}$, let $\mathcal{C}$ be the result of $\sigma$ and let $n$ be the number of steps of $\sigma$. We prove by induction on $n$ that $\mathcal{C}$ being this evident for $n=1$.

   Let $n>1$ and the claim be true for all the SLA's in $\mathcal{L}$ with less steps than $\sigma$. The last step of $\sigma$ is a grafting of a complex $\mathcal{C}_1$ into a complex $\mathcal{C}_2$ at vertices $v_1$ and $v_2$ of $\mathcal{C}_1$ and $\mathcal{C}_2$ respectively. By induction hypothesis  $\mathcal{C}_1$ and $\mathcal{C}_2$ verify the properties of the claim.

By 1 of Theorem \ref{th1} the faces of $\mathcal{C}$ are its stems and corollas and they are planar by the First and Second case. Theor. 27 of \cite{Whitn} says that a graph is planar if its unseparable components are; therefore $\mathcal{C}$ is planar since its unseparable components are its stems and corollas. Ex. 23 in Chapter 0 of \cite{H} says that the union of two contractible 2-complexes is contractible if their intersection is contractible, which is the case for $\mathcal{C}$ being equal to the union of $\mathcal{C}_1$ and $\mathcal{C}_2$ and being the intersection of the latter equal to a single vertex.

The boundary cycle of $\mathcal{C}$ has no repeated edges because it is the insertion of that of $\mathcal{C}_2$ into that $\mathcal{C}_1$ and both have no repeated edges by induction hypothesis. 

All the vertices that are non-simple in $\mathcal{C}_1$ and $\mathcal{C}_2$ are also non-simple in $\mathcal{C}$. This means that the initial vertices of the stems and corollas of $\mathcal{C}$ are non-simple, except at most the initial vertices of $\mathcal{C}_1$ and $\mathcal{C}_2$. Since the intersection of $\mathcal{C}_1$ and $\mathcal{C}_2$ consists only in the initial vertex of $\mathcal{C}_2$, the latter is non-simple and furthermore all the vertices that are simple in $\mathcal{C}_1$ and $\mathcal{C}_2$ are still simple in $\mathcal{C}$. Finally the initial vertex of $\mathcal{C}_1$ coincides with the initial vertex of $\mathcal{C}$.
             \end{proof}

 Theorem \ref{maithe2} says that $L_g$ is the set of all relators, reduced and non-reduced. This implies the following 

\medskip

\textbf{Corollary \ref{vKlemma}} \textit{ The set of relators (reduced and non) of the presentation $\langle \, X \, | \, R \, \rangle$ coincides with set of labels of elements of $\mathcal{L}$.}

\section{Preliminary results I} \label{prel1}

We recall that \textit{corollas} have been introduced in Definition \ref{overlineR}; the set $L_g$ has been introduced in Definition \ref{maiobj}. The words of the form $w w^{-1}$ (where $w$ is a reduced word) are called \textit{stems} (Definition \ref{stem}).

The term SLA without other specifications means SLA \textit{relative to L}. We also recall that a stem or a corolla of an SLA is called \textit{a base element} of that SLA (Definition \ref{basele}). The relation of “coming from" for components of outputs was introduced in Definition \ref{coming from}

\begin{theorem} \label{stecor} Let $l \in L_g$ be computed by $\sigma$. Then every component of $l$ comes from a base element of $\sigma$. \end{theorem}

\begin{proof} Let $n$ be the number of steps of $\sigma$; we prove the claim by induction on $n$. If $n=1$ then $l$ is a base element and the claim is evident. Let $n>1$ and the claim be true for every $n'<n$; $l$ is the insertion of a preceding output $l_1$ into another one $l_2$. Let $\sigma_1$ and $\sigma_2$ be the pSLsA's (Definition \ref{pSLsA}) computing $l_1$ and $l_2$; then every component of $l$ comes from $l_1$ or $l_2$ and every base element of $\sigma$ is a base element of $\sigma_1$ or of $\sigma_2$. The claim follows from induction hypothesis because $\sigma_1$ and $\sigma_2$ have less steps than $\sigma$.
            \end{proof}

\begin{definition} \label{flow} \rm The result of an SLA whose base elements are all stems is called a \textit{stem element}. The insertion (Definition \ref{insert}) of a corolla into a stem $w w^{-1}$ at $w$ is called a \textit{flower}. An insertion of a corolla into a stem element is called a \textit{flower element}. 
\end{definition}

A stem is a stem element and a flower is a flower element. The function $A$ has been introduced in Definition \ref{A}.

\begin{remark} \label{A(s)=0} \rm If $s$ is a stem element then $A(s)=0$. 
\end{remark}

The reduced form of a stem element is $1$. The converse is proved in the following

\begin{proposition} \label{=1}  Let $w$ be a word whose reduced form is $1$. Then  $w$ is a stem element.  \end{proposition}

\begin{proof} Let $w:=x_1 \cdots x_m$ where the $x_i$ are letters; we prove the claim by induction on $m$, being trivial for $m=2$. Let $m>2$ and the claim be true for every $m'<m$. Since $x_1 \cdots x_m=1$ in $\mathcal{F}(X)$, then there exists $i:1\leqslant i\leqslant m-1$ such that $x_{i+1}=x_i^{-1}$ (otherwise $w$ would be reduced and different from $1$). This implies that $w':=x_1 \cdots x_{i-1} x_{i+2} \cdots x_m$ is equal to $1$ in $\mathcal{F}(X)$. By induction hypothesis $w'$ is a stem element and thus there exists a straight line algorithm $\sigma$ whose base elements are all stems and whose result is $w'$. If we add to $\sigma$ a base step equal to the stem $x_i \, x_{i+1}$ and a step equal to the insertion of $x_i \, x_{i+1}$ into $w'$ at $x_{i-1}$, then we have obtained an SLA whose base elements are all stems and whose result is $w$.
      \end{proof}

We recall that $\rho(w)$ (Definition \ref{rho}) denotes the reduced form of $w$.

\begin{proposition} \label{preceding} \begin{enumerate}

  \item An insertion of a stem element into another one is still a stem element;
  
  \item a cyclic conjugate of a stem element is a stem element;
  
  \item the inverse of a stem element is a stem element;
   
   \item a cyclic conjugate of a g. corolla is a g. corolla.

\end{enumerate}   \end{proposition}

\begin{proof}\begin{enumerate}

  \item Trivial by virtue of Proposition \ref{=1}.
  
   \item Let $w$ be a stem element and let $w'$ be a cyclic conjugate of $w$. There exist words $u$ and $v$ such that $w=uv$ and $w'=vu$. Since $\rho(uv)=1$ then $\rho(u)\rho(v)=1$, $\rho(u)=\rho(v)^{-1}$, $\rho(v)\rho(u)=1$ and finally $\rho(vu)=1$. Proposition \ref{=1} implies that $vu$ is a stem element. 
   
    \item Trivial by Proposition \ref{=1}.
   
   \item Trivial because $\overline{R}$ is closed with respect to cyclic conjugation.
    \end{enumerate}
   \end{proof}

   \begin{proposition} \label{mesc} Let $w, u, v_1, \cdots, v_{m-1} \in L_g$ and let $w:=x_1 \cdots x_m, u:=u' u''$. Then 
        $$w':=u' x_1 v_1 \cdots x_{m-1} v_{m-1} x_m u'' \in L.$$
 Let $\sigma, \tau, \tau_1, \cdots, \tau_{m-1}$ be SLA's computing respectively $w, u, v_1, \cdots$, $v_{m-1}$; then there exists an SLA $\sigma'$ computing $w'$ such that 
  \begin{equation} \label{do} A(\sigma')=A(\sigma) + A(\tau) + A(\tau_1) + \cdots + A(\tau_{m-1}). \end{equation}    
In particular
   \begin{equation} \label{re} A(w') \leqslant A(w) + A(u) + A(v_1) + \cdots + A(v_{m-1}). \end{equation} 
Finally, if $M, N, N_1, \cdots, N_{m-1}$ are CMDR for $\sigma, \tau, \tau_1, \cdots, \tau_{m-1}$ respectively (Definition \ref{CMDR}), then $M\cup N \cup N_1\cup \cdots \cup N_{m-1}$ is a CMDR for $\sigma'$.                                                                                                                  \end{proposition}

\begin{proof}  Let $\sigma, \tau, \tau_1, \cdots, \tau_{m-1}$ be SLA's  computing $w, u, v_1, \cdots, v_{m-1}$ respectively. We define $\sigma'$ as the SLA whose steps are all the steps of  $\sigma$, $\tau$, $\tau_1$, $\cdots$, $\tau_{m-1}$, plus the insertions of $v_i$ at $x_i$ for every $i$ and finally the insertion of $x_1 v_1 \cdots x_{m-1} v_{m-1} x_m$ into $u$ at $u'$. $\sigma'$ computes $w'$ and verifies (\ref{do}). Since $A(w') \leqslant A(\sigma')$ by Proposition \ref{A}, to prove (\ref{re}) it is sufficient to take $\sigma, \tau, \tau_1, \cdots, \tau_{m-1}$ such that $A(\sigma)=A(w)$, $A(\tau)=A(u)$, $A(\tau_1)=A(v_1)$, $\cdots, A(\tau_{m-1})=A(v_{m-1})$. Finally by Remark \ref{union}, the union of $M, N, N_1, \cdots, N_{m-1}$ is a CMDR for $\sigma'$ if the latter are CMDR for $\sigma, \tau, \tau_1, \cdots, \tau_{m-1}$ respectively.
        \end{proof}

\begin{corollary} \label{promesc} \begin{enumerate}

  \item Let $\sigma$ and $\sigma'$ be SLA's with results $w$ and $w'$ and let $v$ be an insertion of $w$ into $w'$. Then there exists an SLA $\sigma''$ computing $v$ such that $A(\sigma'')=A(\sigma)+A(\sigma')$; in particular $A(v) \leqslant A(w)+A(w')$. Moreover the union of a CMDR for $\sigma$ and of one for $\sigma'$ is a CMDR for $\sigma''$.

  \item Let $\sigma$ be an SLA with result $w$ and let $x$ be a letter. There exists an SLA $\sigma'$ computing $xwx^{-1}$ and such that $A(\sigma')=A(\sigma)$; in particular $A(xwx^{-1}) \leqslant A(w)$. Moreover a CMDR for $\sigma$ is a CMDR for $\sigma'$. 
  \end{enumerate}
     \end{corollary}

\begin{proof} Follows from Proposition \ref{mesc}. \end{proof}

 \begin{remark} \label{A(f)} \rm Let $f$ be a flower element which is the insertion of a g. corolla $c$ into a stem element $s$. By Proposition \ref{mesc} we have that $A(f)\leqslant A(c)+A(s)$ and since $A(s)=0$ by Remark \ref{A(s)=0}, then $A(f)\leqslant A(c)$. 
\end{remark}

\begin{proposition} \label{cc'} Let $c$ and $c'$ be g. corollas and let $c c'=f_1 z z^{-1} f_2$ (where $f_1$ and $f_2$ are words and $z$ a letter). Then $f_1 f_2$ belongs to $L_g$, in particular it is a stem if $c'=c^{-1}$ or it is an insertion of a stem into a reduced flower (i.e., a reduced word which is a flower) whose g. corolla is $\pi(c,c')$. If moreover $c$ and $c'$ are proper g. corollas (Definition \ref{propcor}) then $A(f_1 f_2)\leqslant A(c)+A(c')$. 
\end{proposition}

\begin{proof} If $c'=c^{-1}$ then $c c'$ is a stem and $z$ is the last letter of $c$ and $z^{-1}$ the first of $c'$ since $c$ and $c'$ are reduced. In this case $c=f_1 z$, $c'=z^{-1} f_2$ and $f_1 f_2$ is a stem. Therefore $A(f_1 f_2)=0\leqslant A(c)+A(c')$.

Let $c'\neq c^{-1}$; then the reduced product of $c$ by $c'$ is of the form $u d u^{-1}$ (which is a reduced flower) where $d:=\pi(c, c')$ is their cyclically reduced product; $d$ is a g. corolla because there is cancellation in $\pi(c, c')$. By Definitions \ref{reducprod} and \ref{scope}, $cc'$ is an insertion of the cancelled part $a a^{-1}$, which is a stem, into $u d u^{-1}$. Since $u d u^{-1}$ is reduced, $z z^{-1}$ is a subword of $a a^{-1}$. In particular since $a$ is reduced then $z$ is the last letter of $a$ and consequently is $z^{-1}$ the first one of $a^{-1}$. Therefore $a=bz$ for some reduced word $b$ and thus $f_1 f_2$ is the insertion of the stem $b b^{-1}$ into $u d u^{-1}$ at the same letter as $a a^{-1}$ is inserted into $u d u^{-1}$ to give $cc'$.
  
By Proposition \ref{mesc} we have that $A(f_1 f_2)\leqslant A(u u^{-1})+A(b b^{-1})+A(d)$. Since $u u^{-1}$ and $b b^{-1}$ are stems then $A(u u^{-1})=A(b b^{-1})=0$ by Remark \ref{A(s)=0} and $A(f_1 f_2)\leqslant A(d)$. Since $d=\pi(c, c')$, then $A(d)\leqslant \eta\big(\pi(c, c')\big)$ by Remark \ref{etaA} and $\eta\big(\pi(c, c')\big) \leqslant \eta(c)+\eta(c')$ by Remark \ref{eta(pi(cc'))}. Thus the inequality $A(f_1 f_2)\leqslant A(c)+A(c')$ follows from the fact that if $c$ and $c'$ are proper g. corollas then $A(c)=\eta(c)$ and $A(c')=\eta(c')$.   
    \end{proof}

\begin{lemma} \label{flower} Let $f$ be an insertion of a stem element into a flower element with g. corolla $c$ and let $f'$ be a cyclic conjugate of $f$. $f'$ is an insertion of two stem elements (possibly empty) into a flower element with g. corolla a cyclic conjugate of $c$.
 \end{lemma}

\begin{proof} $f$ is an insertion of a stem element $u$ into a flower element $s_1c s_2$, where $s_1 s_2$ is a stem element. Therefore $f$ is 
\begin{enumerate}
  
  \item either of the form $s_1' u s_1'' c s_2$, where $s_1' s_1''=s_1$,
  
  \item or of the form $s_1 c_1 u c_2 s_2$, where $c_1 c_2=c$,
  
  \item or of the form $s_1 c s'_2 u s''_2$, where $s_2' s_2''=s_2$.
  
\end{enumerate}

We can suppose that $f$ is of the form $t_1 c_1 v c_2 t_2$ where $c_1 c_2=c$ and $t_1 t_2$ and $v$ are stem elements such that

\begin{enumerate}
  
  \item $t_1=s_1' u s_1''$, $t_2=s_2$ and $v=1$ in case 1,
  
  \item $t_1=s_1$, $t_2=s_2$ and $v=u$ in case 2,
  
  \item $t_1=s_1$, $t_2=s'_2 u s''_2$ and $v=1$ in case 3.
  
\end{enumerate}

Since $f'$ is a cyclic conjugate of $f=t_1 c_1 v c_2 t_2$ then $f'$ is \begin{itemize}

  \item either of the form $t_1'' c_1 v c_2 t_2 t_1'$, where $t_1' t_1''=t_1$, 
  
  \item or of the form $c_1'' v c_2 t_2 t_1 c_1'$, where $c_1'c_1''=c_1$,
  
  \item or of the form $v_2 c_2 t_2 t_1 c_1 v_1$, where $v_1 v_2=v$,
  
  \item or of the form $c_2'' t_2 t_1 c_1 v c_2' $, where $c_2' c_2''=c_2$
  
  \item or of the form $t_2'' t_1 c_1 v c_2 t_2'$, where $t_2' t_2''=t_2$.
  
\end{itemize}

The following are cyclic conjugates of $c$: $c_1 c_2$, $c_1'' c_2 c_1'$, $c_2 c_1$, $c_2'' c_1 c_2'$. By Parts 1 and 2 of Proposition \ref{preceding}, the following are stem elements: $t_1'' t_2 t_1'$, $t_2 t_1$, $v_2 v_1$, $t_2'' t_1 t_2'$.

Therefore in all the cases $f'$ is either a flower element or an insertion of one or two stem elements into a flower element (a g. corolla is a flower element). The g. corollas of those flower elements are cyclic conjugates of $c$.
        \end{proof}

\begin{proposition} \label{bigg. corolla} Let $c$ and $c'$ be g. corollas, let $f$ be an insertion of $c'$ into $c$ and let $f=f_1 z z^{-1} f_2$ (where $f_1$ and $f_2$ are words and $z$ a letter). Then $f_1 f_2$ belongs to $L_g$, in particular it is a stem element if $c'=c^{-1}$ or it is an insertion of a stem element into a flower element with g. corolla a cyclic conjugate of $\pi(d, c')$, where $d$ is a cyclic conjugate of $c$. If moreover $c$ and $c'$ are proper g. corollas then $A(f_1 f_2)\leqslant A(c)+A(c')$.
         \end{proposition}

\begin{proof} There exist words $c_1$ and $c_2$ such that $c=c_1 c_2$ and $f=c_1 c' c_2$. Let $d=c_2 c_1$; then $f$ is a cyclic conjugate of $d c'$, which is the product of two g. corollas ($d$ is a g. corolla by Part 4 of Proposition \ref{preceding}). By Definitions \ref{reducprod} and \ref{scope}, $d c'$ is an insertion of the cancelled part $a a^{-1}$, which is a stem, into the reduced flower $u d' u^{-1}$, where $d'=\pi(d, c')$. By Lemma \ref{flower}, $f$ is an insertion of two stem elements into a flower element with g. corolla a cyclic conjugate of $d'$.   
  
 Since $u d' u^{-1}$ is the reduced form of $d c'$, then either $z z^{-1}$ is a subword of the cancelled part $a a^{-1}$, or $z$ is the first letter of $u$ and consequently $z^{-1}$ is the last one of $u^{-1}$. In both cases $f_1 f_2$ is a cyclic conjugate of an insertion of a stem into a flower with g. corolla $d'$. Therefore by Lemma \ref{flower} it is an insertion of two stem elements into a flower element with g. corolla a cyclic conjugate of $d'$. Let $d''$ be that cyclic conjugate. Thus $A(f_1 f_2)\leqslant A(d'')$ by Proposition \ref{mesc} and by Remark \ref{A(s)=0}; $A(d'')\leqslant \eta(d'')$ by Remark \ref{etaA}; $\eta(d'')=\eta(d')$ by Remark \ref{eta(psi)}; $\eta(d') \leqslant \eta(d)+\eta(c')$ by Remark \ref{eta(pi(cc'))} and $\eta(d)=\eta(c)$ by Remark \ref{eta(psi)}. The final inequality follows from the fact that if $c$ and $c'$ are proper g. corollas then $A(c)=\eta(c)$ and $A(c')=\eta(c')$.  
   \end{proof}

\begin{remark} \label{dopo} \rm Let $c$ and $c'$ be g. corollas and let $c=c_1 c_2$ where $c_1$ and $c_2$ are words. By virtue of Propositions \ref{cc'} and \ref{bigg. corolla}, if $f$ is a word such that $\rho(f)=\rho(cc')$ or $\rho(f)=\rho(c_1 c' c_2)$ then either $c'=c^{-1}$ or there exists an SLA computing $f$ and that has only one g. corolla, which we call $d$, such that $\eta(d)\leqslant \eta(c)+\eta(c')$. The word $d$ is equal to $\pi(c, c')$ if $\rho(f)=\rho(cc')$ or equal to a cyclic conjugate of $\pi(c_2 c_1, c')$ if $\rho(f)=\rho(c_1 c' c_2)$. 

Let $\sigma_1$ be an SLA, let $M_1$ be its multiset of g. corollas and let $c$ and $c'$ be elements of $M$. Suppose that $c'=c^{-1}$ and that $\sigma_2$ is an SLA such that $M\setminus{\{c, c'\}}$ is the multiset of g. corollas of $\sigma_2$\footnote{If $M_1:=(S, \lambda_1)$ then the multiset difference $M\setminus{\{c, c'\}}$ is the multiset $M_2:=(S, \lambda_2)$ where $\lambda_2(c)=\lambda_1(c)-1$, $\lambda_2(c')=\lambda_1(c')-1$ and $\lambda_2$ coincides with $\lambda_1$ in $S\setminus{\{c, c'\}}$}. Then it is obvious that $A(\sigma_2) < A(\sigma_1)$.

Now suppose that $c'\neq c^{-1}$ and let $d$ as above. Let $\sigma_1$ and $\sigma_2$ be SLA's, let $M_1$ and $M_2$ be respectively their multisets of g. corollas and let $c$ and $c'$ be elements of $M_1$. Suppose also that $M_1$ and $M_2$ coincide except that in $M_2$, $c$ and $c'$ are replaced by $d$, that is $M_1\setminus{\{c, c'\}}=M_2\setminus{\{d\}}$. Since $\eta(d)\leqslant \eta(c)+\eta(c')$, then $A(\sigma_2)\leqslant A(\sigma_1)$.

We have that $A(\sigma_2)=A(\sigma_1)$ if and only if $\eta(d)=\eta(c)+\eta(c')$. Let $\tau$ and $\tau'$ be SLA's in $\overline{R}$ computing $c$ and $c'$ respectively and such that $\eta(\tau)=\eta(c)$ and $\eta(\tau')=\eta(c')$. We suppose that $\eta(d)=\eta(c)+\eta(c')$; this implies that there exists an SLA in $\overline{R}$, which we call $\tau''$, which computes $d$ and such that $\eta(\tau'')=\eta(d)$. $\tau''$ is constructed in the following way. If $d$ is the g. corolla of $cc'$ then $\tau''$ is as constructed in Remark \ref{eta(pi(cc'))}. If $d$ is the g. corolla of $c_1 c' c_2$ then we take all the steps of $\tau$ and $\tau'$ and we add: the cyclic conjugation from $c$ to $c_2 c_1$; the cyclically reduced product of $c_2 c_1$ by $c'$ which we denote $d'$; and finally the cyclic conjugation from $d'$ to $d$. We have that $\eta(\tau'')=\eta(\tau)+\eta(\tau')=\eta(d)$ because a cyclic conjugation does not change the value of $\eta$ and with a cyclically reduced product the value of $\eta$ is the sum. 

Finally the equality $\eta(d)=\eta(c)+\eta(c')$ implies that a CMDR for $\sigma_1$ is a CMDR also for $\sigma_2$, because $M_1$ and $M_2$ coincide except that $d$ replaces $c$ and $c'$ in $M_2$ and because the multiset of base elements of $\tau''$ is the union of those of $\tau$ and $\tau'$.
 \end{remark}

   \section{Preliminary results II} \label{prel2}

In this section we continue proving results necessary for the proof of the Lemma \ref{maile}. We define \textit{straight line subalgorithms}, we introduce the intuitive notions of \textit{ramifications} and \textit{surround} and we prove some technical lemmas.

\begin{definition} \label{contain} \rm Let $\sigma$ be an SLA with result $w$ and let $s_1$ and $s_2$ be two steps of $\sigma$; since we identify a step with its output, we consider $s_1$ and $s_2$ as subwords of $w$. We say that \textit{$s_1$ comprises $s_2$} if in $w$ all the letters of $s_2$ are comprised between two consecutive letters of $s_1$. We say that \textit{$s_1$ precedes $s_2$} if in $w$ any letter of $s_1$ precedes any letter of $s_2$.
      \end{definition}

\begin{remark} \label{compr} \rm If $s_2$ is inserted into $s_1$ at a letter which is not the last then $s_1$ comprises $s_2$. Consider the reflexive closure of the relation “being comprised in"; it is a partial order equal to that of Definition \ref{extremal}. Since the set of steps of $\sigma$ is finite, then by Zorn's Lemma \cite{zorn} any step of $\sigma$ either is minimal or comprises a minimal step. This minimal step is thus a contiguous subword of the given step.    \end{remark}

\begin{proposition} \label{consecutive} Given two steps of an SLA such that none of them depends on the other, then one of them comprises or precedes the other. \end{proposition}   
   
\begin{proof} Let $s_1$ and $s_2$ be two steps of an SLA and let $t$ be the first step containing both of them, that is $t$ is the first step such that $s_1$ and $s_2$ are subwords of $t$. We have that $t$ is different from $s_1$ and $s_2$ because by hypothesis none of them depends on the other. This means that $t$ is the insertion of a preceding step $t_2$ into another $t_1$ with $t_2$ depending on $s_2$ and $t_1$ on $s_1$. 

If $t$ is the product of $t_1$ by $t_2$, then every letter of $t_1$ precedes every letter of $t_2$, therefore $s_1$ precedes $s_2$. Suppose on the contrary that $t$ is the insertion of $t_2$ into $t_1$ at a letter $x$ that is not the last one. Let $x'$ be the last letter of $t_1$ coming from $s_1$ and preceding or equal to $x$. If $x'$ is the last letter of $s_1$ then it precedes the first one of $t_2$ and therefore $s_1$ precedes $s_2$. Suppose that $x'$ is not the last letter of $s_1$. Let $x''$ be the first letter of $t_1$ following $x$ and coming from $s_1$; then all the letters of $t_2$ (and thus of $s_2$) are comprised between $x'$ and $x''$, that is between two consecutive letters of $s_1$ and $s_1$ comprises $s_2$. Since the steps following $t$ do not change the relative order of the letters of $s_1$ and $s_2$, we have proved the claim.
                  \end{proof}

\begin{remark} \label{alcap} \rm We prove that every step of an SLA contains as a contiguous subword a base element which is minimal with respect to the order defined in Remark \ref{compr}.

Let $\sigma$ be an SLA and let $s$ be a step of $\sigma$; $s$ is or comprises a minimal step $s'$ and therefore contains it as a contiguous subword. If this step is a base step the claim is proved; if it is not then any base step used by $s'$ is minimal because otherwise $s'$ would not be minimal and therefore this minimal base step is a contiguous subword of $s'$ (and therefore of $s$).

This means that if $w$ is the result of $\sigma$ and if $s$ is a step of $\sigma$, then $s$ is a contiguous subword of $w$ if and only if $s$ is minimal. If $w$ is reduced then no stem of $\sigma$ is minimal because a stem is not reduced and $w$ cannot have contiguous subwords which are non-reduced. Therefore if $w$ is reduced then there is at least a g. corolla of $\sigma$ which is a contiguous subword of $w$.
   \end{remark}

\smallskip

Let $\sigma$ be an SLA and let $M:=(B, \lambda)$ be its multiset of base elements (Definition \ref{basele}). We represent $M$ as the set of pairs $(b, k)$ where $b \in B$ and $k$ is a non-zero natural number less or equal to the multiplicity of $b$. For instance, if the multiplicty of an element $b$ is $3$, then $M$ contains $(b, 1)$, $(b, 2)$, $(b, 3)$ and does not contain $(b, k)$ for $k>3$. $b$ is called \textit{the underlying element of $(b, k)$}. Sometimes we will identify the pair $(b, k)$ with $b$. There is a natural bijection between the base steps of $\sigma$ and $M$, given by sending a base step $s$ to $(b, k)$ if $s$ is the $k$-th step of $\sigma$ equal to $b$.

Let $\sigma_1$ and $\sigma_2$ be two SLA's with multisets of base elements $M_1$ and $M_2$ respectively. An \emph{homomorphism of multisets} is an application $\omega: M_1 \rightarrow M_2$ that sends an element of $M_1$ to an element of $M_2$ with the same underlying element, for instance sends $(b, 3)$ to $(b, 1)$. Since a base element cannot be at the same time a stem and a g. corolla, an homomorphism sends stems to stems and g. corollas to g. corollas. If $\omega$ is injective then for every $b \in B$ the multiplicity of $b$ in $\sigma_1$ is less or equal to the multiplicity in $\sigma_2$; this means in particular that $A(\sigma_1)\leqslant A(\sigma_2)$ because to every g. corolla of $\sigma_1$ corresponds the same g. corolla in $\sigma_2$. If there is an element of $B$ with non-zero multiplicity in $\sigma_1$ and zero multiplicity in $\sigma_2$, then no homomorphism can be defined from $M_1$ to $M_2$.

\begin{definition} \label{corresp} \rm Let $\sigma_1$ and $\sigma_2$  be SLA's with results $w_1$ and $w_2$, with multisets of base elements $M_1$ and $M_2$ respectively and let $\omega: M_1 \rightarrow M_2$ be an homomorphism. One letter of $w_1$ and one of $w_2$ are said to \textit{correspond by $\omega$} if there exists $\mu \in M_1$ such that the two letters come (Definition \ref{coming from}) from the same letter of $\mu$ and $\omega(\mu)$ respectively (we recall that $\mu$ and $\omega(\mu)$ have the same underlying element).
   \end{definition}

We now define a notion of straight line subalgorithm which generalizes that of Remark \ref{pSLsA}.

\begin{definition} \label{SLsA} \rm Let $\sigma_1$ and $\sigma_2$ be SLA's with results $w_1$ and $w_2$ and with multisets of base elements $M_1$ and $M_2$ respectively. $\sigma_1$ is a \textit{straight line subalgorithm (SLsA)} of $\sigma_2$ if $w_1$ is a (not necessarily contiguous) subword of $w_2$ and if there exists an injective homomorphism from $M_1$ to $M_2$ such that every letter of $w_1$ corresponds by $\omega$ to the same letter in $w_2$ (since $w_1$ is a subword of $w_2$, every letter of $w_1$ is also a letter of $w_2$). In this case we say that \textit{$w_1$ is a part of $w_2$}.
   \end{definition}
   
A proper straight line subalgorithm (Remark \ref{pSLsA}) is a straight line subalgorithm.

\begin{remark} \label{N} \rm Let $\sigma$ be an SLA with result $w$ and with multiset of base elements $M$ and let $N$ be a \textit{sub-multiset} of $M$ (denoted $N\subset M$), that is if $M=(B, \lambda)$ then $N=(B, \lambda_0)$ where $\lambda_0(b)\leqslant \lambda(b)$ for every $b\in B$. There is an SLsA of $\sigma$ whose result is the subword of $w$ whose letters come from the elements of $N$; its multiset of base elements is $N$. This SLsA and its result are called \textit{the SLsA and the part determined by $N$}. It is constructed in the following way. If $N$ has only one element, take the SLsA with a single step equal to this element. Let $|N|>1$ and let the construction be done for every $N'$ with less elements than $N$. Let $\nu$ be an element minimal in $N$ with respect to the order of Definition \ref{contain} and let $N'=N\setminus{\{\nu\}}$, that is $N'=(B, \lambda')$ where $\lambda'(\nu)=\lambda_0(\nu)-1$ and $\lambda'$ coincide with $\lambda_0$ on $B\setminus{\{\nu\}}$. Let $\tau'$ be the SLsA of $\sigma$ defined by $N'$ and let $f'$ be its result. Let $f$ be the subword of $w$ whose letters come from elements of $N$. Since $\nu$ is minimal, there exist words $f'_1$ and $f'_2$ such that $f=f'_1 \, \nu \, f'_2$ and $f'=f'_1 \, f'_2$. If we add to $\tau'$ a base step equal to $\nu$ and another one equal to the insertion of $\nu$ into $f'$ at $f'_1$, then we have constructed an SLsA with result $f$ and with multiset of base elements equal to $N$.
                         \end{remark}

 \begin{definition} \label{corresponding}  \rm  Let $\sigma_1$ and $\sigma_2$ be SLA's with results $w_1$ and $w_2$ and with multisets of base elements $M_1$ and $M_2$ respectively and suppose given an homomorphism $\omega$ from $M_1$ to $M_2$. Let $N_1\subset M_1$ and let $f_1$ be the part of $w_1$ determined by $N_1$ (Remark \ref{N}). Let $N_2:=\omega(N_1)$ and let $f_2$ be the part of $w_2$ determined by $N_2$. We say that \textit{$f_1$ corresponds to $f_2$ by $\omega$.}    \end{definition}

\begin{definition} \label{ramsur}  \rm Let $\sigma$ be an SLA with result $w:=x_1 \cdots x_m$, let $s$ be a step\footnote{in particular, since we identify a step with its output, $s$ is a not necessarily contiguous subword of $w$.} of $\sigma$ and let $x_i$ and $x_k$ (with $i<k$) be letters of $w$ coming from two consecutive letters of $s$. The subword $x_{i+1} \cdots x_{k-1}$ is called \textit{a ramification from $s$}. Let $x_f$ and $x_l$ be letters of $w$ coming respectively from the first and the last letters of $s$; the subword $x_1 \cdots x_{f-1} \, x_{l+1} \cdots x_m$ is called \textit{the surround of $s$}. $x_1 \cdots x_{f-1}$ is called \textit{the preceding of $s$} and $x_{l+1} \cdots x_m$ \textit{the following of $s$}.
   \end{definition}

 \begin{remark} \label{cons2}  \rm Given two steps $s$ and $s'$ of an SLA, we have that $s$ comprises $s'$ (Remark \ref{compr}) if and only if $s'$ is contained in a ramification from $s$, if and only if the surround of $s'$ contains $s$; $s$ precedes $s'$ if and only if the following of $s$ contains $s'$, if and only if the preceding of $s'$ contains $s$.
     \end{remark}

\begin{proposition} \label{rs} Let $\sigma$ be an SLA and let $s$ and $s'$ be steps of $\sigma$ such that none of them depends on the other. Then

\begin{enumerate}

  \item if a ramification from $s$ contains a letter of $s'$ then it contains all the letters of $s'$; 
  
  \item if the surround of $s$ contains a letter of $s'$ then it contains all the letters of $s'$.

   \end{enumerate}
          \end{proposition}

\begin{proof} It is a consequence of Proposition \ref{consecutive} in view of Remark \ref{cons2}. If a ramification from $s$ contains a letter of $s'$ then it is not possible that $s$ precedes $s'$, nor that $s'$ precedes $s$, nor that $s'$ comprises $s$. Thus $s$ comprises $s'$, that is all the letters of $s'$ are contained in the given ramification from $s$. Since two different ramifications have no letters in common then Part 1 is proved.

If the surround of $s$ contains a letter of $s'$ then it is not possible that $s$ comprises $s'$, nor that $s$ precedes $s'$. If $s'$ comprises $s$ then by Remark \ref{cons2} the surround of $s$ contains $s'$; if $s'$ precedes $s$ then the preceding of $s$ (and therefore its surround) contains $s'$.
 \end{proof}

\begin{proposition} \label{fam} Let $\sigma$ be an SLA whose result is $w$ and let $s$ be a step of $\sigma$. The surround of $s$ and any ramification from $s$ are parts of $w$. \end{proposition}

\begin{proof} Let $q$ be the surround of $s$; we have to prove that there exists an SLsA of $\sigma$ computing $q$. Let $M$ be the multiset of base elements of $\sigma$ and let $N$ be the sub-multiset of $M$ of the elements which have at least one letter contained in $q$; let $f$ be the result of the SLsA defined by $N$ (Remark \ref{N}). $f$ is a part of $w$ and contains $q$ as a subword. Vice versa we prove that $q$ contains $f$ as a subword. Let $s'$ be a base step which has at least one letter in common with $q$. $s'$ does not use $s$ because a base step does not any step; $s$ does not use $s'$ because otherwise $s'$ must be a subword of $s$ and instead at least one letter of $s'$ is contained in $q$ (which has no letters in common with $s$). Therefore by Proposition \ref{rs} $q$ contains every base element with which it has at least a letter in common; thus $q=f$. 

Analogously we do for a ramification.
     \end{proof}

\begin{definition} \label{def.cons} \rm Let $\tau, \tau_0, \tau_1, \cdots, \tau_{m-1}$ be SLA's with results respectively the words $w, q, r_1, \cdots, r_{m-1}$. Let $w:=x_1 \cdots x_m$ and $q:=q_0 \, q_1$ (with the $x_i$ letters and $q_0$ and $q_1$ words) and let $\iota_0$ be the insertion of $w$ into $q$ at $q_0$; its result is $w_0:= q_0 \, w \, q_1$. Let $\iota_1$ be the insertion of $r_1$ into $w_0$ at $x_1$; its result is $w_1:=q_0 \, x_1 \, r_1 \, x_2 \cdots x_m \, q_1$. $\cdots$ Let $\iota_{m-1}$ be the insertion of $r_{m-1}$ into $w_{m-2}$ at $x_{m-1}$; its result is $q_0 \, x_1 \, r_1 \cdots x_{m-1} \, r_{m-1} \, x_m  \, q_1$.

Then $\sigma:=(\tau, \tau_0, \iota_0, \tau_1, \iota_1, \cdots, \tau_{m-1}, \iota_{m-1})$ is an SLA, $\tau$ is a pSLsA computing $w$, $q$ is the surround of $w$ and $r_1, \cdots, r_{m-1}$ are the ramifications from $w$.

We say that \textit{$\sigma$ defines consecutively the insertions into $w$} and that $\tau_0$, $\tau_1$, $\cdots$, $\tau_{m-1}$ are the \textit{pSLsA's of $\sigma$ which compute respectively the surround and the ramifications from $w$}.
     \end{definition}

Given an SLA $\sigma$ and given a step $s$, we want to prove that there exists an SLA “equivalent" to $\sigma$ (in a sense that we are going to specify) which defines consecutively the insertions into $s$.

\begin{definition} \label{equiSLA} \rm If $\sigma$ and $\sigma'$ are SLA's such that any of the two is an SLsA of the other (Definition \ref{SLsA}), then we say that $\sigma$ and $\sigma'$ are equivalent.
\end{definition}

\begin{remark} \label{nikei} \rm Let $\sigma_1$ and $\sigma_2$ be SLA's with multisets of base elements $M_1$ and $M_2$ respectively. Then $\sigma_1$ and $\sigma_2$ are equivalent if and only if their results are equal (let $w$ be their result) and there exists an isomorphism $\omega: M_1 \rightarrow M_2$ such that any letter of $w$ corresponds (Definition \ref{corresp}) to itself by $\omega$. If $\sigma$ and $\sigma'$ are equivalent then $A(\sigma)=A(\sigma')$.
\end{remark}

\begin{proposition} \label{consec} Let $\sigma$ be an SLA and let $s$ be one of its steps. Then there exists an SLA $\sigma'$ equivalent to $\sigma$ and defining consecutively the insertions into the step of $\sigma'$ corresponding to $s$ (Definition \ref{corresponding}).
    \end{proposition} 

\begin{proof} Let $s:=x_1 \cdots x_m$. By Proposition \ref{fam} there exist SLsA's $\tau_0$, $\tau_1$, $\cdots$, $\tau_{m-1}$ computing respectively the surround $q$ and the ramifications $r_1$, $\cdots$, $r_{m-1}$ from $s$. 

We define the insertions $\iota_0, \iota_1, \cdots, \iota_{m-1}$ in the following way. $\iota_0$ is the insertion of $s$ into $q$ at $q_0$, where $q_0$ is the preceding of $s$; call $w_0$ its result. We have that $w_0=q_0 s q_1$ where $q_1$ is the following of $s$. We define recursively $\iota_j$ for $j=1, \cdots, m-1$ as the insertion of $r_j$ into $w_{j-1}$ at $x_j$.

Let $\tau$ be an SLA computing $s$. Then $\sigma':=(\tau, \tau_0, \iota_0, \tau_1, \iota_1, \cdots, \tau_{m-1}, \iota_{m-1})$ is an SLA, its result is the same of $\sigma$ and there is an evident isomorphism between its multiset of base elements and that of $\sigma$. Moreover $\sigma'$ defines consecutively the insertions into $s$.
    \end{proof}

\begin{lemma} \label{oso} Let $\sigma$ be an SLA whose result is $w:=x_1 \cdots x_m$, let $s:=y_1 \cdots y_p \, y_p^{-1} \cdots y_1^{-1}$ be a stem of $\sigma$, let $h$ and $h'$ be indices such that $x_h=y_n$ and $x_{h'}=y_n^{-1}$ for some $n: 1 \leqslant n \leqslant p$. Then there exist two SLA's $\sigma_1$ and $\sigma_2$ computing respectively 
$$x_1 \cdots x_{h-1} \, x_{h'+1} \cdots x_m \, \, \, \, \,  \, \, \, and  \, \, \, \, \, \, \, \, x_{h+1} \cdots x_{h'-1}$$
and such that $A(\sigma_1)+A(\sigma_2)=A(\sigma)$. Moreover the union of a CMDR for $\sigma_1$ and of one for $\sigma_2$ is a CMDR for $\sigma$.
\end{lemma}

\begin{proof} By Proposition \ref{consec} we can suppose that $\sigma$ defines consecutively the insertions into $s$. Let $q_0$ and $q_1$ be the preceding and the following of $s$ and let $\tau$ be the SLsA computing the surround $q:=q_0 q_1$. Let $r_1, \cdots, r_p$ be the ramifications from $s$ at $y_1, \cdots, y_p$ respectively and $r'_2, \cdots, r'_p$ the ones at $y^{-1}_2, \cdots, y^{-1}_p$. Let $\tau_1, \cdots, \tau_p$ and $\tau'_2, \cdots, \tau'_p$ be the SLsA's computing them. This means that 
       $$w=q_0  \textbf{\emph{y}}_\textbf{1}  r_1 \textbf{\emph{y}}_\textbf{2}  \cdots \textbf{\emph{y}}_\textbf{\emph{p}}  r_p \textbf{\emph{y}}_\textbf{\emph{p}}^{\textbf{-1}}  r'_p \cdots \textbf{\emph{y}}_\textbf{2}^{\textbf{-1}}  r'_2  \textbf{\emph{y}}_\textbf{1}^{\textbf{-1}}  q_1.$$
       (we write in bold the letters $y_i$).
Set 
$$s_1:=y_1 \cdots y_{n-1} \, y_{n-1}^{-1} \cdots y_1^{-1} \, \, \, \, \textrm{and} \, \, \, \,  s_2:=y_{n+1} \cdots y_p \, y_p^{-1} \cdots y_{n+1}^{-1}.$$
$s_1$ and $s_2$ are stems. Set 
$$\sigma_1:=(s_1, \tau, \tau_1, \cdots, \tau_{n-1}, \tau'_2, \cdots, \tau'_n, \iota_0, \iota_1, \cdots, \iota_{n-1}, \iota'_2, \cdots, \iota'_n)$$ 
and 
$$\sigma_2:=(s_2, \tau_n, \cdots, \tau_p, \tau'_{n+1}, \cdots, \tau'_p, \iota_n, \cdots, \iota_p, \iota'_{n+1}, \cdots, \iota'_p),$$
where $\iota_0$ is the insertion of $s_1$ into $q$ at $q_0$ and for $j\neq 1$, $\iota_j$ and $\iota'_j$ are the insertions of $r_j$ and of $r'_j$ at $y_j$ and $y_j^{-1}$ respectively. The results of $\sigma_1$ and $\sigma_2$ are  
          $$w_1:=q_0  \textbf{\emph{y}}_\textbf{1}  r_1 \cdots \textbf{\emph{y}}_\textbf{\emph{n}-1}  r_{n-1} r'_n \textbf{\emph{y}}_\textbf{\emph{n}-1}^{\textbf{-1}} r'_{n-1} \cdots \textbf{\emph{y}}_\textbf{2}^{\textbf{-1}} r'_2  \textbf{\emph{y}}_\textbf{1}^{\textbf{-1}}  q_1$$
    and
$$w_2:=r_n \textbf{\emph{y}}_\textbf{\emph{n}+1}  r_{n+1}  \cdots \textbf{\emph{y}}_\textbf{\emph{p}}  r_p \textbf{\emph{y}}_\textbf{\emph{p}}^{\textbf{-1}}  r'_p \cdots  \textbf{\emph{y}}_\textbf{\emph{n}+1}^{\textbf{-1}} r'_{n+1}$$
respectively and $w_1=x_1 \cdots x_{h-1} \, x_{h'+1} \cdots x_m$, $w_2= x_{h+1} \cdots x_{h'-1}$.  Finally, the equality $A(\sigma_1)+A(\sigma_2)=A(\sigma)$ and the last claim follow from Remark \ref{union}.
       \end{proof}

\begin{lemma} \label{faso} Let $\sigma$ be an SLA and let $w:=x_1 \cdots x_m$ be its result. If the first [respectively the last] letter of $w$ comes (Definition \ref{coming from}) from a g. corolla of $\sigma$, then there exists an SLA $\sigma'$ whose result is $x_2 \cdots x_m \, x_1$ [respectively $x_m \, x_1 \cdots x_{m-1}$], such that $A(\sigma')=A(\sigma)$ and such that a CMDR for $\sigma$ is a CMDR also for $\sigma'$.  \end{lemma}

\begin{proof} 
  Let $c:=y_1 \cdots y_p$ be the g. corolla from which comes $x_1$ [respectively $x_m$]. By Proposition \ref{consec} we can suppose that $\sigma$ defines consecutively the insertions into $c$. Since $x_1=y_1$ [respectively $x_m=y_p$] then the preceding [respectively the following] of $c$ is empty, therefore $w$ is equal to $\textbf{\emph{y}}_\textbf{1} r_1 \cdots r_{p-1}  \textbf{\emph{y}}_\textbf{\emph{p}} q$ [respectively to $q \textbf{\emph{y}}_\textbf{1} r_1 \cdots r_{p-1}  \textbf{\emph{y}}_\textbf{\emph{p}}$] where $q$ is the surround and the $r_j$ are the ramifications from $c$. This means that $\sigma$ is of the form 
    $$(c, \tau, \iota, \tau_1, \iota_1, \cdots, \tau_{p-1}, \iota_{p-1})$$
     where $\tau$ computes $q$, $\tau_j$ computes $r_j$, $\iota$ is the product $c q$ [respectively the product $q c$] and $\iota_j$ is the insertion of $r_j$ at $y_j$.

Set $d:=y_2 \cdots y_p \, y_1$ [respectively $d:=y_p \, y_1 \cdots y_{p-1}$] and 
    $$\sigma':=(d, \tau, \tau_1, \iota', \iota'_1, \tau_2, \iota_2, \cdots, \tau_{p-1}, \iota_{p-1})$$
    where $\iota'$ is the product $r_1 d$ [respectively $\iota'=\iota_1$] and $\iota'_1$ is the insertion of $q$ at $y_p$. $\sigma'$ is an SLA whose result is $x_2 \cdots x_m \, x_1$ [respectively $x_m \, x_1 \cdots x_{m-1}$].

The g. corollas of $\sigma$ and those of $\sigma'$ coincide except that $d$ takes the place of $c$ in $\sigma'$; that is, the multiplicity of $d$ in the multiset of g. corollas of $\sigma'$ is greater by one than that in $\sigma$ (consequently the multiplicity of $c$ is less by one in $\sigma'$ than in $\sigma$). This implies that $A(\sigma')=A(\sigma)$ because $\eta(c)=\eta(d)$ by Remark \ref{eta(psi)}. Moreover if $\tau$ is an SLA in $\overline{R}$ computing $c$ and such that $\eta(\tau)=\eta(c)$ then if we add to $\tau$ a step equal to the cyclic conjugation of $c$ (which gives $d$), we obtain an SLA (which we call $\tau'$) computing $d$, such that $\eta(d)=\eta(\tau')$ and such that the multiset of base elements of $\tau$ coincides with that of $\tau'$. This implies that a CMDR for $\sigma$ is a CMDR also for $\sigma'$.
     \end{proof}

\chapter{The Main Theorem} \label{tmt}

In this chapter we prove the Main Theorem of the thesis. The proof is split in two parts, the first part in Section \ref{casebycase} in a specific case and the second part in Section \ref{conc} in the general case.

\section{The proof of the Main Theorem: a case by case analysis} \label{casebycase}

Let $X$ be a set of letters, let $X^{-1}$ be the set of inverses of elements of $X$ and let $R$ be a set of cyclically reduced non-empty words in $X\cup X^{-1}$ such that $R^{-1}\subset R$. Let $\mathcal{M}(X\cup X^{-1})$ be the free monoid on $X\cup X^{-1}$ and let $L_g$ (Definition \ref{maiobj}) be the subset of $\mathcal{M}(X\cup X^{-1})$ recursively defined by g. corollas and stems and by the operation of insertion. Let $\mathcal{F}(X)$ be the free group on $X$ and let $\mathcal{N}$ be the normal closure of $R$ in $\mathcal{F}(X)$; in particular $\mathcal{N}$ is the set of (reduced) relators of the group presentation $\langle \, X \, | \, R \, \rangle$. Let the functions Area and $A$ as in Definitions \ref{defarea} and \ref{A}. Complete multisets of defining relators have been introduced in Definition \ref{CMDR}.

In this section and in the next we will show that the following result holds

\begin{theorem} \label{maithe} Let $\langle \, X \, | \, R \, \rangle$ be a group presentation and let $\mathcal{N}$ be the set of reduced relators. Then $\mathcal{N}$ coincides with the subset of $L$ consisting of reduced words. Let $w$ be the reduced form of $f_1 r_1 f_1^{-1} \cdots f_n r_n f_n^{-1}$, where $r_i\in R$; then there exist a submultiset $M$ of $\{r_1, \cdots, r_n\}$ and a straight line algorithm $\sigma$ computing $w$ which has $M$ as a complete multiset of defining relators (CMDR) and such that $A(\sigma)\leqslant n$. If $n=\textrm{Area}(w)$ then $A(\sigma)=n$, every corolla of $\sigma$ is a proper corolla and $\{r_1, \cdots, r_n\}$ is a CMDR for $\sigma$. Finally $A(w)=\textrm{Area}(w)$.   
       \end{theorem}

The equality $\textrm{Area}(w)=A(w)$ is very interesting because it gives an alternative way to define the area of a relator and therefore the Dehn function of a presentation.

Theorem \ref{maithe} implies the following

\begin{corollary} \label{maicor} Let $w \in \mathcal{N}$ and let $\sigma$ be an SLA computing $w$ and such that $A(\sigma)=A(w)$. Then the area of $w$ is equal to the sum of the areas of the g. corollas of $\sigma$, that is if $c_1, \cdots, c_m$ are the g. corollas of $\sigma$ then $\textrm{Area}(w)=\displaystyle \sum_{i=1}^m \textrm{Area}(c_i)$ and $|w|\geqslant \sum_{i=1}^m |c_i|$. \end{corollary}

We have seen in Remark \ref{alcap} that any reduced element of $L_g$ has a contiguous subword equal to a proper g. corolla. By Definition \ref{propcor} and Theorem \ref{maithe} we have that if $c$ is a proper g. corolla then $\eta(c)=\textrm{Area}(c)$. If $m=\textrm{Area}(c)$ then $c\in \overline{R}_m$; if $R$ is finite then $\overline{R}_m$ is finite by Theorem \ref{algo}. Thus we have

\begin{corollary} \label{cont} A relator of area $m$ has a contiguous subword equal to a proper g. corolla whose area is less or equal to $m$. If the presentation is finite then there are finitely many of such proper g. corollas.
     \end{corollary}

In Section \ref{statres} we have given the proof of 

\begin{corollary} \label{hyper2} The presentation $\langle \, X \, | \, R \, \rangle$ is hyperbolic (Definition \ref{hyper}) if and only if there exists a positive real constant $\alpha$ such that $\textrm{Area}(c) \leqslant \alpha |c|$ for every proper g. corolla $c$.
   \end{corollary}
   
   Corollary \ref{hyper2} says that to verify if a group is hyperbolic it is sufficient to verify the inequality $\textrm{Area}(w) \leqslant \alpha |w|$ only on proper g. corollas instead of all relators. This is a very interesting result because the set of g. corollas is a proper subset of the set of all relators, in particular it contains only cyclically reduced words.

For every $n$ let $\Delta_0(n)$ be the maximal area of proper g. corollas of length at most $n$; obviously $\Delta_0(n)\leqslant \Delta(n)$. In Section \ref{statres} we have proved

\begin{corollary} \label{hyper3} Let $\alpha$ be a positive real number; the Dehn function $\Delta$ is bounded by the linear function $\alpha n$ if and only if the $\Delta_0$ is. \end{corollary}

As we have seen in Section \ref{statres}, to prove the Main Theorem \ref{maithe} it is sufficient proving the following

\begin{lemma} \label{maile} Let $l:=l_1 z z^{-1} l_2$ (where $l_1$ and $l_2$ are words and $z$ a letter) be an element of $L_g$; then $l_1 l_2\in L_g$. In particular if $\sigma$ is a straight line algorithm computing $l$ and if $M$ is a CMDR for $\sigma$ (Definition \ref{CMDR}), then there exists a straight line algorithm $\sigma'$ computing $l_1 l_2$, such that $A(\sigma')\leqslant A(\sigma)$ and such that if $M'$ is a CMDR for $\sigma'$ then $M'\subset M$. Moreover if $A(\sigma)=\textrm{Area}\big(\rho(l)\big)$ then $A(\sigma')=A(\sigma)$ and $M'=M$.
                \end{lemma}
                
Lemma \ref{maile} is proved in this and in the following sections. Lemma \ref{maile} implies 
 
\begin{lemma} \label{maile2} \begin{enumerate}

  \item Let $l \in L$ and let $w \in \mathcal{M}(X\cup X^{-1})$ be such that $\rho(w)=\rho(l)$. Then $w \in L$.
  
  \item $L$ contains the cyclic conjugate of any of its elements.

\end{enumerate} 
     \end{lemma}
     
Finally Lemma \ref{maile2} implies the following interesting result

\begin{theorem} \label{maithe2} $L$ is the subset of $\mathcal{M}(X\cup X^{-1})$ of words whose reduced form belongs to $\mathcal{N}$, i.e., if $\rho:\mathcal{M}(X\cup X^{-1}) \rightarrow \mathcal{F}(X)$ is the function “reduced form" (Definition \ref{rho}) then $L=\rho^{-1}(\mathcal{N})$.       \end{theorem}

Lemma \ref{maile2} and Theorem \ref{maithe2} have been proved in Section \ref{statres} having assumed that Lemma \ref{maile} is true. Theorem \ref{maithe2} says in particular that $L_g$ is the set of all relators (reduced and non). In Section \ref{statres} we have proved the following

\begin{corollary} \label{vKlemma} The set of relators (reduced and non) of the presentation $\langle \, X \, | \, R \, \rangle$ coincides with the set of labels of the elements of $\mathcal{L}$ (Definition \ref{mathcalL}). \end{corollary}

In this section we will prove Lemma \ref{maile} under the following hypothesis: \textbf{the only output of $\sigma$ containing the subword $z z^{-1}$ of $l_1 z z^{-1} l_2$ is the last one}. We fix the notation until the end of Subsection \ref{III}: $l$ will denote the result of $\sigma$; $t:=x_1 \cdots x_m$ and $t':=y_1 \cdots y_p$ will denote the steps of $\sigma$ such that $l$ is the insertion of $t'$ into $t$; $\tau$ and $\tau'$ will denote the pSLsA's of $\sigma$ computing $t$ and $t'$ respectively. By Remark \ref{union} we have that $A(\sigma)=A(\tau)+A(\tau')$.

By the hypothesis assumed for this section, the letters $z$ and $z^{-1}$ of $z z^{-1}$ do not come both from $t$ or $t'$; therefore since $l$ contains $z z^{-1}$ and since $l$ is the insertion of $t'$ into $t$, then $z$ comes from $t$ and $z^{-1}$ from $t'$ or vice versa $z$ comes from $t'$ and $z^{-1}$ from $t$. Furthermore the insertion of $t'$ into $t$ makes $z$ and $z^{-1}$ consecutive. This means that there exists $n: 1\leqslant n \leqslant m$ such that $l=x_1 \cdots x_n \, y_1 \cdots y_p \, x_{n+1} \cdots x_m$ and: \begin{itemize}

  \item either $x_n=z$ and $y_1=z^{-1}$ (we call it \textit{subcase $\alpha$});
  
  \item or $y_p=z$ and $x_{n+1}=z^{-1}$ (we call it \textit{subcase $\beta$}).
        
        \end{itemize}

Until the end of the Subsection \ref{III} we also denote $s$ and $s'$ the base steps of $\sigma$ such that one of them contains the letter $z$ of $z z^{-1}$ and the other one contains $z^{-1}$ and such that $t$ depends on $s$, $t'$ depends on $s'$. In the subcase $\alpha$, $s$ contains $z$ and $s'$ contains $z^{-1}$; in the subcase $\beta$, $s'$ contains $z$ and $s$ contains $z^{-1}$. We can assume that $s\neq s'$ because $s=s'$ implies that $s$ contains $z z^{-1}$, therefore $s$ is the last step of $\sigma$ in view of our hypothesis. Since $s$ is a base step, it cannot use preceding steps and this means that $s$ is the only step of $\sigma$. $s$ cannot be a g. corolla because it contains $z z^{-1}$ as a subword and g. corollas are reduced. Indeed $s$ would be a stem and $z$ would be the last letter of its first half ($z^{-1}$ would be the first letter of the second half.) For this situation Lemma \ref{maile} is trivially true.
 
Four cases are then possible: I) $s$ and $s'$ are stems; II) $s$ is a stem and $s'$ a g. corolla; III) $s$ is a g. corolla and $s'$ a stem; IV) $s$ and $s'$ are g. corollas.

Let $s$ be a stem (Cases I and II); since $s$ is the product of a word by its inverse and since $s$ contains either the letter $z$ or the letter $z^{-1}$ of $z z^{-1}$, then two of its opposite letters (and therefore two letters of $t$) are equal to $z$ and $z^{-1}$, one (and only one) of which is of the subword $zz^{-1}$ of $l_1 zz^{-1} l_2$. We let $h, h': 1\leqslant h < h' \leqslant m$ be such that $x_h$ and $x_{h'}$ are those letters of $t$, that is $\{x_h, x_{h'}\}=\{z, z^{-1}\}$. We call \textit{subcase $1$} when $x_h=z$ and $x_{h'}=z^{-1}$, we call \textit{subcase 2} when $x_h=z^{-1}$ and $x_{h'}=z$. Therefore in the subcase $1\alpha$ we have $n=h$ and $x_{h'}=y_1=z^{-1}$; in the subcase $1\beta$ we have $n+1=h'$ and $x_h=y_p=z$; in the subcase $2\alpha$ we have $n=h'$ and $x_h=y_1=z^{-1}$; in the subcase $2\beta$ we have $n+1=h$ and $x_{h'}=y_p=z$.

Finally we let $j, j': 1 \leqslant j < j' \leqslant p$ be such that the letters $y_j$ and $y_{j'}$ of $t'$ are equal respectively to the first and the last letter of $s'$. In the subcase $\alpha$ we have $j=1$ and therefore $y_1=z^{-1}$; in the subcase $\beta$ we have $j'=p$ and $y_p=z$. If $s'$ is a stem (Cases I and III) then in the subcase $\alpha$ we have $y_{j'}=z$, in the subcase $\beta$ we have $y_j=z^{-1}$.

 We recall that we have denoted $\tau$ the proper straight line subalgorithm of $\sigma$ computing $t$.

\begin{lemma} \label{cl2} Let $s$ be a stem (Cases I and II) and let $v_1, v_2$ and $v$ be the following subwords of $t$: 
  $$v_1=x_1 \cdots x_{h-1}, \, \, \, \, \, \, \, \, v_2=x_{h'+1} \cdots x_m, \, \, \, \,\, \, \, \,  v=x_{h+1} \cdots x_{h'-1}.$$
 There exist two SLA's of $\sigma$, denoted $\sigma_1$ and $\sigma_2$, which compute $v_1 v_2$ and $v$ respectively and such that $A(\sigma_1)+A(\sigma_2)=A(\tau)$. Moreover the union of a CMDR for $\sigma_1$ and of one of $\sigma_2$ is a CMDR for $\tau$. 
           \end{lemma}

\begin{proof}  Follows from Lemma \ref{oso}.
                                                                         \end{proof}

 We recall that we have denoted $\tau'$ the proper straight line subalgorithm of $\sigma$ computing $t'$.

\begin{lemma} \label{cl3} Let $s'$ be a stem (Cases I and III) and let $w_1, w_2, w'_1$ and $w'_2$ be the following subwords of $t'$: 
  $$w_1=y_1 \cdots y_{j-1}, \, \, \, \, \, \,  w_2=y_{j+1} \cdots y_{p-1},   \, \, \, \, \, \, w'_1=y_2 \cdots y_{j'-1},  \, \, \, \, \, \, w'_2=y_{j'+1} \cdots y_p.$$
  There exist two SLA's of $\sigma$, denoted $\tau_1$ and $\tau_2$, such that: \begin{enumerate}

  \item in the subcase $\alpha$, $\tau_1$ and $\tau_2$ compute $w'_1$ and $w'_2$ respectively and $A(\tau_1)+A(\tau_2)=A(\tau')$;
  
  \item in the subcase $\beta$, $\tau_1$ and $\tau_2$ compute $w_1$ and $w_2$ respectively and $A(\tau_1)+A(\tau_2)=A(\tau')$.

\end{enumerate}  

Finally the union of a CMDR for $\tau_1$ and of one for $\tau_2$ is a CMDR for $\tau$. 

\end{lemma}

\begin{proof} \begin{enumerate}

  \item Follows from Lemma \ref{oso} because in the subcase $\alpha$, $y_1=z^{-1}$ and $y_{j'}=z$.

  \item Follows from Lemma \ref{oso} because in the subcase $\beta$, $y_j=z^{-1}$ and $y_p=z$. 
  
  \end{enumerate}
                 \end{proof}
                 
 \begin{remark} \label{ora} \rm Let $\sigma'$ be an SLA computing $l_1 l_2$ and such that $A(\sigma') \leqslant A(\sigma)$ and let $A(\sigma)=\textrm{Area}\big(\rho(l)\big)$. We have that $\textrm{Area}\big(\rho(l)\big)=\textrm{Area}\big(\rho(l_1 l_2)\big)$ since $\rho(l)=\rho(l_1 l_2)$ and Area$\big(\rho(l_1 l_2)\big) \leqslant A(\sigma')$ by Proposition \ref{impo}. These inequalities imply that $A(\sigma') = A(\sigma)$.
                \end{remark}

\subsection{Case I} $s$ and $s'$ are stems. As in Lemmas \ref{cl2} and \ref{cl3} we let $v_1=x_1\cdots x_{h-1}$, $v_2=x_{h'+1} \cdots x_m$, $v=x_{h+1} \cdots x_{h'-1}$, $w_1=y_1\cdots  y_{j-1}$, $w_2=y_{j+1} \cdots y_{p-1}$, $w'_1=y_2 \cdots  y_{j'-1}$ and $w'_2=y_{j'+1} \cdots y_p$.

\begin{remark} \label{cl4} \rm  By Lemmas \ref{cl2} and \ref{cl3} there exist SLA's $\sigma_1$ and $\sigma_2$ computing $v_1 v_2$ and $v$ and SLA's $\tau_1$ and $\tau_2$ computing $w'_1$ and $w'_2$ in the subcase $\alpha$, $w_1$ and $w_2$ in the subcase $\beta$, such that $A(\tau)=A(\sigma_1)+A(\sigma_2)$ and $A(\tau')=A(\tau_1)+A(\tau_2)$. Furthermore, since $A(\tau)+A(\tau')=A(\sigma)$ by Remark \ref{union}, then
 $$A(\sigma_1)+A(\sigma_2)+A(\tau_1)+A(\tau_2)=A(\sigma).$$
Finally the union of CMDR's for $\sigma_1, \sigma_2, \tau_1$ and $\tau_2$ is a CMDR for $\sigma$.
      \end{remark} 

\smallskip  \smallskip  \smallskip

\textbf{Subcase 1$\alpha$.} We have $n=h$, $x_h=y_{j'}=z$ and $x_{h'}=y_1=z^{-1}$.
Therefore
      $$l=x_1\cdots x_{h-1} \, z\,( z^{-1} y_2 \cdots  y_{j'-1} \, z \, y_{j'+1} \cdots y_p) \, x_{h+1} \cdots x_{h'-1} \, z^{-1} x_{h'+1} \cdots x_m=$$
      $$v_1 z ( z^{-1} w'_1 \, z \, w'_2) \, v \, z^{-1} v_2$$
and $l_1= v_1$, $l_2=w'_1 \, z \, w'_2 \, v \, z^{-1} v_2$.     
   By Remark \ref{cl4}, Proposition \ref{mesc} and Part 2 of Corollary \ref{promesc} there exists an SLA $\sigma'$ computing 
           $$v_1 \, w'_1 \, z \, w'_2 \, v \, z^{-1} v_2=l_1 l_2$$
such that $A(\sigma')=A(\sigma_1)+A(\sigma_2)+A(\tau_1)+A(\tau_2)=A(\sigma)$ and such that a CMDR for $\sigma$ is a CMDR also for $\sigma'$.

\smallskip  \smallskip  \smallskip

\textbf{Subcase 1$\beta$.} We have that $n+1=h'$, that $x_h=y_p=z$ and that $x_{h'}=y_j=z^{-1}$.
Therefore
    $$l=x_1\cdots x_{h-1} \, z  \, x_{h+1} \cdots x_{h'-1} (y_1\cdots  y_{j-1} \, z^{-1} y_{j+1} \cdots y_{p-1} \, z) \, z^{-1} \,  x_{h'+1} \cdots x_m=$$
       $$v_1 \, z  \, v (w_1 \, z^{-1} w_2 \, z) \, z^{-1} \,  v_2$$
and $l_1= v_1 \, z  \, v w_1 \, z^{-1} w_2$, $l_2=v_2$.
        By Remark \ref{cl4}, Proposition \ref{mesc} and Part 2 of Corollary \ref{promesc} there exists an SLA $\sigma'$ computing
         $$v_1 \, z \, v \, w_1 \, z^{-1} \, w_2 \, v_2=l_1 l_2$$
such that $A(\sigma')=A(\sigma_1)+A(\sigma_2)+A(\tau_1)+A(\tau_2)=A(\sigma)$ and such that a CMDR for $\sigma$ is a CMDR also for $\sigma'$.

\smallskip  \smallskip  \smallskip

\textbf{Subcase 2$\alpha$.} We have that $n=h'$, that $x_{h'}=y_{j'}=z$ and that $x_h=y_1=z^{-1}$.
Therefore
          $$l=x_1\cdots x_{h-1} \, z^{-1} \, x_{h+1} \cdots x_{h'-1} \, z \, (z^{-1}\, y_2\cdots  y_{j'-1} \, z \, y_{j'+1} \cdots y_p) \, x_{h'+1} \cdots x_m=$$
          $$v_1 \, z^{-1} \, v \, z \, (z^{-1}\, w'_1 \, z \, w'_2) \, v_2$$
and $l_1=v_1 \, z^{-1} \, v$, $l_2=w'_1 \, z \, w'_2 \, v_2$.
   By Remark \ref{cl4}, Proposition \ref{mesc} and Part 2 of Corollary \ref{promesc} there exists an SLA $\sigma'$ computing     
       $$v_1 \, z^{-1} \, v \, w'_1 \, z \, w'_2 \, v_2=l_1 l_2$$
such that $A(\sigma')=A(\sigma_1)+A(\sigma_2)+A(\tau_1)+A(\tau_2)=A(\sigma)$ and such that a CMDR for $\sigma$ is a CMDR also for $\sigma'$.

\smallskip  \smallskip  \smallskip

\textbf{Subcase 2$\beta$.} We have that $n+1=h$, that $x_{h'}=y_p=z$ and that $x_h=y_j=z^{-1}$.
         Therefore
$$l=x_1\cdots x_{h-1} \, (y_1\cdots  y_{j-1} \, z^{-1} \, y_{j+1} \cdots y_{p-1} \, z) \, z^{-1} \, x_{h+1} \cdots x_{h'-1} \, z \, x_{h'+1} \cdots x_m=$$
    $$v_1 \,  w_1 \, z^{-1} \, w_2 \, z \, z^{-1} \, v \, z \, v_2$$
and $l_1=v_1 \,  w_1 \, z^{-1} \, w_2$, $l_2=v \, z \, v_2$.
  By Remark \ref{cl4}, Proposition \ref{mesc} and Part 2 of Corollary \ref{promesc} there exists an SLA $\sigma'$ computing     
      $$v_1 \,  w_1 \, z^{-1} \, w_2 \, v \, z \, v_2=l_1 l_2$$
such that $A(\sigma')=A(\sigma_1)+A(\sigma_2)+A(\tau_1)+A(\tau_2)=A(\sigma)$ and such that a CMDR for $\sigma$ is a CMDR also for $\sigma'$.

\subsection{Case II}

$s$ is a stem and $s'$ a g. corolla. As in Lemmas \ref{cl2} and \ref{cl3} we let $v_1=x_1\cdots x_{h-1}$, $v_2=x_{h'+1} \cdots x_m$, $v=x_{h+1} \cdots x_{h'-1}$, $w_1=y_1\cdots  y_{j-1}$, $w_2=y_{j+1} \cdots y_{p-1}$, $w'_1=y_2 \cdots  y_{j'-1}$ and $w'_2=y_{j'+1} \cdots y_p$.

\begin{remark} \label{anc} \rm Set $u:=y_1\cdots y_{p-1}$ and $u':=y_2 \cdots y_p$. In the subcase $\alpha$ we have $t'=z^{-1} u'$ and the first letter of $t'$ comes from $s'$; in the subcase $\beta$ we have $t'= u \, z$ and the last letter of $t'$ comes from $s'$. By Lemma \ref{faso} there exists an SLA $\tau'_1$ computing $u'  z^{-1}$ in the subcase $\alpha$, computing $z \, u$ in the subcase $\beta$ such that $A(\tau'_1)=A(\tau')$ and a CMDR for $\tau'$ is a CMDR also for $\tau'_1$.

By Lemma \ref{cl2} there exist SLA's $\sigma_1$ and $\sigma_2$ computing $v_1 v_2$ and $v$ and such that $A(\sigma_1)+A(\sigma_2)=A(\tau)$. Furthermore, since $A(\tau)+A(\tau')=A(\sigma)$ by Remark \ref{union}, then
 $$A(\sigma_1)+A(\sigma_2)+A(\tau'_1)=A(\sigma).$$
Finally the union of CMDR's for $\sigma_1, \sigma_2$ and $\tau'_1$ is a CMDR for $\sigma$.
          \end{remark}

\textbf{Subcase 1$\alpha$.} We have that $n=h$, that $x_h=z$ and that $x_{h'}=y_1=z^{-1}$.
Therefore
$$l=x_1\cdots x_{h-1} \, z \, (z^{-1} \, y_2\cdots y_p) \, x_{h+1} \cdots x_{h'-1} \, z^{-1} \, x_{h'+1} \cdots x_m=$$
   $$v_1 \, z \, (z^{-1} \, u') \, v \, z^{-1} \, v_2$$
and $l_1=v_1$, $l_2=u' \, v \, z^{-1} \, v_2$.

By Remark \ref{anc} and Proposition \ref{mesc} there exists an SLA $\sigma'$ computing 
                $$v_1 \, u' \, v \, z^{-1} \, v_2=l_1 l_2$$
such that $A(\sigma')=A(\sigma_1)+A(\sigma_2)+A(\tau'_1)=A(\sigma)$ and a CMDR for $\sigma$ is a CMDR also for $\sigma'$. 

\smallskip  \smallskip  \smallskip

\textbf{Subcase 1$\beta$.} We have that $n+1=h'$, that $x_h=y_p=z$ and that $x_{h'}=z^{-1}$. Therefore
          $$l=x_1\cdots x_{h-1} \, z \, x_{h+1} \cdots x_{h'-1} \, (y_1\cdots y_{p-1} \, z )\, z^{-1} \, x_{h'+1} \cdots x_m=$$
          $$v_1 \, z \, v \, (u \, z )\, z^{-1} \, v_2$$
and $l_1=v_1 \, z \, v \, u$, $l_2=v_2$.          
          
By Remark \ref{anc} and Proposition \ref{mesc} there exists an SLA $\sigma'$ computing  
                $$v_1 \, z \, v \, u \, v_2=l_1 l_2$$
such that $A(\sigma')=A(\sigma_1)+A(\sigma_2)+A(\tau'_1)=A(\sigma)$ and a CMDR for $\sigma$ is a CMDR also for $\sigma'$.

\smallskip  \smallskip   \smallskip

\textbf{Subcase 2$\alpha$.} We have that $n=h'$, that $x_{h'}=z$ and that $x_h=y_1=z^{-1}$. Therefore
        $$l=x_1\cdots x_{h-1} \, z^{-1} \, x_{h+1} \cdots x_{h'-1} \, z \, (z^{-1} \, y_2\cdots y_p) \, x_{h'+1} \cdots x_m=$$
        $$v_1 \, z^{-1} \, v \, z \, (z^{-1} \, u') \, v_2$$
and $l_1=v_1 \, z^{-1} \, v$, $l_2=u' \, v_2$.        
        By Remark \ref{anc} and Proposition \ref{mesc} there exists an SLA $\sigma'$ computing  
                $$v_1 \, z^{-1} \, v \, u' \, v_2=l_1 l_2$$
such that $A(\sigma')=A(\sigma_1)+A(\sigma_2)+A(\tau'_1)=A(\sigma)$ and a CMDR for $\sigma$ is a CMDR also for $\sigma'$.

\smallskip  \smallskip  \smallskip

\textbf{Subcase 2$\beta$.} We have that $n+1=h$, that $x_{h'}=y_p=z$ and that $x_h=z^{-1}$. Therefore
           $$l=x_1\cdots x_{h-1} \, (y_1\cdots y_{p-1} \, z) \, z^{-1} \, x_{h+1} \cdots x_{h'-1} \, z \, x_{h'+1} \cdots x_m=$$
           $$v_1 \, (u \, z) \, z^{-1} \, v \, z \, v_2$$
and $l_1=v_1 \, u$, $l_2=v \, z \, v_2$.           
           
By Remark \ref{anc} and Proposition \ref{mesc} there exists an SLA $\sigma'$ computing  
                $$v_1 \, u \, v \, z \, v_2=l_1 l_2$$
such that $A(\sigma')=A(\sigma_1)+A(\sigma_2)+A(\tau'_1)=A(\sigma)$ and a CMDR for $\sigma$ is a CMDR also for $\sigma'$.

\subsection{Case III}   \label{III}
$s$ is a g. corolla and $s'$ a stem. 

\smallskip  \smallskip

\textbf{Subcase $\alpha$.} We have that $y_1=z^{-1}$ and that $x_n=y_{j'}=z$. Set $u_1:=x_1\cdots x_{n-1}$ and $u_2:=x_{n+1} \cdots  x_m$; thus $t=u_1  z u_2$ and as said at the beginning of the section, $t$ is computed by the pSLsA $\tau$. As in Lemma \ref{cl3}, let $w'_1=y_2 \cdots y_{j'-1}$ and $w'_2=y_{j'+1} \cdots y_p$. Therefore
     $$l=x_1\cdots x_{n-1} \, z \, (z^{-1} \, y_2\cdots  y_{j'-1} \, z \, y_{j'+1} \cdots y_p) \, x_{n+1} \cdots  x_m=$$
        $$u_1 \, z \, (z^{-1} \, w'_1 \, z \, w'_2) \, u_2$$
and $l_1=u_1$, $l_2=w'_1\, z \, w'_2 \, u_2$.        
By Lemma \ref{cl3} and by Proposition \ref{mesc} there exists an SLA $\sigma'$ computing
                $$u_1 \, w'_1 \, z \, w'_2 \, u_2=l_1 l_2$$  
and such that $A(\sigma')=A(\tau)+A(\tau')$. By Remark \ref{union}, $A(\tau)+A(\tau')=A(\sigma)$, thus $A(\sigma')=A(\sigma)$, and a CMDR for $\sigma$ is a CMDR also for $\sigma'$.

\smallskip  \smallskip  \smallskip

\textbf{Subcase $\beta$.} We have that $x_{n+1}=y_j=z^{-1}$ and that $y_p=z$. Set $u_1:=x_1\cdots x_n$ and $u_2:=x_{n+2} \cdots  x_m$; thus $t=u_1 \, z^{-1} \, u_2$ and as said at the beginning of the section, $t$ is computed by the pSLsA $\tau$. As in Lemma \ref{cl3}, let $w_1=y_1 \cdots y_{j-1}$ and $w_2=y_{j+1} \cdots y_{p-1}$. Therefore
      $$l=x_1\cdots x_n \, (y_1\cdots  y_{j-1} \, z^{-1} \, y_{j+1} \cdots y_{p-1} \, z) \, z^{-1} \, x_{n+2} \cdots x_m=$$
      $$u_1 \, (w_1 \, z^{-1} \, w_2 \, z) \, z^{-1} \, u_2$$
and $l_1=u_1 \, w_1 \, z^{-1} \, w_2$, $l_2=u_2$.        
By Lemma \ref{cl3} and by Proposition \ref{mesc} there exists an SLA $\sigma'$ computing
          $$u_1 \,  w_1 \, z^{-1} \, w_2 \, u_2=l_1 l_2$$
and such that $A(\sigma')=A(\tau)+A(\tau')$. By Remark \ref{union}, $A(\tau)+A(\tau')=A(\sigma)$, thus $A(\sigma')=A(\sigma)$, and a CMDR for $\sigma$ is a CMDR also for $\sigma'$.

\subsection{Case IV}  \label{IV}

There exist two g. corollas $c$ and $c'$ such that one of them contains the letter $z$ of $z z^{-1}$ and the other one contains $z^{-1}$. Let $l:=z_1 \cdots z_e$ be the result of $\sigma$, let $c:=x_1 \cdots x_m$, $c':=y_1 \cdots y_p$ and let $h_1, \cdots, h_m, i_1, \cdots, i_p$ be indices such that 
   $$z_{h_1}=x_1, \cdots, \, z_{h_m}=x_m,  \, z_{i_1}=y_1, \cdots, \, z_{i_p}=y_p.$$ 
 There are two possibilities: either there exists $k: 1\leqslant k \leqslant m$ such that $x_k=z$, $y_1=z^{-1}$ and $h_k+1=i_1$ (we call it \textit{subcase $\alpha$}) or there exists $k: 1 < k \leqslant m$ such that $y_p=z$, $x_k=z^{-1}$ and $i_p+1=h_k$ (we call it \textit{subcase $\beta$})\footnote{the case $y_p=z$ and $x_1=z^{-1}$ is analogous to the case $x_m=z$ and $y_1=z^{-1}$.}. We call \textit{subcase $\alpha_1$} the subcase $\alpha$ with $1\leqslant k<m$, \textit{subcase $\alpha_2$} the subcase $\alpha$ with $k=m$.

Let $r_1, \cdots, r_{m-1}$ be the ramifications (Definition \ref{ramsur}) from $c$ at $x_1$, $\cdots$, $x_{m-1}$ respectively; let $r'_1, \cdots,$ $r'_{p-1}$ be the ones from $c'$ at $y_1, \cdots, y_{p-1}$. Let $q_0$ be the preceding of $c$ and $q_1$ its following (Definition \ref{ramsur}), that is $q_0 q_1$ is the surround of $c$.

\smallskip  \smallskip  \smallskip

\textbf{Subcase $\alpha_1$.} We have that $x_k=z$ and $y_1=z^{-1}$. Since $h_k+1=i_1$ then 
                      $$h_1<\cdots <h_k<i_1<\cdots<i_p<h_{k+1}<\cdots<h_m$$
   and
          $$l=z_1 \cdots z_{h_1} \cdots z_{h_k} \, (z_{i_1} \cdots z_{i_p}) \, z_{i_p+1} \cdots z_{h_{k+1}-1} \, z_{h_{k+1}} \cdots z_{h_m} \cdots z_e.$$
   By Proposition \ref{consec} we can suppose that $\sigma$ defines consecutively the insertions into $c'$. Since $z_{i_1}=y_1$ and $z_{i_p}=y_p$, then the surround of $c'$ is
         $$l':=z_1 \cdots z_{h_k} \, z_{i_p+1} \cdots z_{h_{k+1}-1} \, z_{h_{k+1}} \cdots z_e.$$    
     By Proposition \ref{fam} there exists an SLsA of $\sigma$, denoted $\tau$, which computes $l'$ and $c$ is one of its g. corollas. By Proposition \ref{consec} we can suppose that $\tau$ defines consecutively the insertions into $c$. Since $z_{h_k}=x_k$ and $z_{h_{k+1}}=x_{k+1}$, then the ramification of $l'$ from $c$ at $x_k$ is $v:=z_{i_p+1} \cdots z_{h_{k+1}-1}$ and by Proposition \ref{fam} there exists an SLsA of $\tau$ (and therefore of $\sigma$) which computes it.

We have that:
$$z_1 \cdots z_{h_k-1}=u,\, \, \, \textrm{where} \, \, \,  u=q_0 \, \textbf{\emph{x}}_\textbf{1} \, r_1 \cdots \textbf{\emph{x}}_\textbf{{k-1}} \, r_{k-1}; \, \, \, \, \, \,  \, \, \, z_{h_k}=\textbf{\emph{z}};$$
   $$z_{i_1}=\textbf{\emph{z}}^\textbf{-1}; \, \, \,  \, \, \, \, \, \,  \, \, \, z_{i_1+1} \cdots z_{i_p}=u', \, \, \, \textrm{where} \, \, \,  u'=r'_1 \, \textbf{\emph{y}}_\textbf{2} \cdots r'_{p-1} \, \textbf{\emph{y}}_\textbf{p};$$
   $$z_{h_{k+1}} \cdots z_e=u'', \, \, \, \textrm{where} \, \, \, u''=\textbf{\emph{x}}_\textbf{k+1} \, r'_{k+1} \cdots r'_{m-1} \, \textbf{\emph{x}}_\textbf{m} \, q_1$$
 (we write in bold the letters coming from $c$ and $c'$).
   
Thus 
                                        $$l=u \, z \, (z^{-1} \, u') \, v \, u''=$$
 $$q_0 \, \textbf{\emph{x}}_\textbf{1} \, r_1 \cdots \textbf{\emph{x}}_\textbf{{k-1}} \, r_{k-1} \, \textbf{\emph{z}} \, (\textbf{\emph{z}}^\textbf{-1} \, r'_1 \, \textbf{\emph{y}}_\textbf{2} \cdots r'_{p-1} \, \textbf{\emph{y}}_\textbf{p}) \, v \, \textbf{\emph{x}}_\textbf{k+1} \, r'_{k+1} \cdots r'_{m-1} \, \textbf{\emph{x}}_\textbf{m} \, q_1$$
   and $l_1=u$, $l_2=u' v  u''$.

We have that $c=x_1 \cdots x_{k-1} \, z \, x_{k+1} \cdots x_m$ and $c'=z^{-1} \, y_2 \cdots y_p$.

As in the proof of Proposition \ref{bigg. corolla} we set 
$$c_1:=x_1 \cdots x_{k-1}, \, \, \, \, \, c_2:=x_{k+1} \cdots x_m, \, \, \, \, \, d':=\pi(c_2 c_1, c')$$
  $$f_1:=x_1 \cdots x_{k-1}, \, \, \, \, \, f_2:=y_2 \cdots y_p \, x_{k+1} \cdots x_m.$$ 
By Proposition \ref{bigg. corolla}, either $f_1 f_2$ is a stem element or an insertion of two stem elements (possibly empty) into a flower element with g. corolla a cyclic conjugate of $d'$, which we denote $d''$; in this case $\eta(d'')\leqslant \eta(c)+\eta(c')$.

We modify $\sigma$ by replacing $c$ and $c'$ with an SLA $\tau'$ defined in the following way: if $f_1 f_2$ is a stem element then $\tau'$ has only one step which is equal to $f_1 f_2$; otherwise $\tau'$ is an SLA computing $f_1 f_2$ and having only one g. corolla, equal to $d''$, therefore $A(\tau')=\eta(d'')$. We consider the insertions at a letter of $c$ or $c'$ as insertions at the same letter of $f_1 f_2$. An insertion at $z^{-1}$ is replaced by an insertion at the last letter that in the output of the same step of $\sigma$ was preceding $z$. With this modification we obtain an SLA $\sigma'$ computing $u \, u' \, v \, u''=l_1 l_2$. 

Let $N$ and $N'$ be the multisets of base elements of $\sigma$ and $\sigma'$ respectively. For any g. corolla of $\sigma$ let us fix an SLA in $\overline{R}$ computing that g. corolla and let $M$ be the union of the multisets of the base elements of the SLA's in $\overline{R}$ for the elements of $N$. $M$ is a CMDR for $\sigma$. If $f_1 f_2$ is a stem then $N'=N\setminus{\{c, c'\}}$ and thus $A(\sigma') < A(\sigma)$. Then $M$ minus the base elements of the chosen SLA's in $\overline{R}$ for $c$ and $c'$ is a CMDR for $\sigma'$. Let $f_1 f_2$ be not a stem; then $N'\setminus{\{d'\}}=N\setminus{\{c, c'\}}$, which means that the g. corollas of $\sigma$ and $\sigma'$ coincide except that $c$ and $c'$ are replaced in $\sigma'$ by $d''$; therefore the inequality $\eta(d'')\leqslant \eta(c)+\eta(c')$ implies that $A(\sigma')\leqslant A(\sigma)$. Moreover $M$ is a CMDR for $\sigma'$. 

If $A(\sigma)=\textrm{Area}\big(\rho(l)\big)$ then $A(\sigma')=A(\sigma)$ by Remark \ref{ora}, which implies by Remark \ref{dopo} that a CMDR for $\sigma$ is a CMDR also for $\sigma'$.

\smallskip  \smallskip  \smallskip

\textbf{Subcase $\beta$.} We have that $y_p=z$ and $x_k=z^{-1}$. Since $i_p+1=h_k$, then
               $$h_1<\cdots <h_{k-1}<i_1<\cdots<i_p<h_k<\cdots<h_m$$
and
      $$l=z_1 \cdots z_{h_1} \cdots z_{h_{k-1}} \cdots z_{i_1-1} \, (z_{i_1} \cdots z_{i_p}) \, z_{h_k} \cdots z_{h_m} \cdots z_e.$$
By Proposition \ref{consec} we can suppose that $\sigma$ defines consecutively the insertions into $c'$. Since $z_{i_1}=y_1$ and $z_{i_p}=y_p$, then the surround of $c'$ is 
   $$l'=z_1 \cdots z_{h_1} \cdots z_{h_{k-1}} \cdots z_{i_1-1}  \, z_{h_k} \cdots z_{h_m} \cdots z_e.$$
By Proposition \ref{fam} there exists an SLsA of $\sigma$, denoted $\tau$, which computes $l'$ and $c$ is one of its g. corollas. By Proposition \ref{consec} we can suppose that $\tau$ defines consecutively the insertions into $c$. Since $z_{h_{k-1}}=x_{k-1}$ and $z_{h_k}=x_k$, then the ramification of $l'$ from $c$ at $x_{k-1}$ is $v:=z_{h_{k-1}+1} \cdots z_{i_1-1}$ and by Proposition \ref{fam} there exists an SLsA of $\tau$ (and therefore of $\sigma$) which computes it.

We have that
        $$z_1 \cdots z_{i_1-1}=u, \, \, \, \textrm{where} \, \, \, u=q_0 \, \textbf{\emph{x}}_\textbf{1} \, r_1 \cdots \textbf{\emph{x}}_\textbf{{k-1}} \, v;$$
        $$z_{i_1} \cdots z_{i_p-1}=u', \, \, \, \textrm{where} \, \, \, u'=\textbf{\emph{y}}_\textbf{1} \, r'_1 \, \cdots \textbf{\emph{y}}_\textbf{{p-1}} \, r'_{p-1}; \, \, \,  \, \, \, \, \, \,  \, \, \, z_{i_p}=\textbf{\emph{z}}; \, \, \,  \, \, \, \, \, \,  \, \, \, z_{h_k}=\textbf{\emph{z}}^\textbf{-1};$$
        $$z_{h_{k+1}} \cdots z_e=u'',  \, \, \, \textrm{where} \, \, \, u''=\textbf{\emph{x}}_\textbf{k+1} \, r'_{k+1} \cdots r'_{m-1} \, \textbf{\emph{x}}_\textbf{m} \, q_1.$$
    Thus
                $$l=u \, v \, (u' \, z) \, z^{-1} \, u''=$$
$$q_0 \, \textbf{\emph{x}}_\textbf{1} \, r_1 \cdots \textbf{\emph{x}}_\textbf{{k-1}} \, v \, (\textbf{\emph{y}}_\textbf{1} \, r'_1 \, \cdots \textbf{\emph{y}}_\textbf{{p-1}} \, r'_{p-1} \,\textbf{\emph{z}}) \, \textbf{\emph{z}}^\textbf{-1} \,\textbf{\emph{x}}_\textbf{k+1} \, r'_{k+1} \cdots r'_{m-1} \, \textbf{\emph{x}}_\textbf{m} \, q_1$$
   and $l_1=u v u'$, $l_2=u''$.
   
We have that $c=x_1 \cdots x_{k-1} \, z^{-1} \, x_{k+1} \cdots x_m$ and $c'=y_1 \cdots y_{p-1} z$. As in the proof of Proposition \ref{bigg. corolla} we set 
$$c_1:=x_1 \cdots x_{k-1}, \, \, \, \, \, c_2:=x_{k+1} \cdots x_m, \, \, \, \, \, d':=\pi(c_2 c_1, c')$$
  $$f_1:=x_1 \cdots x_{k-1} y_1 \cdots y_{p-1}, \, \, \, \, \, f_2:=x_{k+1} \cdots x_m.$$ 
By Proposition \ref{bigg. corolla}, either $f_1 f_2$ is a stem element or an insertion of two stem elements (possibly empty) into a flower element with g. corolla a cyclic conjugate of $d'$, which we denote $d''$; in this case $\eta(d'')\leqslant \eta(c)+\eta(c')$.

We modify $\sigma$ by replacing $c$ and $c'$ with an SLA $\tau'$ defined in the following way: if $f_1 f_2$ is a stem element then $\tau'$ has only one step which is equal to $f_1 f_2$; otherwise $\tau'$ is an SLA computing $f_1 f_2$ and having only one g. corolla, equal to $d''$, therefore $A(\tau')=\eta(d'')$. We consider the insertions at a letter of $c$ or $c'$ as insertions at the same letter of $f_1 f_2$. An insertion at $z^{-1}$ is replaced by an insertion at the last letter that in the output of the same step of $\sigma$ was preceding $z$. With this modification we obtain an SLA $\sigma'$ computing $u v u' u''=l_1 l_2$. 

Let $N$ and $N'$ be the multisets of base elements of $\sigma$ and $\sigma'$ respectively. For any g. corolla of $\sigma$ let us fix an SLA in $\overline{R}$ computing that g. corolla and let $M$ be the union of the multisets of the base elements of the SLA's in $\overline{R}$ for the elements of $N$. $M$ is a CMDR for $\sigma$. If $f_1 f_2$ is a stem then $N'=N\setminus{\{c, c'\}}$ and thus $A(\sigma') < A(\sigma)$. Then $M$ minus the base elements of the chosen SLA's in $\overline{R}$ for $c$ and $c'$ is a CMDR for $\sigma'$. Let $f_1 f_2$ be not a stem; then $N'\setminus{\{d'\}}=N\setminus{\{c, c'\}}$, which means that the g. corollas of $\sigma$ and $\sigma'$ coincide except that $c$ and $c'$ are replaced in $\sigma'$ by $d''$; therefore the inequality $\eta(d'')\leqslant \eta(c)+\eta(c')$ implies that $A(\sigma')\leqslant A(\sigma)$. Moreover $M$ is a CMDR for $\sigma'$.

If $A(\sigma)=\textrm{Area}\big(\rho(l)\big)$ then $A(\sigma')=A(\sigma)$ by Remark \ref{ora}, which implies by Remark \ref{dopo} that a CMDR for $\sigma$ is a CMDR also for $\sigma'$.

\smallskip  \smallskip  \smallskip

\textbf{Subcase $\alpha_2$.} We have that $x_m=z$ and $y_1=z^{-1}$. Since $h_m+1=i_1$, then
                $$h_1<\cdots <h_m<i_1<\cdots<i_p$$
     and
  $$l=z_1 \cdots z_{h_1} \cdots z_{h_m} \, (z_{i_1} \cdots z_{i_p}) \, z_{i_p+1} \cdots z_e.$$
By Proposition \ref{consec} we can suppose that $\sigma$ defines consecutively the insertions into $c'$. Since $z_{h_1}=x_1$ and $z_{h_m}=x_m$, then the surround of $s$ is 
      $$l'=z_1 \cdots z_{h_1-1} \, z_{i_1} \cdots z_{i_p} \, z_{i_p+1} \cdots z_e.$$ 
   By Proposition \ref{fam} there exists an SLsA of $\sigma$, denoted $\tau$, which computes $l'$ and $c$ is one of its g. corollas. By Proposition \ref{consec} we can suppose that $\tau$ defines consecutively the insertions into $c$.

We have that
$$z_1 \cdots z_{h_m-1}=u, \, \, \, \textrm{where} \, \, \, u=q_0\, \textbf{\emph{x}}_\textbf{1} \, r_1 \, \cdots \textbf{\emph{x}}_\textbf{m-1} \, r_{m-1};  \, \, \, \, \, \, z_{h_m}=\textbf{\emph{z}}; \, \, \, \, \, \, z_{i_1}=\textbf{\emph{z}}^\textbf{-1};$$
    $$z_{i_1+1} \cdots z_{i_p}=u', \, \, \, \textrm{where} \, \, \, u'=r'_1 \, \textbf{\emph{y}}_\textbf{2} \cdots r'_{p-1} \, \textbf{\emph{y}}_\textbf{p}.$$
  Thus
                               $$l=u \, z \, (z^{-1} \, u') \, q_1=q_0 \, \textbf{\emph{x}}_\textbf{1} \, r_1 \, \cdots \textbf{\emph{x}}_\textbf{m-1} \, r_{m-1} \, \textbf{\emph{z}} \, (\textbf{\emph{z}}^\textbf{-1} \, r'_1 \, \textbf{\emph{y}}_\textbf{2} \cdots r'_{p-1} \, \textbf{\emph{y}}_\textbf{p}) \, q_1$$
  and $l_1=u$, $l_2=u' \, q_1$.

We have that $c=x_1 \cdots x_{m-1} \, z$ and $c'=z^{-1} \, y_2 \cdots y_p$. As in the proof of Proposition \ref{cc'} we set 
$$d:=\pi(c, c'), \, \, \, \, \, f_1:=x_1 \cdots x_{m-1}, \, \, \, \, \, f_2:=y_2 \cdots y_p.$$ 
By Proposition \ref{bigg. corolla}, either $f_1 f_2$ is a stem or an insertion of a stem into a flower element with g. corolla $d$ and $\eta(d)\leqslant \eta(c)+\eta(c')$.

We modify $\sigma$ by replacing $c$ and $c'$ with an SLA $\tau'$ defined in the following way: if $f_1 f_2$ is a stem then $\tau'$ has only one step which is equal to $f_1 f_2$; otherwise $\tau'$ is an SLA computing $f_1 f_2$ and having only one g. corolla, equal to $d$, therefore $A(\tau')=\eta(d)$. We consider the insertions at a letter of $c$ or $c'$ as insertions at the same letter of $f_1 f_2$. An insertion at $z^{-1}$ is replaced by an insertion at the last letter that in the output of the same step of $\sigma$ was preceding $z$. With this modification we obtain an SLA $\sigma'$ computing $u u' q_1=l_1 l_2$. 

Let $N$ and $N'$ be the multisets of base elements of $\sigma$ and $\sigma'$ respectively. For any g. corolla of $\sigma$ let us fix an SLA in $\overline{R}$ computing that g. corolla and let $M$ be the union of the multisets of the base elements of the SLA's in $\overline{R}$ for the elements of $N$. $M$ is a CMDR for $\sigma$. If $f_1 f_2$ is a stem then $N'=N\setminus{\{c, c'\}}$ and thus $A(\sigma') < A(\sigma)$. Then $M$ minus the base elements of the chosen SLA's in $\overline{R}$ for $c$ and $c'$ is a CMDR for $\sigma'$. Let $f_1 f_2$ be not a stem; then $N'\setminus{\{d\}}=N\setminus{\{c, c'\}}$, which means that the g. corollas of $\sigma$ and $\sigma'$ coincide except that $c$ and $c'$ are replaced in $\sigma'$ by $d$; therefore the inequality $\eta(d)\leqslant \eta(c)+\eta(c')$ implies that $A(\sigma')\leqslant A(\sigma)$. Moreover $M$ is a CMDR for $\sigma'$.

If $A(\sigma)=\textrm{Area}\big(\rho(l)\big)$ then $A(\sigma')=A(\sigma)$ by Remark \ref{ora}, which implies by Remark \ref{dopo} that a CMDR for $\sigma$ is a CMDR also for $\sigma'$.

\section{The proof of the Main Theorem: conclusion}  \label{conc}

In this section we prove Lemma \ref{maile} in the general case. Let $l:=l_1 z z^{-1} l_2$ be an element of $L_g$ computed by a straight line algorithm $\sigma$; we construct an SLA $\sigma'$ computing $l_1 l_2$ and such that $A(\sigma')\leqslant A(\sigma)$. Moreover if $A(\sigma)=\textrm{Area}\big(\rho(l)\big)$ then $A(\sigma')=A(\sigma)$ and a CMDR for $\sigma$ is a CMDR also for $\sigma'$.

\begin{lemma} \label{tol} Let $\sigma$ be an SLA, let $l$ be its result and let $A(\sigma)=\textrm{Area}\big(\rho(l)\big)$. If $\tau$ is a pSLsA of $\sigma$ with result $w$ then $A(\tau)=\textrm{Area}\big(\rho(w)\big)$. 
\end{lemma}

\begin{proof} We prove the claim by induction of the number of steps of $\sigma$. If $\sigma$ has only one step the claim is obvious because the only pSLsA is $\sigma$ itself. Let the number of steps of $\sigma$ be greater than one and the claim be true for every SLA with less steps than $\sigma$; let $\tau$ be a pSLsA of $\sigma$ and let $w$ be its result. The last step of $\sigma$ is the insertion of a word $l_2$ into a word $l_1$. By Part 2 of Proposition \ref{mudepends}, $l$ uses $w$ and since $l$ uses directly $l_1$ and $l_2$, then by Part 4 of the same proposition either $l_1$ or $l_2$ uses $w$. This implies that if $\sigma_1$ and $\sigma_2$ are the pSLsA's of $\sigma$ computing $l_1$ and $l_2$, then $\tau$ is a pSLsA of $\sigma_1$ or $\sigma_2$. The claim follows thus by induction hypothesis because $\sigma_1$ and $\sigma_2$ have less steps than $\sigma$.
       \end{proof}

Let $s$ be the first step of $\sigma$ whose output contains as subword the subword $z z^{-1}$ of $l_1 z z^{-1} l_2$. Let $w:=x_1 \cdots x_m$ be the output of $s$ and let $k: 1 \leqslant k < m$ be such that $x_k x_{k+1}=z z^{-1}$. By Remark \ref{lambda} we can reorder the steps of $\sigma$ in such a way that $s$ depends on every step preceding it. By Remark \ref{pSLsA} these steps form a pSLsA $\tau$ whose result is $w$. By the results of the preceding section we have that $w':=x_1 \cdots x_{k-1} x_{k+2} \cdots x_m$ belongs to $L_g$ and that there exists an SLA $\tau'$ computing $w'$ and such that $A(\tau')\leqslant A(\tau)$. Furthermore, by Lemma \ref{tol}, if $A(\sigma)=\textrm{Area}\big(\rho(l)\big)$ then $A(\tau)=\textrm{Area}\big(\rho(w)\big)$ and thus by the results of the preceding section, $A(\tau')=A(\tau)$ and a CMDR for $\tau$ is a CMDR also for $\tau'$.

If $s$ is the last step of $\sigma$, then we are in the situation of the preceding section. Suppose that $s$ is not the last step; for every step $t$ in $\sigma$ that follows $s$ we define a new step $t'$ in the following way. Let $t$ be the first step that follows $s$; $t$ cannot be an insertion because it can use directly only $s$ (the steps preceding $s$ are already used by $s$ by the observation made earlier), therefore $t$ is a base step. We set $t':=t$. Let $n>1$ and let $t$ be the $n$-th step that follows $s$. If $t$ does not use $s$ then we set $t':=t$. If $t$ does, $t$ is an insertion of a step $t_2$ into a step $t_1$ and by Part 4 of Proposition \ref{mudepends}, one and only one between $t_1$ and $t_2$ uses $s$. We can assume by induction hypothesis that we have already defined $t'_1$ and $t'_2$. We let $t'$ be the insertion of $t'_2$ into $t'_1$ at the same letter as $t_2$ is inserted into $t_1$. This letter cannot be equal to $x_k$, because either $t_1$ does not contain it or does contain both $x_k$ and $x_{k+1}$ consecutively. If $t$ is the insertion of $t_2$ into $t_1$ at $x_{k+1}$, then we let $t'$ be the insertion of $t'_2$ into $t'_1$ at the letter of $t_1$ that precedes $x_k$.

For every step $t$, the output of $t'$ is equal to the output of $t$ if $t$ does not use $s$; if $t$ uses $s$ then the output of $t'$ is the word obtained by cancelling $z z^{-1}$ from the output of $t$. If $t$ is the last step, then the output is $l_1 l_2$. Replacing $\tau$ with $\tau'$ and every step $t$ that follows $s$ with the corresponding step $t'$ we obtain an SLA $\sigma'$ computing $l_1 l_2$. 

The base steps of $\sigma$ following $s$ coincide with those of $\sigma'$ following $s'$. This implies that since $A(\tau')\leqslant A(\tau)$ then $A(\sigma')\leqslant A(\sigma)$. Moreover we have seen that if $A(\sigma)=\textrm{Area}\big(\rho(l)\big)$ then $A(\tau')=A(\tau)$ and a CMDR for $\tau$ is a CMDR also for $\tau'$; this implies that $A(\sigma')=A(\sigma)$ and a CMDR for $\sigma$ is a CMDR also for $\sigma'$.

 \chapter{Applications and complexity} \label{aac}

In Section \ref{worpro} we illustrate some applications to the Word Problem of the results of this thesis. In Section \ref{simpli} we show how, given an SLA in $L$ whose result is a reduced word, one can find an SLA with the same result and whose intermediate outputs are reduced words. Finally in Section \ref{comple} we show that the algorithm presented in this thesis to compute the relators has a better complexity than the other methods currently known in the literature.

\section{Applications to the Word Problem} \label{worpro}

Let $\langle \, X \, | \, R \, \rangle$ be a group presentation and let $R$ be finite. Consider the following functions introduced in Definition \ref{propcor}:
   $$\Delta_0(n):=\max\{\textrm{Area}(w) : \textrm{$w$ is a proper corolla} \,\,\textrm{and} \, \, |w|\leqslant n\},$$
   $$\Omega_0(n):=\max\{\textrm{Work}(w) : \textrm{$w$ is a proper corolla} \,\,\textrm{and} \, \, |w|\leqslant n\}.$$
  It is obvious that $\Delta_0\leqslant \Delta$ and  $\Omega_0\leqslant \Omega$ where $\Delta$ is the Dehn function (Definition \ref{Dehn}) and $\Omega$ has been introduced in Definition \ref{omega}.

  By Proposition \ref{gersho} a finite presentation has a solvable Word Problem if and only if its Dehn function is bounded above by a computable function if and only if its Dehn function is computable. Let us improve this result.
    
Let $n$ be a natural number and let $w$ be a relator of length at most $n$. By Corollary \ref{maicor} there exist proper corollas $c_1, \cdots, c_m$ such that $\textrm{Area}(w)= \sum_{i=1}^m \textrm{Area}(c_i)$ and $|w|\geqslant |c_i|$ for every $i$. Thus 
      $$\textrm{Area}(w)\leqslant \Delta_0(|c_1|)+\cdots +\Delta_0(|c_m|),$$
 that is     
      $$\Delta(n) \leqslant \max\{\Delta_0(k_1)+\cdots +\Delta_0(k_m) : k_1+\cdots +k_m\leqslant n\}\leqslant n \Delta_0(n)$$
 and then we have the following inequality:
  \begin{equation} \label{jona} \Delta_0(n) \leqslant \Delta(n) \leqslant n \Delta_0(n).\end{equation}

This means that if $\Delta_0$ is bounded by a sequence $\{k_n\}_{n\in \mathbb{N}^*}$ then $\Delta$ is bounded by $\{n k_n\}_{n\in \mathbb{N}^*}$. As we have said, this bound can be improved if $k_n$ is linear in $n$ in which case $\Delta$ is bounded by $k_n$ (see Corollary \ref{hyper3}).

If $\{k_n\}_{n\in \mathbb{N}^*}$ is computable then also $\{n k_n\}_{n\in \mathbb{N}^*}$ is computable. Thus if $\Delta_0$ is bounded above by a computable function then so is $\Delta$ and the Word Problem is solvable. 

We prove that in this case $\Delta_0$ is computable. Since the presentation is finite we can compute the results of all the $L-SLA$'s $\tau$ such that $A(\tau)\leqslant k_n$ and which compute a word of length at most $n$. In this way we determine all the proper corollas of length at most $n$ and their area, thus we compute $\Delta_0$. Thus

 \begin{theorem} \label{aska} A finite presentation has a solvable Word Problem if and only if the function $\Delta_0$ is bounded by a computable function if and only if the function $\Delta_0$ is computable. Moreover $\Delta$ is computable or bounded by a computable function if and only if the same is true for $\Delta_0$.
\end{theorem}

In the same way one sees that
    $$\Omega_0(n) \leqslant \Omega(n) \leqslant n \Omega_0(n)$$
and thus we improve part 1 of Proposition \ref{risGI} with the following

\begin{theorem} A finitely generated decidable presentation has a solvable Word Problem if and only if the function $\Omega_0$ is bounded by a computable function if and only if the function $\Omega_0$ is computable. \end{theorem}

We now exhibit an explicit algorithm for solving the Word Problem for a finite presentation for which $\Delta_0$ is bounded by a computable function. This gives a direct proof of Theorem \ref{aska}, i.e., a proof not relying on Proposition \ref{gersho}.

The algorithm is the following. Let $\langle \, X \, | \, R \, \rangle$ be a finite presentation and for every $k$ let $C_k$ (Definition \ref{kg. corolla}) be the set of $k$-corollas. As we have seen in Remark \ref{Roverlinen}, there is an explicit algorithm for computing $C_k$. Suppose that $\Delta_0$ is bounded by a computable function, that is $\Delta_0(n) \leqslant k_n$ for every $n$ where $k_n$ is computable by a finite algorithm and let $w$ be a reduced word in the alphabet $X\cup X^{-1}$.

 By Corollary \ref{cont}, if $w$ is a relator of $\langle \, X \, | \, R \, \rangle$ then $w$ has a contiguous subword equal to a proper corolla whose area is at most equal to the area of $w$. Let $n$ be the length of $w$; if no element of $C_h$ for $h\leqslant k_n$ is a contiguous subword of $w$, then $w$ is not a relator. Suppose on the contrary that there exists an element $c$ of $C_h$ for some $h\leqslant k_n$ which is a contiguous subword of $w$; this means that there exist reduced words $w'$ and $w''$ such that $w=w' c w''$. Let $w_1$ be the reduced form of $w' w''$; if $w_1=1$ then $w$ is a relator. Let $w_1\neq 1$; $w$ is a relator if and only $w_1$ is a relator. Moreover the length of $w_1$ is less than that of $w$.

If $w_1$ is a relator, its area is equal at most to the minimum between $k_n-h$ and $k_{|w_1|}$. Let $m_1$ be that minimum; if no element of $C_h$ for $h \leqslant m_1$ is a contiguous subword of $w_1$ then $w_1$ is not a relator and neither $w$; if yes, we repeat for $w_1$ the same procedure done for $w$ and we find a word $w_2$ of length less than thath of $w_1$. This algorithm stops after at most $n$ steps: if for some $i=1, 2, ...$, we have that $w_i=1$, then $w$ is a relator; if no contiguous subword of $w_i$ is equal to an element of $C_h$ for $h\leqslant m_i$ then $w$ is not a relator. This concludes the algorithm.

\medskip

For every $h$ set $P_h:=C_h\setminus{\{C_1\cup \cdots \cup C_{h-1}\}}$ and let $m_h$ be the minimal length of words of $P_h$. If $X$ is finite and if $m_h\neq 0$ for infinitely many $h$ then the sequence $\{m_h\}_{h\in \mathbb{N}^*}$ diverges to infinity. Indeed suppose that $\{m_h\}_{h\in \mathbb{N}^*}$ does not diverge to infinity. Then there would exist a natural number $n$ and an infinite number of indices $h$ such that $m_h\leqslant n$. Since the $P_h$ are disjoint two by two, there would be infinitely many words of length not greater than $n$ and this is impossible because their number is bounded by $(2|X|)^n$.

\begin{lemma} \label{nebr} For any $n$ let $f(n)$ be the least natural number such that $m_h>n$ for every $h>f(n)$. Then $\Delta_0(n)\leqslant f(n)$.
\end{lemma}

\begin{proof} Let $h>f(n)$; since $m_h>n$ then if $w \in P_h$ we have that $|w|>n$. Let $w$ be a proper corolla of length not greater than $n$; then $w$ does not belong to $P_h$ for every $h>f(n)$, that is it belongs to $P_1\cup \cdots \cup P_{f(n)}$. By Definition \ref{propcor} and the Main Theorem \ref{maithe}, $\textrm{Area}(w)=\eta(w)$ and thus the area of $w$ cannot be greater than $f(n)$, which implies that $\Delta_0(n)\leqslant f(n)$.
 \end{proof}

 Lemma \ref{nebr} implies:

 \begin{theorem} \label{linc} The Word Problem is solvable for the finite presentation $\langle \, X \, | \, R \, \rangle$ if and only if the function $f$ defined in Lemma \ref{nebr} is bounded above by a computable function.
\end{theorem}

\begin{proof} If $f$ is bounded above by a computable function then by Lemma \ref{nebr} also $\Delta_0(n)$ is bounded above by the same computable function and the Word Problem for $\langle \, X \, | \, R \, \rangle$ is solvable by Theorem \ref{aska}.

Let the Word Problem for $\langle \, X \, | \, R \, \rangle$ be solvable; by Theorem \ref{aska} there exist computable $k_n$ for every $n$ such that $\Delta_0(n) \leqslant k_n$. That is, if $w$ is a proper corolla and $|w|\leqslant n$ then $\textrm{Area}(w) \leqslant k_n$. This means that if $w$ is a proper corolla and $\textrm{Area}(w) > k_n$ then $|w|>n$. Then $m_h>n$ for every $h> k_n$, thus $k_n\geqslant f(n)$ and $f$ is bounded above by a computable function. \end{proof}

Theorem \ref{linc} says that the Word Problem is solvable for $\langle \, X \, | \, R \, \rangle$ if and only if for every $n$ we can compute a $k_n$ such that every proper corolla of area greater than $k_n$ has length greater than $n$.

 \begin{corollary} \label{} Let $\{a_h\}_{h\in \mathbb{N}^*}$ be a computable non-decreasing and divergent sequence. If $m_h\geqslant a_h$ for every $h$ then the Word Problem for $\langle \, X \, | \, R \, \rangle$ is solvable. \end{corollary}

\begin{proof} Let $n$ be a natural number; then there exists a natural number $k_n$ such that $a_h>n$ for every $h>k_n$. We have that $f(n)\leqslant k_n$ where $f$ is the function defined in Lemma \ref{nebr}. Let us compute $a_1, a_2, \cdots$ until we find an $a_i$ greater than $n$. Then $k_n$ is equal to this $i$ and there exists a finite algorithm computing $k_n$. The Word Problem for $\langle \, X \, | \, R \, \rangle$ is solvable by Theorem \ref{linc}.\end{proof}

  Corollary \ref{hyper2} says that if $\Delta_0(n)$ is bounded by a linear function in $n$ then also $\Delta$ is bounded by the same linear function and the presentation is hyperbolic. We now generalize this fact.

We recall (Definition \ref{equiv}) that we write $f \preceq g$ if there exists a positive constant $\alpha$ such that $f(n)\leqslant \alpha g(\alpha n) + \alpha n$ for every $n\in \mathbb{N}^*$. Obviously we have that $\Delta \preceq g$ [respectively $\Delta_0 \preceq g$] if and only if there exists a positive constant $\alpha$ such that 
  \begin{equation} \label{azer} \textrm{Area}(w) \leqslant \alpha g(\alpha |w|)+ \alpha |w| \end{equation}
 for every relator $w$ [respectively for every proper corolla $w$].

 \begin{lemma} \label{pott} Let $g: \mathbb{R}_+ \rightarrow \mathbb{R}_+$ be a non-decreasing function and suppose that there exists $\alpha_0>0$ such that for every $\alpha \geqslant \alpha_0$ there exists $\beta \in \mathbb{R}_+$ such that for every $y_1, \cdots, y_m, z \in  \mathbb{N}^*$ such that $y_1+ \cdots + y_m\leqslant z$ we have
                \begin{equation} \label{alfbet} \sum_{j=1}^m g(\alpha y_j) \leqslant g(\beta z). \end{equation}
If $\Delta_0 \preceq g$ then $\Delta \preceq g$. \end{lemma}

\begin{proof} If $\Delta_0 \preceq g$ then there exists $\alpha \in \mathbb{R}_+$ such that for every proper corolla $c$ we have that $\textrm{Area}(c) \leqslant \alpha g(\alpha |c|) + \alpha |c|$. By Remark \ref{lis} we can suppose that $\alpha \geqslant \alpha_0$.

Let $n$ be a natural number and let $w$ be a relator of length at most $n$. By Corollary \ref{maicor} there exist proper corollas $c_1, \cdots, c_m$ such that $\textrm{Area}(w)= \sum_{j=1}^m \textrm{Area}(c_j)$ and $\sum_{j=1}^m |c_j|\leqslant |w|$. Then 
$$\textrm{Area}(w) \leqslant \alpha \sum_{j=1}^m \big[g(\alpha |c_j|) + |c_j|\big]$$
 and by (\ref{alfbet}) there exists $\beta$ such that 
 $$\textrm{Area}(w) \leqslant \alpha [g(\beta |w|)+ |w|].$$
    Since $g$ is non-decreasing, if we set $\beta':=\max\{\alpha, \beta\}$, then 
 $$\textrm{Area}(w) \leqslant \alpha [g(\beta |w|)+ |w|]\leqslant \beta' [g(\beta' |w|)+ |w|]$$
  and the claim follows from (\ref{azer}).
\end{proof}

 We recall (Definition \ref{equiv}) that we write $f \simeq g$ if $f \preceq g$ and $g \preceq f$.

 \begin{theorem} Let $g: \mathbb{R}_+ \rightarrow \mathbb{R}_+$ be a non-decreasing function satisfying (\ref{alfbet}) and let $\Delta_0 \simeq g$; then $\Delta \simeq g$ and thus $\Delta_0 \simeq \Delta$.
\end{theorem}

\begin{proof} We have that $\Delta_0 \preceq g$ and $g \preceq \Delta_0$. By Lemma \ref{pott} we have that $\Delta \preceq g$. Conversely we have that $g \preceq \Delta$ because $g \preceq \Delta_0$ and $\Delta_0 \leqslant \Delta$. \end{proof}

We now show that polynomial and exponential functions satisfy (\ref{alfbet}).

\begin{lemma} \label{fout} \begin{enumerate}

  \item Let $r>1$ and let $h(x)=(x+1)^r-x^r-1$. Then $h(x)>0$ for every $x>0$.
  
  \item Let $r>1$ and let $y_1, \cdots, y_m> 0$. Then $(y_1+ \cdots + y_m)^r>y_1^r + \cdots + y_m^r$.

\end{enumerate}

\end{lemma}

\begin{proof} \begin{enumerate}

  \item We have that $h(0)=0$ and that the derivative of $h$ is   
      $$h'(x)=r(x+1)^{r-1}-rx^{r-1}=r[(x+1)^{r-1}-x^{r-1}]$$
  which is positive for every $x\geqslant 0$.

  \item  It is enough to prove the claim for $m=2$. Since $y_1, y_2 > 0$ then there exists $c>0$ such that $y_2=cy_1$, therefore the inequality is equivalent to $y^r + c^r y^r > (c+1)^r y^r$ for every $y, c > 0$ and the latter is equivalent to $(1+c)^r>1+c^r$ for every $c > 0$. By Part 1, $(1+c)^r-c^r-1>0$ for every $c>0$ and thus we have the claim.
  \end{enumerate}
\end{proof}

\begin{remark} \label{brock} \rm For $r=1$ we have that $(y_1+ \cdots + y_m)^r=y_1^r + \cdots + y_m^r$; then by Part 2 of Lemma \ref{fout} we have that if $r\geqslant 1$ then $(y_1+ \cdots + y_m)^r\geqslant y_1^r + \cdots + y_m^r$. \end{remark}

\begin{proposition} \label{vega} If $r$ is a real number and $r\geqslant 1$ then the polynomial functions $x^r$ satisfy (\ref{alfbet}).
\end{proposition}

\begin{proof} We will prove that for any $y_1, \cdots, y_m, z \in \mathbb{N}^*$ such that $y_1+ \ldots + y_m \leqslant z$ we have
    $$ \sum_{j=1}^m (\alpha y_j)^r \leqslant (\alpha z)^r.$$

We have that 
$$\sum_{j=1}^m (\alpha y_j)^r \leqslant \Big(\sum_{j=1}^m \alpha y_j\Big)^r \leqslant (\alpha z)^r$$  
 when the first inequality follows from Remark \ref{brock}.
  \end{proof}

\begin{remark} \label{margo} \rm Let $r_1, \cdots, r_n\geqslant 1$ be real numbers; we show that there exists $1<a_0\leqslant 2$ such that $x^{r_1+ \cdots+ r_n}\geqslant x^{r_1}+ \cdots+ x^{r_n}$ for every $x\geqslant a$. It is sufficient to show the claim for $n=2$. Let $\varphi(x)=x^{r_1+ r_2}-x^{r_1}-x^{r_2}$ and suppose that $r_1\geqslant r_2$. Let us show that $\varphi(2)\geqslant 0$; we have that 
    $$\varphi(2)=2^{r_1+ r_2}-2^{r_1}-2^{r_2}=2^{r_1}(2^{r_2}-1-2^{r_2-r_1})$$
and $2^{r_2}-1-2^{r_2-r_1} \geqslant 0$ because $2^{r_2} \geqslant 2$ (being $r_2\geqslant 1$) and $2^{r_2-r_1} \leqslant 1$ (being $r_2-r_1 \leqslant 0$). Let us show that $\varphi'(x)\geqslant 0$; we have that $$\varphi'(x)=(r_1+r_2)x^{r_1+ r_2-1}-r_1x^{r_1-1}-r_2x^{r_2-1}$$
$$=r_1(x^{r_1+ r_2-1}-x^{r_1-1})+r_2(x^{r_1+ r_2-1}-x^{r_2-1})\geqslant 0.$$
 \end{remark}

 \begin{proposition} \label{ange} Let $a_0$ as in Remark \ref{margo}. Then for every $a\geqslant a_0$ the functions $a^x$ satisfy (\ref{alfbet}). \end{proposition}

\begin{proof} Let $y_1, \cdots, y_m, z\in \mathbb{N}^*$ be such that $y_1+ \cdots+ y_m\leqslant z$ and let $\alpha\geqslant 1$. Since $\alpha y_j\geqslant 1$ for $j=1, \cdots, m$ then by Remark \ref{margo} we have that 
$$\sum_{j=1}^m a^{\alpha y_j} \leqslant a^{\alpha z}.$$
    \end{proof}

We now show that for every $b>1$ the function $\log_bx$ does not satisfy (\ref{alfbet}). Let $\alpha, \beta>0$ and let $x, y$ be non-zero natural numbers. Then there exists a rational number $c$ such that $y=cx$. We have that 
   $$\log_b(\alpha x)+\log_b(\alpha y) \leqslant \log_b\big(\beta (x+y)\big)$$
if and only if
    $$\log_b(\alpha x)+\log_b(\alpha cx) \leqslant \log_b\big(\beta (c+1)x\big)$$
which is equivalent to  
 $$\log_b(\alpha^2c)+2\log_b(x) \leqslant \log_b\big(\beta (c+1)\big)+\log_b(x)$$
 and thus to 
  $$\log_b(x) + \log_b(\gamma)\leqslant 0$$
where $\gamma=\frac{\alpha^2c}{\beta (c+1)}$.
 But this is false because $\log_b(\gamma)$ is constant and $\log_b(x)$ goes to infinity.

\section{Straight line algorithms for reduced words} \label{simpli}

Let $\sigma$ be an SLA whose result is a reduced word; in this section we will show how to find an SLA equivalent in the sense of Definition \ref{equiSLA} to $\sigma$ and such that the output of every step is a reduced word.  
First we show that any SLA is equivalent to one in which the mid-vertex of a stem does not coincide with the initial vertex of another stem.

Consider the complex 

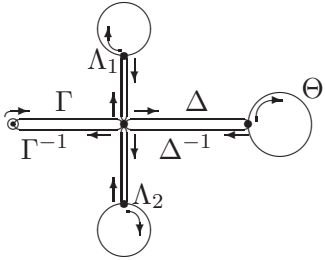
\begin{picture}(180,150)(40,-35)


\put(145,44){\oval(10,2)[tl]}   

\put(145,45){\vector(1,0){4}}   

 
\put(181,43){\vector(0,1){9}}

 
\put(189,65){\vector(0,-1){9}}


\put(184,73){\oval(10,10)[bl]}   
 
\put(179,73){\vector(0,1){4}}   


\put(189,45){\vector(1,0){9}}

 
\put(241,41){\oval(12,16)[tl]}   

\put(241,49){\vector(1,0){4}}   


\put(232,36){\vector(-1,0){9}}


\put(180,36){\vector(-1,0){9}}

 
\put(189,36){\vector(0,-1){9}}


\put(186,2){\oval(10,10)[tr]}   
 
\put(191,2){\vector(0,-1){4}}   

 
\put(181,11){\vector(0,1){9}}


\put(164,40){\oval(41,4)}  

\put(209,40){\oval(47,4)} 

\put(185,54){\oval(2,26)} 

\put(185,76){\circle{20}}  

\put(244,40){\circle{25}}  

\put(185,24){\oval(2,30)} 

\put(185,0){\circle{20}}  


\put(143,40){\circle*{2}}   

\put(143,40){\circle{4}}    

\put(185,40){\circle*{3}}       

\put(185,66){\circle*{3}}       

\put(232,40){\circle*{3}}      

\put(185,10){\circle*{3}}       


\put(171,61){$\Lambda_1$}

\put(159,45){$\Gamma$}

\put(208,45){$\Delta$}

\put(252,50){$\Theta$}

\put(198,27){$\Delta^{-1}$}

\put(146,27){$\Gamma^{-1}$}

\put(188,10){$\Lambda_2$}

\end{picture}

Its label is $\Gamma \, \Lambda_1 \,  \Delta \, \Theta \, \Delta^{-1} \, \Lambda_2 \, \Gamma^{-1}$, where: $\Gamma$ is the label of the first half of the stem on the left and $\Delta$ the label of that on the right; $\Theta$ is the label of the corolla on the right; $\Lambda_1$ and $\Lambda_2$ the labels of the upper and lower flowers respectively. The mid-vertex of the stem on the left coincides with the initial vertex of the stem on the right.

The next complex has the same label as the preceding one and the latter situation does not happen.

 
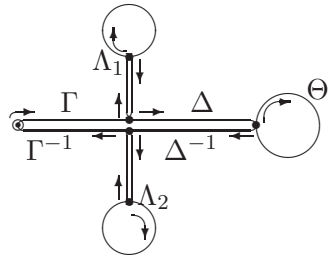
\begin{picture}(180,145)(32,-30)


\put(145,44){\oval(10,2)[tl]}   

\put(145,45){\vector(1,0){4}}   

 
\put(181,43){\vector(0,1){9}}

 
\put(189,65){\vector(0,-1){9}}


\put(184,73){\oval(10,10)[bl]}   
 
\put(179,73){\vector(0,1){4}}   


\put(189,45){\vector(1,0){9}}

 
\put(242,41){\oval(12,16)[tl]}   

\put(242,49){\vector(1,0){4}}   


\put(232,36){\vector(-1,0){9}}


\put(180,36){\vector(-1,0){9}}

 
\put(189,36){\vector(0,-1){9}}


\put(186,1){\oval(10,10)[tr]}   
 
\put(191,1){\vector(0,-1){4}}   

 
\put(181,12){\vector(0,1){9}}


\put(188,40){\oval(90,4)}  

\put(185,55){\oval(2,24)} 

\put(185,76){\circle{20}}  

\put(245,40){\circle{25}}  

\put(185,24){\oval(2,28)} 

\put(185,1){\circle{20}}  


\put(143,40){\circle*{2}}   

\put(143,40){\circle{4}}    

\put(185,42){\circle*{3}}       

\put(185,38){\circle*{3}}       

\put(185,66){\circle*{3}}       

\put(233,40){\circle*{3}}      

\put(185,11){\circle*{3}}       


\put(171,61){$\Lambda_1$}

\put(159,45){$\Gamma$}

\put(208,45){$\Delta$}

\put(252,50){$\Theta$}

\put(198,27){$\Delta^{-1}$}

\put(146,27){$\Gamma^{-1}$}

\put(188,10){$\Lambda_2$}

\end{picture}

The two stems labeled by $\Gamma  \Gamma^{-1}$ and $\Delta  \Delta^{-1}$ have been replaced by one labeled by $\Gamma \Delta  \Delta^{-1} \Gamma^{-1}$.

Let us formalize this situation. Let $\sigma$ be an SLA, let $w:=x_1\cdots x_m$ be its result and let $b_1:=y_1 \cdots y_p$ and $b_2$ be two stems such that the mid-vertex of $b_1$ coincides with the initial vertex of $b_2$. We formalize the notion of “sharing a vertex" and then show how to modify $\sigma$ in order to avoid it.

Let $b_1$ comprise $b_2$. By Proposition \ref{consec} we can suppose that $\sigma$ defines consecutively the insertions into $b_1$. Since $b_2$ is comprised in $b_1$ then all the letters of $b_2$ are comprised between two consecutive letters $y_h=x_{i_1}$ and $y_{h+1}=x_{j_1}$ of $b_1$, where $1 \leqslant h < p$ and $1\leqslant i_1<j_1\leqslant m$. Therefore by Definition \ref{def.cons}, $b_2$ belongs to the pSLsA of $\sigma$ which computes the ramification (Definition \ref{ramsur}) from $b_1$ at $y_h$. Call $\tau_h$ and $r_h$ this pSLsA and this ramification respectively. By Proposition \ref{consec} we can suppose that $\tau_h$ defines consecutively the insertions into $b_2$.

Let $x_{i_2}$ and $x_{j_2}$ (for some indices $i_2$ and $j_2$) be respectively the first and the last letters of $b_2$; then we have that $1\leqslant i_1<i_2<j_2<j_1\leqslant m$. Set $v_1:=x_{i_1+1} \cdots x_{i_2-1}$ and $v_2:=x_{j_2+1} \cdots x_{j_1-1}$. We call $v_1$ and $v_2$ \textit{the subwords comprised between $b_1$ and $b_2$}.

Suppose that $y_h$ is the last letter of the first half of $b_1$; we say that \textit{the mid-vertex of $b_1$ coincides with the initial vertex of $b_2$} if the subwords comprised between $b_1$ and $b_2$ are parts of $w$ (Definition \ref{SLsA}).

 $v_1 v_2$ is the surround of $b_2$ in the ramification $r_h$ and by Definition \ref{def.cons} we can suppose that $v_1 v_2$ is the result of a pSLsA of $\tau_h$; call $\tau'_0$ this pSLsA. In particular $v_1$ and $v_2$ are respectively the preceding and the following of $b_2$ in $r_h$ and if we have supposed that they are parts of $w$, we can assume that $\tau'_0$ is formed by two pSLsA's computing $v_1$ and $v_2$ respectively, followed by the product of $v_1$ by $v_2$. We call $t$ the step of the product of $v_1$ by $v_2$.

We modify $\sigma$ by replacing $b_1$ and $b_2$ with $s:=y_1 \cdots y_h \, b_2 \, y_{h+1} \cdots y_p$ which is a stem because $y_h$ is the last letter of the first half of $b_1$. Furthermore we replace the step $t$ by the insertions of $v_1$ at $y_h$ and of $v_2$ at the last letter of $b_2$ and we consider the insertions at a letter of $b_1$ or $b_2$ as insertions at the corresponding letters of $f$ or $s$ respectively.

We obtain an SLA $\sigma'$ with result $w$ (the same of $\sigma$). $\sigma$ and $\sigma'$ have the same corollas, therefore $A(\sigma')=A(\sigma)$ and a CMDR for $\sigma$ is a CMDR also for $\sigma'$.

With the notion of 'equivalence' of SLA's of Definition \ref{equiSLA}, $\sigma$ and $\sigma'$ cannot be equivalent but it is intuitive that $\sigma$ and $\sigma'$ are equivalent in a more general sense. Let us modify this notion of equivalence.

Let $\tau_1$ and $\tau_2$ be SLA's with the same result (let $w$ be their result), with multisets of base elements $M_1$ and $M_2$ and let $\omega: M_1 \rightarrow M_2$ be an application such that: \begin{enumerate}

  \item the restriction of $\omega$ to the corollas of $M_1$ is an isomorphism onto the corollas of $M_2$;
  
  \item the restriction of $\omega$ to the stems of $M_1$ is surjective onto the stems of $M_2$;
  
  \item let $vv^{-1}$ be a stem of $\tau_2$ and let $u_1u_1^{-1}, \cdots, u_n u_n^{-1}$ be the stems of $\tau_2$ whose image by $\omega$ is $vv^{-1}$ and let the first letter of $u_i u_i^{-1}$ precedes the first of $u_{i+1} u_{i+1} ^{-1}$ for $i=1, \cdots, n-1$. Then $vv^{-1}=u_1 \cdots u_n u_n^{-1} \cdots u_1^{-1}$;
  
  \item any letter of $w$ comes (Definition \ref{coming from}) from the same letter of $\mu$ and $\omega(\mu)$ respectively (we recall that $\mu$ is a subword of $\omega(\mu)$).
\end{enumerate}

We say that $\tau_2$ is a \textit{refinement} of $\tau_1$.

Let $\sigma$ be an SLA, let $w$ be its result and let $A(\sigma)=A(w)$. If we apply repeatedly this procedure, we obtain an SLA $\sigma''$ with the same result of $\sigma$, such that $A(\sigma'')=A(w)$, such that a CMDR for $\sigma$ is a CMDR also for $\sigma''$ and in which the mid-vertex of any stem does not coincide with the initial vertex of another stem. Moreover $\sigma''$ is a refinement of $\sigma$.

Any element $w$ of $L_g$ (and thus any relator of $\langle X  | R \rangle$ by Theorem \ref{maithe}) is the result of an SLA $\sigma$, that is it is obtained by means of the operation of grafting starting with base elements (stems and corollas). By what said above, we can suppose that the mid-vertex of any stem of $\sigma$ does not coincide with the initial vertex of another stem.

Suppose that $w$ is reduced and suppose that one of the steps of $\sigma$ is a stem $uu^{-1}$. We can suppose that $\sigma$ defines consecutively the insertion in $uu^{-1}$. Since $uu^{-1}$ is non-reduced and since the result of $\sigma$ is a reduced word then the ramification $r$ from the last letter of $u$ is non-empty. Since the mid-vertex of $uu^{-1}$ does not coincide with the initial vertex of another stem, then the first and the last letter of this ramification come from corollas. Let $x$ be the last letter of $u$, that is $x^{-1}$ is the first letter of $u^{-1}$, and let $c$ be the corolla from which comes the first letter of $r$. The first letter of $c$ is different from $x^{-1}$. If the last letter of $c$ is different from $x$ then the flower $ucu^{-1}$ is reduced; otherwise, the subword comprised between $c$ and $u^{-1}$ is non-empty, the first letter is different from $x^{-1}$ and comes from a corolla. By repeating the same procedure, we find that there exists a corolla $c'$ such that $uc'u^{-1}$ is reduced and the subwords comprised between $u$ and $c'$ and between $c'$ and $u^{-1}$ are parts of $w$. That is in the multiset of base elements of $\sigma$ we can replace $uu^{-1}$ and $c'$ with $uc'u^{-1}$ and then we can assume that the base elements are not stems and corollas but reduced flowers and corollas. Since a corolla can be considered a flower with trivial stem, then if $w$ is reduced then we assume that the base elements are the reduced flowers.

We now prove that if $w$ is reduced then then $\sigma$ is equivalent to an SLA the outputs of whose steps are all reduced. Suppose that $\sigma$ has only one step; then this step is reduced since it is a reduced flower. Let the claim be true for every SLA with less steps than $\sigma$. Let $f$ be the flower from which comes the first letter of $w$. The ramifications from $f$ are reduced words since they are contiguous subwords of a reduced word; the surround of $f$ coincides with its following, therefore it is reduced since it is a contiguous subword of $w$. By induction hypothesis they are defined by SLA's whose outputs are all reduced and finally inserting the ramifications and the surround into $f$ give also reduced outputs.

We denote $D_1$ the set of words which are reduced and which are of the form $u w u^{-1}$ where $w \in C_1$. Suppose by induction to have defined $D_{k-1}$ and define $D_k$ as the set of reduced words which are either of the form $u w u^{-1}$ where $w \in C_k$ or are insertions of a word of $D_m$ into one of $D_n$ for $m+n=k$. We denote $D$ the union of all the $D_k$. Then by what said above, given an $L_g$-straight line algorithm $\sigma$ there exists a straight line algorithm $\sigma'$ whose base elements are reduced flowers, whose operations are insertions of words, such that the output of any step of $\sigma'$ is reduced and such that $\sigma'$ is a refinement of $\sigma$.

 \section{Comparison of the complexities} \label{comple}

Let $\mathcal{P}:=\langle \, X \, |  \, R  \, \rangle$ be a finite group presentation, let $k$ and $n$ be natural numbers and suppose that we want to find all the relators of $\mathcal{P}$ of length at most $n$ and which are products of at most $k$ conjugates of defining relators. In Section \ref{enume} we have illustrated the two methods currently known in the literature; both use the van Kampen diagrams associated with relators. We have also computed (\ref{schi}) and (\ref{frot}), the complexities of these methods and shown that the better is (\ref{frot}), which is a product of three exponentials, one in $k$, one in $kn$ and one in $k^2$. Namely this complexity is equal to $q^k a^{kn} s^{k^2}$ where $q=p (a+1)/(a-1)$, $p=|R|$, $a=2|X|-1$, $s=a^m$ and $m$ is the maximal length of an element of $R$. We have the following result

\begin{theorem} \label{furish} The complexity of the method illustrated in this thesis is bounded above by $n^3 \, \beta^n \, \gamma^k$ if $a>4$, where $\beta=2\sqrt{a}$ and $\gamma=4p\sqrt{2}\,^m$; it is bounded above by $n^3 \, 4^n \, \gamma^k$ if $a <4$. This complexity is less than those of the two classical methods (illustrated in Section \ref{enume}) by a factor of at least $a^{k(k-c_1)} a^{n(k-c_2)}$, with $c_1=m/2+2$ and $c_2=2$. \end{theorem}

Obviously $a$ cannot be equal to 4 since it is odd. The proof of Theorem \ref{furish} is given in Section \ref{fuzumi}.

Theorem \ref{furish} proves that an upper bound for the complexity of our method is a product of a cubic polynomial times an exponential in $n$ and one in $k$. This is much better than the complexities with the methods currently known in the literature. Moreover the value given in Theorem \ref{furish} is far from being an optimal bound for the complexity of our method, but only a value obtained with not too complicated calculations. With non-trivial arguments this result could be probably much improved. Instead (\ref{frot}) is the real complexity of the classical method, that is it is the exact number of words to be considered with this method to find all the relators of length up to $n$ and area up to $k$.

  \smallskip

Finally let us say some words about the comparison of  the method of this thesis with the first classical method illustrated in Section \ref{enume}, i.e., the method consisting in enumerating the van Kampen diagrams with up to $k$ faces and boundary of up to $n$ edges. As shown in Section \ref{ggior}, with every SLA we can associate a van Kampen diagram. 

In the diagrams associated with an SLA, the biconnected components correspond to the corollas. Thus the diagrams constructed with our method verify the following property: if $\mathcal{B}$ is a biconnected component of the diagram then one can name the faces of $\mathcal{B}$ as $F_1, \cdots, F_n$ in such a way that two $F_i$ are distinct and $F_1 \cup F_2 \cup \cdots \cup F_i$ is biconnected and has a cyclically reduced boundary label for every $i=1, \cdots, n$.

In the classical method one has to construct all the diagrams and there are diagrams not verifying the property above, in particular those whose boundary label is not cyclically reduced.

The central result of this thesis is that for finding all the relators one does not need to construct all the van Kampen diagrams, but only those associated with SLA's. The set of diagrams associated with SLA's is a proper subset of the set of all the diagrams, anyway it contains all the “important information” about the set of relators. 

Let us explain the last sentence. Let $a_1, \cdots, a_m$ be words and $r_1, \cdots, r_m$ be defining relators, let $w:=a_1 r_1 a_1^{-1} \cdots a_m r_m a_m^{-1}$ and let $v$ be the reduced form of $w$. There is a van Kampen diagram with boundary label equal to $v$ for every sequence of cancellations from $w$ to $v$; and two of these diagrams are not necessarily equal. In Chapter \ref{tmt} we have proved that there is a diagram for an SLA for at least one of these sequences of cancellations, thus by considering the set of SLA diagrams instead of that of all van Kampen diagrams no information about the relators is lost.

\appendix
         
\chapter{} \label{app}

In Section \ref{fuzumi} we compute the complexity of the algorithm presented in this thesis for computing the relators. In Section \ref{rdsaslp} we prove a result which is not needed in this thesis; in particular we prove that a recursively defined set can be defined as the set of results of a family of straight line algorithms. Finally in Section \ref{alc} we show a result necessary for the proof of the equality between the set of corollas with that of the generalized corollas.

 \section{Computation of the complexity} \label{fuzumi}

 Suppose given a finite presentation, let $n$ and $k$ be natural numbers and let us compute all the relators of length up to $n$ and area up to $k$. In this section we show that an upper bound for the complexity of the method presented in this thesis is a product of a cubic polynomial times an exponential in $n$ and one in $k$. This is much better than (\ref{frot}), the complexity with the method currently known in the literature, which is a product of three exponentials, one in $k$, one in $kn$ and one in $k^2$.

 \smallskip
 
 Let $\mathcal{P}:= \langle X \, | \, R \rangle$ be a finite group presentation. Let $m$ be the maximal length of elements of $R$, let $p:=|R|$ and let $a:=2|X|-1$.

 Let $n$ and $k$ be natural numbers and let $w$ be a relator of length at most $n$ and which is product of at most $k$ conjugates of defining relators. By Theorem \ref{maithe} and Section \ref{simpli}, $w$ is obtained by insertions of at most $n$ flowers. That is there exist $h\leqslant k$ flowers $f_1, \cdots, f_h$ such that $|f_1|+\cdots+|f_h|\leqslant n$ and $w$ is the result of an SLA in $D$ whose base elements are $f_1, \cdots, f_h$. This means that $f_2$ is inserted into $f_1$ giving a word $f'_2$, then $f_3$ is inserted into $f'_2$ giving a word $f'_3$ and so on until we have $f'_h=w$.
 
 Let $u:=u_1 u_2$ and $v$ be words; the insertion $u_1 v u_2$ of $v$ into $u$ is the reduced form of $u u_2^{-1} v u_2$, therefore the number of possible insertions of $v$ into $u$ is equal to $|u|$. Let $n_i=|f_i|$; since $|f'_i|=|f_1|+\cdots+|f_i|$, then the number of possible insertions of words $f_1, \cdots, f_h$ is 
      $$n_1+(n_1+n_2)+(n_1+n_2+n_3)+\cdots+(n_1+n_2+n_{h-1})=\sum_{i=1}^{h-1} (h-i)n_i.$$
  A flower of length $n_i$ is obtained by the insertion of a corolla $c_i$ of length $n'_i$ into a stem of length $n_i-n'_i$ for some $n'_i=1, \cdots, n_i$. Moreover $n_i-n'_i$ must be even. The number of stems of length $n_i-n'_i$ is equal to the number of reduced words of length $(n_i-n'_i)/2$, which is 1 if $n_i-n'_i=0$, is $a+1$ if $n_i-n'_i=2$ (we recall that $a=2|X|-1$) and is equal to $(a+1)a^{(n_i-n'_i)/2-1}$ otherwise.

Let $k_i=\eta(c_i)$; then $k_1+\cdots + k_h\leqslant k$. A corolla is determined by a sequence of cyclically reduced products of 1-corollas, each cyclically reduced product followed by a cyclic conjugation.

A cyclic conjugation of a word $w$ by $i$ letters is equivalent to conjugate $w$ by the prefix of $w$ of length $i$ and then taking the reduced form. We have seen in Remark \ref{annie} that the total length of the conjugating elements for a corolla $c$ of length at most $n'_i$ and such that $\eta(c)=k_i$ is no more than $(mk_i-1)/2+n'_i-1$ ($m$ is the maximal length of elements of $R$). This means that in constructing such a corolla the sum of the lengths of cyclic conjugations has been at most $(mk_i-1)/2+n'_i-1$.

Let $q\leqslant (mk_i-1)/2+n'_i-1$; suppose that we make $k_i$ cyclically reduced products of defining relators, each followed by a cyclic conjugation of $q_i$ letters for $i=1, \cdots, k_i$ and let $q=q_1 + \cdots + q_{k_i}$. Then $(q_1, \cdots, q_{k_i})$ is a \textit{weak compositions} of $q$ in $k_i$ parts (see \cite{compo}). The number of weak compositions of $q$ in $k_i$ parts is equal to 
               $$\binom {q+k_i-1} {k_i-1}.$$
Finally if $p=|R|$, then there are obviously $p^{k_i}$ products of $k_i$ defining relators.

Therefore if we want to find with the method presented in this thesis all the relators of length at most $n$ and which are products of at most $k$ conjugates of defining relators, then an upper bound for the number of relators which we compute is
   \begin{equation} \label{billie} \sum_{h=1}^{\min(k,n)} \, \, \, \sum_{n_1+\cdots+n_h\leqslant n} \, \, \sum_{j=1}^{h-1} (h-j)n_j \sum_{k_1+\cdots+k_h\leqslant k}    \, \end{equation}
$$\prod_{i=1}^h p^{k_i}  \, \Bigg(\sum_{n'_i=1}^{n_i}   (a+1)a^{(n_i-n'_i)/2-1} \sum_{q=1}^{(mk_i-1)/2+n'_i-1} \binom {q+k_i-1} {k_i-1}\Bigg).$$
     Let us make some observations about this bound. The term $p^{k_i}$ is unrealistically too big because it means that in finding all the corollas $c$ such that $\eta(c)=k_i$ we make all possible cyclically reduced products, whereas by Section \ref{overes} it is enough to take much less (in particular we need to make a product if there is a cancellation).

The term $\sum_{j=1}^{h-1} (h-j)n_j$ too is unrealistically too big because it considers that we make all possible insertions of the words $f_1, \cdots, f_h$ whereas by the result of Section \ref{simpli}, the insertions to be considered are only those such that the resulting word is reduced.

Let us find an upper bound for (\ref{billie}). Since $\binom {q+k_i-1} {k_i-1}<2^{q+k_i-1}$ then

$$\sum_{q=1}^{(mk_i-1)/2+n'_i} \binom {q+k_i-1} {k_i-1} < \sum_{q=1}^{(mk_i-1)/2+n'_i-1} 2^{q+k_i-1}=2^{k_i-1} \, \sum_{q=1}^{(mk_i-1)/2+n'_i-1} 2^{q}=$$

$$=2^{k_i-1} \,  \big(2^{\frac{mk_i-1}{2}+n'_i} -2\big) <2^{k_i-1} \,  \big(2^{\frac{mk_i-1}{2}+n'_i}\big) = \frac{\sqrt{2}^{(m+2)k_i} \, 2^{n'_i}}{2 \sqrt{2}}$$

since we have the equality $\sum_{i=1}^ns^i=(s^{n+1}-s)/(s-1)$ for $s>1$. Thus

$$\sum_{n'_i=1}^{n_i}   (a+1)a^{(n_i-n'_i)/2-1} \, \sum_{q=1}^{(mk_i-1)/2+n'_i-1} \binom {q+k_i-1} {k_i-1}<$$

 \bigskip

 $$\sum_{n'_i=1}^{n_i}   (a+1)a^{(n_i-n'_i)/2-1} \, \frac{\sqrt{2}\,^{(m+2)k_i} \, 2^{n'_i}}{2 \sqrt{2}} = 
   \frac{(a+1)  \, \sqrt{2}\,^{(m+2)k_i}}{2 \sqrt{2} \, a} \sum_{n'_i=1}^{n_i}   a^{(n_i-n'_i)/2} \, 2^{n'_i} =$$
   
   \bigskip

\begin{equation} \label{sandor} = \frac{(a+1)  \, \sqrt{2}\,^{(m+2)k_i} \, \sqrt{a}\, ^{n_i} }{2 \sqrt{2} \, a} \sum_{n'_i=1}^{n_i}   \Big(\frac{2}{\sqrt{a}}\Big)^{n'_i} \end{equation} 
   
   \bigskip

If $s$ is a number greater than 1, we know that 
      $$\sum_{i=1}^n (1/ s)^i=( s^n-1)/ s^n( s-1).$$
 If $ s= s_1/ s_2$ with $ s_1> s_2$ then 
                    $$\sum_{i=1}^n (1/ s)^i=( s_1^n- s_2^n) s_2/ s_1^n( s_1- s_2).$$

 Let $a>4$; thus (\ref{sandor}) is equal to

$$ \frac{(a+1)  \, \sqrt{2}\,^{(m+2)k_i} \, \sqrt{a}\, ^{n_i} }{\sqrt{2} \, a} \, \frac{\sqrt{a}\,^{n_i}-2^{n_i}}{\sqrt{a}\,^{n_i} (\sqrt{a}-2)} = \frac{(a+1) } {\sqrt{2} \, a \, (\sqrt{a}-2)} \, (\sqrt{a}\,^{n_i}-2^{n_i})  \, \sqrt{2}\,^{(m+2)k_i}.$$

 \medskip

 Let us compute now the term
 
 \begin{equation} \label{maris}  \prod_{i=1}^h p^{k_i}  \, \Bigg(\sum_{n'_i=1}^{n_i}   (a+1)a^{(n_i-n'_i)/2-1} \sum_{q=1}^{(mk_i-1)/2+n'_i-1} \binom {q+k_i-1} {k_i-1}\Bigg). \end{equation} 
 
 \medskip
 
 We have that (\ref{maris}) is less than 
  
 $$\frac{(a+1) } {\sqrt{2} \, a \, (\sqrt{a}-2)} \, \prod_{i=1}^h  (\sqrt{a}\,^{n_i}-2^{n_i})  \, p^{k_i}  \, \sqrt{2}\,^{(m+2)k_i}=$$
 
 $$\frac{(a+1) } {\sqrt{2} \, a \, (\sqrt{a}-2)} \, (\sqrt{a}\,^{n_1}-2^{n_1})  \, \cdots  (\sqrt{a}\,^{n_h}-2^{n_h})  \, \big(p  \sqrt{2}\,^{m+2}\big) ^{k_1+\cdots + k_h}=$$

 $$\frac{(a+1) } {\sqrt{2} \, a \, (\sqrt{a}-2)} \, (\sqrt{a}\,^{n_1}-2^{n_1})  \, \cdots  (\sqrt{a}\,^{n_h}-2^{n_h})  \, c^{k_1+\cdots + k_h},$$

 where $c= p\sqrt{2}\,^{m+2}=4p\sqrt{2}\,^m$.

   \medskip

 Let us compute the term

  $$\sum_{k_1+\cdots+k_h\leqslant k}  \, \,  \prod_{i=1}^h   p^{k_i}$$
 
 \begin{equation} \label{mjfan} \Bigg(\sum_{n'_i=1}^{n_i}   (a+1)a^{(n_i-n'_i)/2-1} \, \, \sum_{q=1}^{(mk_i-1)/2+n'_i-1} \binom {q+k_i-1} {k_i-1}\Bigg).\end{equation}

 \medskip
 
 We have that (\ref{mjfan}) is less than 
 
  $$\sum_{k_1+\cdots+k_h\leqslant k}  \, \,  \frac{(a+1) } {\sqrt{2} \, a \, (\sqrt{a}-2)} \, (\sqrt{a}\,^{n_1}-2^{n_1})  \, \cdots  (\sqrt{a}\,^{n_h}-2^{n_h})  \,  c ^{k_1+\cdots + k_h}=$$

 $$ \frac{(a+1) } {\sqrt{2} \, a \, (\sqrt{a}-2)} \, (\sqrt{a}\,^{n_1}-2^{n_1})  \, \cdots  (\sqrt{a}\,^{n_h}-2^{n_h})  \, \sum_{k_1+\cdots+k_h\leqslant k}  \, \, c ^{k_1+\cdots + k_h}.$$
 
 We have that
 
 $$\sum_{k_1+\cdots+k_h\leqslant k}  \, \, c ^{k_1+\cdots + k_h}=\sum_{\alpha=1}^ k \, \sum_{k_1+\cdots+k_h=\alpha}  \, \, c ^{k_1+\cdots + k_h}=$$

 $$=\sum_{\alpha=h}^ k \, \binom {\alpha-1} {h-1}  \, c ^{\alpha} < \binom {k-1} {h-1}  \, \sum_{\alpha=h}^ k \,  c ^{\alpha}=$$

 $$=\binom {k-1} {h-1}  \,  \frac{c ^{k+1} - c^h}{c-1}.$$

Thus (\ref{mjfan}) is less than

    $$ \frac{(a+1) } {\sqrt{2} \, a \, (\sqrt{a}-2)}  \,  (\sqrt{a}\,^{n_1}-2^{n_1})  \, \cdots  (\sqrt{a}\,^{n_h}-2^{n_h})  \, \binom {k-1} {h-1}  \,  \frac{c ^{k+1} - c^h}{c-1}.$$

 \medskip
 
We have that
 
$$\sum_{j=1}^{h-1} (h-j)n_j < (h-1)^2 (n-h+1)$$
 
 since every $n_j$ is less or equal to $n-h+1$ because $n_1+ \cdots + n_h \leqslant n$.
 
 We have also that 
 
$$\sum_{n_1+\cdots+n_h\leqslant n} (\sqrt{a}\,^{n_1}-2^{n_1})  \cdots  (\sqrt{a}\,^{n_h}-2^{n_h})<\sum_{n_1+\cdots+n_h\leqslant n} \sqrt{a}\,^{n_1}  \cdots \sqrt{a}\,^{n_h}=$$
 
 $$\sum_{l=1}^n \, \sum_{n_1+\cdots+n_h= l} \sqrt{a}\,^{n_1+ \cdots +n_h}=\sum_{l=h}^n \,  \binom {l-1} {h-1} \sqrt{a}\,^l<\binom {n-1} {h-1} \frac{\sqrt{a}\,^n - \sqrt{a}\,^h }{\sqrt{a} -1}.$$
 
 \medskip

 Now we compute 
 
 \begin{equation} \label{nabi} \sum_{n_1+\cdots+n_h\leqslant n} \, \, \sum_{j=1}^{h-1} (h-j)n_j \sum_{k_1+\cdots+k_h\leqslant k}    \, \end{equation}
$$\prod_{i=1}^h p^{k_i}  \, \Bigg(\sum_{n'_i=1}^{n_i}   (a+1)a^{(n_i-n'_i)/2-1} \sum_{q=1}^{(mk_i-1)/2+n'_i-1} \binom {q+k_i-1} {k_i-1}\Bigg).$$

We have that (\ref{nabi}) is less than

 $$d (h-1)^2 (n-h+1) \binom {n-1} {h-1}  \, \binom {k-1} {h-1} (\sqrt{a}\,^n - \sqrt{a}\,^h)    \Big(c ^k - c^{h-1}\Big),$$

 where $$d=\frac{(a+1)c } {\sqrt{2} \, a \, (\sqrt{a}-2) (\sqrt{a} -1) (c-1)}.$$
 
This implies that (\ref{billie}) is bounded above by 
 
   $$d \sqrt{a}\,^n c ^k  \sum_{h=2}^{\min(k,n)}   (h-1)^2 (n-h+1) \binom {n-1} {h-1}  \, \binom {k-1} {h-1}.$$

 Using the equality
 
$$\sum_{h=2}^n (h-1)^2 \binom {n-1} {h-1}=\big(n-1+(n-1)^2\big) 2^{n-3}=n (n-1) 2^{n-3}$$
 
 (see \cite{bino}, 6b), we have that (\ref{nabi}) is less than
 
  $$d \sqrt{a}\,^n c ^k  n \sum_{h=2}^n   (h-1)^2 \binom {n-1} {h-1}  \, \sum_{h=2}^k \binom {k-1} {h-1}<$$

 $$d \sqrt{a}\,^n c ^k  n^2 (n-1) 2^{n-3} \, 2^{k-1}=$$
 
  \begin{equation} \label{earlgrey} \frac{d}{16} \, n^2 (n-1)  \sqrt{a}\,^n c ^k \, 2^{n+k}=d' \, n^2 (n-1)  \beta \,^n \gamma^k,\end{equation}
  
  \smallskip
  
  where $\beta=2\sqrt{a}=2\sqrt{2|X|-1}$, \, $\gamma=2c=4p\sqrt{2}\,^m=4\sqrt{2}\,^m |R|$ (we recall that $m$ is the maximal length of the elements of $R$) and 
  
$$d'= \frac{d}{16}=\frac{(a+1)c } {16\sqrt{2} \, a \, (\sqrt{a}-2) (\sqrt{a} -1) (c-1)}\approx \frac{1}{16\sqrt{2} \, a}.$$

For $a <4$ with analogous calculations we prove that (\ref{nabi}) is less than

$$b \, n^2 (n-1)  4^n \gamma^k,$$ 
 
 where
 
$$b\approx \frac{1}{16\sqrt{2} \, (2-\sqrt{a})}.$$

We recall that since $a=2|X|-1$, then $a$ cannot be equal to 4 since it is odd.

We have thus proved that an upper bound for the complexity of our method is a product of a cubic polynomial times an exponential in $n$ and one in $k$. This is much better than (\ref{frot}), the complexity with the method currently known in the literature, which is a product of three exponentials, one in $k$, one in $kn$ and one in $k^2$. 

Let us compute the ratio of the complexity of the classical method over that of ours; for $a>4$ it is 

\begin{equation} \label{rotc} \frac{\big(\frac{a+1}{a-1} \big)^k p^k  a^{mk^2+kn}} {n^2(n-1) \, p^k \,  a^{n/2} \,  4^k \,  2^{mk/2}}. \end{equation}

We can assume that $\big(\frac{a+1}{a-1} \big)^k > n^2(n-1)$ and since $a>4$ then (\ref{rotc}) is greater than

$$\frac{a^{mk^2+kn}}{a^{n/2 + k + mk/2}}=a^{k(mk-m/2-1)+n(k-1/2)}>a^{k(k-c_1)+n(k-c_2)},$$

where $c_1=m/2+1$ and $c_2=1/2$. For $a<4$ we obtain the same result but with  $c_1=m/2+2$ and $c_2=2$.

Moreover (\ref{earlgrey}) is far from being an optimal bound for the complexity of our method, but only a value which is obtainable with not too complicated calculations. Instead (\ref{frot}) is the real complexity of the classical method, that is it is the exact number of words ti be considered with this method to find all the relators of length up to $n$ and area up to $k$.

\section{Recursively defined sets and straight line algorithms} \label{rdsaslp}

Let $U$ be a set and let $\Phi$ be a family of functions $\varphi : U^{n_\varphi} \rightarrow U$ (where $n_\varphi$ is a given non-zero natural number depending on $\varphi$) with codomain $U$ and with domain some Cartesian power of $U$. Let $T\subset U$; we say that \textit{$T$ is $\Phi$-closed} if for every $\varphi \in \Phi$ and for every $t_1, \cdots, t_{n_\varphi} \in T$ we have that $\varphi(t_1, \cdots, t_{n_\varphi}) \in T$.

\begin{definition} \label{first} \rm Let $B\subset U$; the intersection of all $\Phi$-closed subsets of $U$ containing $B$ is called \textit{the subset of $U$ recursively defined by $B$ and $\Phi$}.                                     \end{definition}
                                   
 $B$ is called the \textit{base set} and $\Phi$ \textit{the set of operations}. Since an intersection of $\Phi$-closed sets is still $\Phi$-closed, \textit{the subset of $U$ recursively defined by $B$ and $\Phi$ is the least $\Phi$-closed subset of $U$ containing $B$}.

Let $C$ be the set of the results of the straight line algorithms (Definition \ref{SLA}) relative to $(U, B, \Phi)$. We will now show that $C$ coincides with the subset of $U$ recursively defined by $B$ and $\Phi$.

\begin{proposition}\label{chiusura} $C$ is $\Phi$-closed. \end{proposition}

\begin{proof} Let $\varphi\in \Phi$ and let $c_1, \cdots, c_{n_\varphi} \in C$; we have to prove that 
       $$\varphi(c_1, \cdots, c_{n_\varphi}) \in C.$$
Let $\sigma_1, \cdots, \sigma_{n_\varphi}$ be $SLA$'s computing respectively $c_1, \cdots, c_{n_\varphi}$. The algorithm whose steps are the steps of all the $\sigma_i$ plus a last step equal to the application of $\varphi$ to $c_1, \cdots, c_{n_\varphi}$ is an $SLA$ and its result is $\phi( c_1, \cdots, c_{n_\varphi})$.
                         \end{proof}

Let $b\in B$; $b$ is the result of the $SLA$ with a single step equal to $b$, that is $b \in C$ and therefore $B\subset C$. Since $C$ is $\Phi$-closed and contains $B$, then by Definition \ref{first}, $C$ contains the subset of $U$ recursively defined by $B$ and $\Phi$. The following theorem implies the opposite inclusion.

 \begin{proposition} \label{T} $C$ is contained in any $\Phi$-closed subset of $U$ containing $B$.\end{proposition}

\begin{proof} Let $T$ be a $\Phi$-closed subset of $U$ containing $B$. If $c\in C$ then there exists an $SLA$ $\sigma$ whose result is $c$. We prove that the output of every step of $\sigma$ (and in particular $c$) belongs to $T$.

Since the first step of $\sigma$ is a base step, the first output is an element of $B$, which is contained in $T$. Suppose that the first $k-1$ outputs of $\sigma$ belong to $T$. If the $k$-th is a base step, then its output is an element of $B$ and therefore of $T$; if not, there exist $\varphi \in \Phi$ and $c_1, \cdots, c_{n_\varphi}$ outputs of preceding steps such that the output of the $k$-th step is $\varphi(c_1, \cdots, c_{n_\varphi})$. Since $c_1, \cdots, c_{n_\varphi}$ belong to $T$ by induction hypothesis and since $T$ is $\Phi$-closed, then $\varphi(c_1, \cdots, c_{n_\varphi})\in T$.
              \end{proof}

We have thus proved the following result

\begin{theorem} \label{crucial} Let $C$ be the set of the results of the straight line algorithms relative to $(U, B, \Phi)$. Then $C$ coincides with the subset of $U$ recursively defined by $B$ and $\Phi$.  \end{theorem}

\begin{remark} \label{pr.ind} \rm Proposition \ref{T} is a generalization of the \textit{Principle of mathematical induction} (see \cite{pomi}); it can be stated in the following way:

Let $U$ be a set and let $\Phi$ be a family of functions $\varphi : U^{n_\varphi} \rightarrow U$ (where $n_\varphi$ is a given non-zero natural number depending on $\varphi$) with codomain $U$ and with domain some Cartesian power of $U$. Let $\mathcal{P}$ be a proposition defined on a set $U$, let $B$ be a subset of $U$ and let $\mathcal{P}$ be such that: \begin{itemize}

  \item $\mathcal{P}$ is true for every element of $B$;
  
  \item if $\varphi \in \Phi$ and if $\mathcal{P}$ is true for $c_1, \cdots, c_{n_\varphi}$ then $\mathcal{P}$ is true also for $\varphi(c_1, \cdots, c_{n_\varphi})$.
 
\end{itemize}

Then $\mathcal{P}$ is true for every element of the subset of $U$ recursively defined by $B$ and $\Phi$.

 \end{remark}

\section{A long calculation} \label{alc}
     
  \begin{remark} \label{wonder} \rm Let $x$ and $w$ be cyclically reduced words and let $w$ be a cyclic conjugate of $\pi(u, v)$, where $u$ and $v$ are cyclically reduced words. Let $y:=\pi(x, w)$; we show that there exist words $x'$, $u'$, $v'$ and $y'$ which are cyclic conjugates of $x$, $u$, $v$ and $y$ respectively such that
  
  \begin{itemize}
  
  \item either $y'=\pi(p, u')$ where $p$ is a cyclic conjugate of $\pi(x', v')$
  
  \item or $y'=\pi(p, v')$ where $p$ is a cyclic conjugate of $\pi(x', u')$.

\end{itemize} \end{remark}

\begin{proof} There exist words $w_1$ and $w_2$ such that $w=w_1 w_2$ and $\pi(u, v)=w_2 w_1$. Thus by Remark \ref{shirv}, four cases are possible:

  \begin{enumerate}

  \item $u=u_1 a$, \, $v=a^{-1} t w_2 w_1 t^{-1} u_1^{-1}$. 
  
  \item $u=t w_2 w_4 a$, \, $v=a^{-1} w_3 t^{-1}$ and  $w_4 w_3= w_1$.
  
  \item $u=t w_3 a$, \, $v=a^{-1} w_4 w_1 t^{-1}$ and $w_3 w_4= w_2$.

  \item $u=v_1^{-1} t w_2 w_1  t^{-1} a$, \, $v=a^{-1} v_1$. 
    
\end{enumerate} 

\medskip

Again by Remark \ref{shirv}, three cases are possible:

 \begin{enumerate} [(a)]

  \item $x=x_1 b$, \, $w=b^{-1} s y s^{-1} x_1^{-1}$. 
  
  \item $x=s y_1 b$, \, $w=b^{-1} y_2 s^{-1}$.
  
  \item $x=w_5^{-1} s y  s^{-1} b$, \, $w=b^{-1} w_5$. 
    
\end{enumerate}

\bigskip

Cases (a), (b) and (c) split again. (a) splits in five cases:

 \begin{enumerate} [I]

\item There exist words $b_1$ and $b_2$ such that $b^{-1}=b_2^{-1} b_1^{-1}$ and $w_1=b_2^{-1}$, $w_2=b_1^{-1} s y s^{-1} x_1^{-1}$.

\item There exist words $s_1$ and $s_2$ such that $s=s_1 s_2$ and $w_1=b^{-1} s_1$, $w_2=s_2 \, y \, s_2^{-1} s_1^{-1} x_1^{-1}$.

\item There exist words $y_1$ and $y_2$ such that $y=y_1 y_2$ and $w_1=b^{-1} s y_1$, $w_2=y_2 \, s^{-1} x_1^{-1}$.

\item There exist words $s_1$ and $s_2$ such that $s=s_1 s_2$ and $w_1=b^{-1} s_1 \, s_2 \, y  \, s_2^{-1}$, $w_2= s_1^{-1} x_1^{-1}$.

\item There exist words $x_2$ and $x_3$ such that $x_1^{-1} = x_3^{-1} \, x_2^{-1}$ and $w_1=b^{-1} s \, y \, s^{-1} x_3^{-1}$, $w_2= x_2^{-1}$.

\end{enumerate}

\bigskip

(b) splits in three cases: 

\begin{enumerate} [I]

\item There exist words $b_1$ and $b_2$ such that $b^{-1}=b_2^{-1} b_1^{-1}$ and $w_1=b_2^{-1}$, $w_2=b_1^{-1} y_2 s^{-1}$.

\item There exist words $y_3$ and $y_4$ such that $y_2=y_3 \, y_4$ and $w_1=b^{-1} y_3$, $w_2=y_4 \, s^{-1}$.

\item There exist words $s_1$ and $s_2$ such that $s=s_1 s_2$ and $w_1=b^{-1} y_2  \, s_2^{-1}$, $w_2= s_1^{-1} $.

\end{enumerate}

\bigskip

(c) splits in two cases: 

\begin{enumerate} [I]

\item There exist words $b_1$ and $b_2$ such that $b^{-1}=b_2^{-1} b_1^{-1}$ and $w_1=b_2^{-1}$, $w_2=b_1^{-1} w_5$.

\item There exists a word $w_6$ such that $w_5=w_6 w_2$ and $w_1=b^{-1} w_6$.

\end{enumerate}

Let us examine them all.

\bigskip

\textbf{1 a I.}  

$x=x_1 b_1 b_2$, \, $u=u_1 a$,  \, $v=a^{-1} t b_1^{-1} s y s^{-1} x_1^{-1} b_2^{-1} t^{-1} u_1^{-1}$,  \,  

$w=w_1 w_2=b_2^{-1} b_1^{-1} s y s^{-1} x_1^{-1}$.

 $$y=\pi(x, w)=\pi \big(x_1 b_1 b_2, \, b_2^{-1} b_1^{-1} s y s^{-1} x_1^{-1} \big).$$

$$x=x_1 b_1 b_2 \rightarrow b_1 b_2 x_1,$$ 

$$v=a^{-1} t b_1^{-1} s y s^{-1} x_1^{-1} b_2^{-1} t^{-1} u_1^{-1} \rightarrow x_1^{-1} b_2^{-1} t^{-1} u_1^{-1} a^{-1} t b_1^{-1} s y s^{-1}$$

$$\pi \big(b_1 b_2 x_1, \, x_1^{-1} b_2^{-1} t^{-1} u_1^{-1} a^{-1} t b_1^{-1} s y s^{-1} \big) = b_1 t^{-1} u_1^{-1} a^{-1} t b_1^{-1} s y s^{-1}$$

because $b_1  t^{-1}$ and $s^{-1} b_1$ are subwords of $v^{-1}$.

$$b_1 t^{-1} u_1^{-1} a^{-1} t b_1^{-1} s y s^{-1} \rightarrow  t b_1^{-1} s y s^{-1} b_1 t^{-1} u_1^{-1} a^{-1},$$

$$u=u_1 a \rightarrow a u_1$$

and

$$\pi \big(t b_1^{-1} s y s^{-1} b_1 t^{-1} u_1^{-1} a^{-1}, \, a u_1 \big) = y.$$

\bigskip

\textbf{1 a II.}  

$x=x_1 b$, \, $u=u_1 a$, \, $v=a^{-1} t s_2 \, y \, s_2^{-1} s_1^{-1} x_1^{-1} b^{-1} s_1 t^{-1} u_1^{-1}$,

 $w=w_1 w_2=b^{-1} s_1 s_2 y s_2^{-1} s_1^{-1} x_1^{-1}$

 $$y=\pi(x, w)=\pi \big(x_1 b, \, b^{-1} s_1 s_2 y s_2^{-1} s_1^{-1} x_1^{-1} \big).$$

$$x=x_1 b \rightarrow b x_1,$$ 

$$v=a^{-1} t s_2 \, y \, s_2^{-1} s_1^{-1} x_1^{-1} b^{-1} s_1 t^{-1} u_1^{-1} \rightarrow x_1^{-1} b^{-1} s_1 t^{-1} u_1^{-1} a^{-1} t s_2 \, y \, s_2^{-1} s_1^{-1}$$

$$\pi \big(b x_1, \, x_1^{-1} b^{-1} s_1 t^{-1} u_1^{-1} a^{-1} t s_2 \, y \, s_2^{-1} s_1^{-1} \big) = t^{-1} u_1^{-1} a^{-1} t s_2 \, y \, s_2^{-1}$$

because $s_2^{-1} t^{-1}$ is a subword of $v^{-1}$.

$$t^{-1} u_1^{-1} a^{-1} t s_2 \, y \, s_2^{-1} \rightarrow  t s_2 \, y \, s_2^{-1} t^{-1} u_1^{-1} a^{-1},$$

$$u=u_1 a \rightarrow a u_1$$

and

$$\pi \big(t s_2 \, y \, s_2^{-1} t^{-1} u_1^{-1} a^{-1}, \, a u_1 \big) = y.$$

\bigskip

\textbf{1 a III.}  

$x=x_1 b$, \, $u=u_1 a$,  \, $v=a^{-1} t y_2 \, s^{-1} x_1^{-1} b^{-1} s y_1 t^{-1} u_1^{-1}$,

$w=w_1 w_2=b^{-1} s y_1 y_2 s^{-1} x_1^{-1}$.

 $$y=\pi(x, w)=\pi \big(x_1 b, \, b^{-1} s y_1 y_2 s^{-1} x_1^{-1} \big).$$

$$x=x_1 b \rightarrow b x_1,$$ 

$$v=a^{-1} t y_2 \, s^{-1} x_1^{-1} b^{-1} s y_1 t^{-1} u_1^{-1} \rightarrow x_1^{-1} b^{-1} s y_1 t^{-1} u_1^{-1} a^{-1} t y_2 \, s^{-1}$$

and

$$\pi \big(b x_1, \, x_1^{-1} b^{-1} s y_1 t^{-1} u_1^{-1} a^{-1} t y_2 \, s^{-1} \big) = y_1 t^{-1} u_1^{-1} a^{-1} t y_2$$

because $y_2 y_1$ is a cyclic conjugation of $y$ which is cyclically reduced.

$$y_1 t^{-1} u_1^{-1} a^{-1} t y_2 \rightarrow t y_2 y_1 t^{-1} u_1^{-1} a^{-1},$$

$$u=u_1 a \rightarrow a u_1$$

and

$$\pi \big(t y_2 y_1 t^{-1} u_1^{-1} a^{-1}, \, a u_1 \big) = y_2 y_1$$

and finally

$$y_2 y_1 \rightarrow y.$$

\bigskip

\textbf{1 a IV.}  

 $x=x_1 b$, \, $u=u_1 a$, $v=a^{-1} t s_1^{-1} x_1^{-1} b^{-1} s_1 \, s_2 \, y  \, s_2^{-1} t^{-1} u_1^{-1}$,

 $w=w_1 w_2=b^{-1} s_1 s_2 y s_2^{-1} s_1^{-1} x_1^{-1}$.

 $$y=\pi(x, w)=\pi \big(x_1 b, \, b^{-1} s_1 s_2 y s_2^{-1} s_1^{-1} x_1^{-1} \big).$$

$$x=x_1 b \rightarrow b x_1,$$ 

$$v=a^{-1} t s_1^{-1} x_1^{-1} b^{-1} s_1 \, s_2 \, y  \, s_2^{-1} t^{-1} u_1^{-1} \rightarrow x_1^{-1} b^{-1} s_1 \, s_2 \, y  \, s_2^{-1} t^{-1} u_1^{-1} a^{-1} t s_1^{-1}$$

and

$$\pi \big(b x_1, \, x_1^{-1} b^{-1} s_1 \, s_2 \, y  \, s_2^{-1} t^{-1} u_1^{-1} a^{-1} t s_1^{-1} \big) =  s_2 \, y  \, s_2^{-1} t^{-1} u_1^{-1} a^{-1} t$$

because $t s_2$ is a subword of $v^{-1}$.

$$s_2 \, y  \, s_2^{-1} t^{-1} u_1^{-1} a^{-1} t \rightarrow  t s_2 \, y  \, s_2^{-1} t^{-1} u_1^{-1} a^{-1},$$

$$u=u_1 a \rightarrow a u_1$$

and

$$\pi \big(t s_2 \, y  \, s_2^{-1} t^{-1} u_1^{-1} a^{-1}, \, a u_1 \big) = y.$$

\bigskip

\textbf{1 a V.}  

$x=x_2 x_3 b$, \, $u=u_1 a$,  \, $v=a^{-1} t x_2^{-1} b^{-1} s \, y \, s^{-1} x_3^{-1} t^{-1} u_1^{-1}$,

$w=w_1 w_2=b^{-1} s y s^{-1} x_3^{-1} \, x_2^{-1}$.

 $$y=\pi(x, w)=\pi \big(x_2 x_3 b, \, b^{-1} s y s^{-1} x_3^{-1} \, x_2^{-1} \big).$$

$$x=x_2 x_3 b \rightarrow x_3 b x_2,$$ 

$$v=a^{-1} t x_2^{-1} b^{-1} s \, y \, s^{-1} x_3^{-1} t^{-1} u_1^{-1} \rightarrow x_2^{-1} b^{-1} s \, y \, s^{-1} x_3^{-1} t^{-1} u_1^{-1} a^{-1} t$$

$$\pi \big(x_3 b x_2, \, x_2^{-1} b^{-1} s \, y \, s^{-1} x_3^{-1} t^{-1} u_1^{-1} a^{-1} t \big) = x_3 s \, y \, s^{-1} x_3^{-1} t^{-1} u_1^{-1} a^{-1} t$$

because $x_3  s$ and $t x_3$ are subwords of $v^{-1}$.

$$x_3 s \, y \, s^{-1} x_3^{-1} t^{-1} u_1^{-1} a^{-1} t \rightarrow  t x_3 s \, y \, s^{-1} x_3^{-1} t^{-1} u_1^{-1} a^{-1},$$

$$u=u_1 a \rightarrow a u_1$$

and

$$\pi \big( t x_3 s \, y \, s^{-1} x_3^{-1} t^{-1} u_1^{-1} a^{-1}, \, a u_1 \big) = y.$$

\bigskip

\textbf{1 b I.} 

$x=s y_1 b_1 b_2$, \, $u=u_1 a$, \, $v=a^{-1} t b_1^{-1} y_2 s^{-1} b_2^{-1} t^{-1} u_1^{-1}$,

$w=w_1 w_2=b_2^{-1} b_1^{-1} y_2 s^{-1}$.

 $$y=\pi(x, w)=\pi \big(s y_1 b_1 b_2, \, b_2^{-1} b_1^{-1} y_2 s^{-1} \big).$$

$$x=s y_1 b_1 b_2 \rightarrow y_1 b_1 b_2 s,$$ 

$$v=a^{-1} t b_1^{-1} y_2 s^{-1} b_2^{-1} t^{-1} u_1^{-1} \rightarrow  s^{-1} b_2^{-1} t^{-1} u_1^{-1} a^{-1} t b_1^{-1} y_2$$

$$\pi \big(y_1 b_1 b_2 s, \, s^{-1} b_2^{-1} t^{-1} u_1^{-1} a^{-1} t b_1^{-1} y_2 \big) = y_1 b_1 t^{-1} u_1^{-1} a^{-1} t b_1^{-1} y_2$$

because $b_1  t^{-1}$ is a subword of $v^{-1}$ and $y_2 y_1$ is a cyclic conjugation of $y$ which is cyclically reduced.

$$y_1 b_1 t^{-1} u_1^{-1} a^{-1} t b_1^{-1} y_2 \rightarrow  t b_1^{-1} y_2 y_1 b_1 t^{-1} u_1^{-1} a^{-1},$$

$$u=u_1 a \rightarrow a u_1$$

and

$$\pi \big(t b_1^{-1} y_2 y_1 b_1 t^{-1} u_1^{-1} a^{-1}, \, a u_1 \big) = y_2 y_1$$

and finally

$$y_2 y_1 \rightarrow y.$$

\bigskip

\textbf{1 b II.} 

$x=s y_1 b$, \, $u=u_1 a$, \, $v=a^{-1} t y_4 \, s^{-1} b^{-1} y_3 t^{-1} u_1^{-1}$,

$w=w_1 w_2 =b^{-1} y_3 \, y_4 s^{-1}$.

 $$y=\pi(x, w)=\pi \big(s y_1 b, \, b^{-1} y_3 \, y_4 s^{-1} \big).$$

$$x=s y_1 b \rightarrow y_1 b s,$$ 

$$v=a^{-1} t y_4 \, s^{-1} b^{-1} y_3 t^{-1} u_1^{-1} \rightarrow s^{-1} b^{-1} y_3 t^{-1} u_1^{-1} a^{-1} t y_4$$

$$\pi \big(y_1 b s, \, s^{-1} b^{-1} y_3 t^{-1} u_1^{-1} a^{-1} t y_4 \big) = y_1 y_3 t^{-1} u_1^{-1} a^{-1} t y_4$$

because $y_1 y_3$ and $y_4 y_1$ are subwords of $y$ and of a cyclic conjugate of $y$.

$$y_1 y_3 t^{-1} u_1^{-1} a^{-1} t y_4 \rightarrow  t y_4 y_1 y_3 t^{-1} u_1^{-1} a^{-1},$$

$$u=u_1 a \rightarrow a u_1$$

and

$$\pi \big(t y_4 y_1 y_3 t^{-1} u_1^{-1} a^{-1}, \, a u_1 \big) = y_4 y_1 y_3$$

and finally

$$y_4 y_1 y_3 \rightarrow y.$$

\bigskip

\textbf{1 b III.} 

$x=s_1 s_2 y_1 b$, \, $u=u_1 a$, \, $v=a^{-1} t s_1^{-1} b^{-1} y_2  \, s_2^{-1} t^{-1} u_1^{-1}$,

$w=w_1 w_2=b^{-1} y_2 s_2^{-1} s_1^{-1}$.

 $$y=\pi(x, w)=\pi \big(s_1 s_2 y_1 b, \, b^{-1} y_2 s_2^{-1} s_1^{-1} \big).$$

$$x=s_1 s_2 y_1 b \rightarrow s_2 y_1 b s_1,$$ 

$$v=a^{-1} t s_1^{-1} b^{-1} y_2  \, s_2^{-1} t^{-1} u_1^{-1} \rightarrow s_1^{-1} b^{-1} y_2  \, s_2^{-1} t^{-1} u_1^{-1} a^{-1} t$$

$$\pi \big(s_2 y_1 b s_1, \, s_1^{-1} b^{-1} y_2  \, s_2^{-1} t^{-1} u_1^{-1} a^{-1} t \big) = s_2 y_1  y_2  \, s_2^{-1} t^{-1} u_1^{-1} a^{-1} t$$

because $y_1  y_2=y$ is reduced and $t s_2$ is a subword of $v^{-1}$.

$$s_2 y  \, s_2^{-1} t^{-1} u_1^{-1} a^{-1} t \rightarrow  t s_2 y  \, s_2^{-1} t^{-1} u_1^{-1} a^{-1},$$

$$u=u_1 a \rightarrow a u_1$$

and

$$\pi \big(t s_2 y  \, s_2^{-1} t^{-1} u_1^{-1} a^{-1}, \, a u_1 \big) = y.$$

\bigskip

\textbf{1 c I.} 

$x=w_5^{-1} s y  s^{-1} b_1 b_2$,  \, $u=u_1 a$,  \, $v=a^{-1} t b_1^{-1} w_5 b_2^{-1} t^{-1} u_1^{-1}$,

$w=w_1 w_2=b_2^{-1} b_1^{-1} w_1$.

 $$y=\pi(x, w)=\pi \big(w_5^{-1} s y  s^{-1} b_1 b_2, \, b_2^{-1} b_1^{-1} w_1 \big).$$

$$x=w_5^{-1} s y  s^{-1} b_1 b_2 \rightarrow s y  s^{-1} b_1 b_2 w_1^{-1},$$ 

$$v=a^{-1} t b_1^{-1} w_5 b_2^{-1} t^{-1} u_1^{-1} \rightarrow w_1 b_2^{-1} t^{-1} u_1^{-1} a^{-1} t b_1^{-1}$$

and

$$\pi \big(s y  s^{-1} b_1 b_2 w_5^{-1}, \, w_5 b_2^{-1} t^{-1} u_1^{-1} a^{-1} t b_1^{-1} \big) = s y  s^{-1} b_1 t^{-1} u_1^{-1} a^{-1} t b_1^{-1}$$

because $b_1  t^{-1}$ and $b_1 s$ are subwords of $v^{-1}$ and $x^{-1}$ respectively.

$$s y  s^{-1} b_1 t^{-1} u_1^{-1} a^{-1} t b_1^{-1} \rightarrow t b_1^{-1} s y  s^{-1} b_1 t^{-1} u_1^{-1} a^{-1},$$

$$u=u_1 a \rightarrow a u_1$$

and

$$\pi \big(t b_1^{-1} s y  s^{-1} b_1 t^{-1} u_1^{-1} a^{-1}, \, a u_1 \big) = y.$$

\bigskip

\textbf{1 c II.} 

$x=w_2^{-1} w_6^{-1} s y  s^{-1} b$, \, $u=u_1 a$, \,$v=a^{-1} t w_2 b^{-1} w_6 t^{-1} u_1^{-1}$,

$w=w_1 w_2=b^{-1} w_6 w_2$.

 $$y=\pi(x, w)=\pi \big(w_2^{-1} w_6^{-1} s y  s^{-1} b, \, b^{-1} w_6 w_2 \big).$$

$$x=w_2^{-1} w_6^{-1} s y  s^{-1} b \rightarrow w_6^{-1} s y  s^{-1} b w_2^{-1},$$ 

$$v=a^{-1} t w_2 b^{-1} w_6 t^{-1} u_1^{-1} \rightarrow w_2 b^{-1} w_6 t^{-1} u_1^{-1} a^{-1} t$$

and

$$\pi \big(w_6^{-1} s y  s^{-1} b w_2^{-1}, \, w_2 b^{-1} w_6 t^{-1} u_1^{-1} a^{-1} t \big) = w_6^{-1} s y  s^{-1}  w_6 t^{-1} u_1^{-1} a^{-1} t$$

because $s^{-1}  w_6$ and $t w_6^{-1}$ are subwords of $x^{-1}$ and $v^{-1}$ respectively.

$$w_6^{-1} s y  s^{-1}  w_6 t^{-1} u_1^{-1} a^{-1} t \rightarrow  t w_6^{-1} s y  s^{-1}  w_6 t^{-1} u_1^{-1} a^{-1},$$

$$u=u_1 a \rightarrow a u_1$$

and

$$\pi \big(t w_6^{-1} s y  s^{-1}  w_6 t^{-1} u_1^{-1} a^{-1}, \, a u_1 \big) = y.$$

\bigskip

\textbf{2 a I.} 

$x=x_1 b_1 w_3^{-1} w_4^{-1}$, \, $u=t b_1^{-1} s y s^{-1} x_1^{-1} w_4 a$, \, $v=a^{-1} w_3 t^{-1},$

$w_4 w_3= b_2^{-1}$, \, $w=w_1 w_2=w_4 w_3 b_1^{-1} s y s^{-1} x_1^{-1}$.

 $$y=\pi(x, w)=\pi \big(x_1 b_1 b_2, \, b_2^{-1} b_1^{-1} s y s^{-1} x_1^{-1}\big).$$

 $$x=x_1 b_1 w_3^{-1} w_4^{-1} \rightarrow b_1 w_3^{-1} w_4^{-1} x_1,$$ 

$$u= t b_1^{-1} s y s^{-1} x_1^{-1} w_4 a \rightarrow x_1^{-1} w_4 a t b_1^{-1} s y s^{-1}$$

and

$$\pi \big(b_1 w_3^{-1} w_4^{-1} x_1, \, x_1^{-1} w_4 a t b_1^{-1} s y s^{-1} \big) = b_1 w_3^{-1}  a t b_1^{-1} s y s^{-1}$$

because $w_3^{-1} a$ and $s^{-1} b_1$ are subwords of $v^{-1}$ and $w^{-1}$ respectively.

$$b_1 w_3^{-1}  a t b_1^{-1} s y s^{-1} \rightarrow b_1^{-1} s y s^{-1} b_1 w_3^{-1}  a t,$$

$v=a^{-1} w_3 t^{-1} \rightarrow t^{-1} a^{-1} w_3,$

and

$$\pi \big(b_1^{-1} s y s^{-1} b_1 w_3^{-1}  a t, \, t^{-1} a^{-1} w_3 \big) = y.$$

\bigskip

\textbf{2 a II.} 

$x=x_1 b$, \, $u=t s_2 \, y \, s_2^{-1} s_1^{-1} x_1^{-1} w_4 a$, \, $v=a^{-1} w_3 t^{-1}$ 

 $w_4 w_3= b^{-1} s_1$, \, $w=w_4 w_3 w_2=b^{-1} s_1 s_2 y s_2^{-1} s_2^{-1} x_1^{-1}$.

  $$y=\pi(x, w)=\pi \big(x_1 b, \, b^{-1} s_1 s_2 y s_2^{-1} s_2^{-1} x_1^{-1} \big).$$

 $$x=x_1 b \rightarrow b x_1,$$ 

$$u= t s_2 \, y \, s_2^{-1} s_1^{-1} x_1^{-1} w_4 a \rightarrow  x_1^{-1} w_4 a t s_2 \, y \, s_2^{-1} s_1^{-1} \equiv x_1^{-1} b^{-1} s_1 w_3^{-1} a t s_2 \, y \, s_2^{-1} s_1^{-1}$$

since $w_4 \equiv b^{-1} s_1 w_3^{-1}$. We observe that the last word of the preceding 

line is not necessarily reduced.

$$\pi \big(b x_1, \, x_1^{-1} b^{-1} s_1 w_3^{-1} a t s_2 \, y \, s_2^{-1} s_1^{-1} \big) = w_3^{-1} a t s_2 \, y \, s_2^{-1}$$

because $w_3^{-1} a$ is a subword of $v^{-1}$; moreover $w_3 s_2$ is reduced because the last letter of $w_3$ is equal to the last of $s_1$ which is not opposite to the first of $s_2$.

$$w_3^{-1} a t s_2 \, y \, s_2^{-1} \rightarrow s_2 \, y \, s_2^{-1} w_3^{-1} a t,$$

$$v=a^{-1} w_3 t^{-1} \rightarrow t^{-1} a^{-1} w_3$$

and

$$\pi \big(s_2 \, y \, s_2^{-1} w_3^{-1} a t, \, t^{-1} a^{-1} w_3 \big) = y.$$

\bigskip

\textbf{2 a III.}

$x=x_1 b$, \, $u=t y_2 \, s^{-1} x_1^{-1} w_4 a$, \, $v=a^{-1} w_3 t^{-1}$ and  $w_4 w_3= b^{-1} s y_1$

$w=w_4 w_3 w_2=b^{-1} s y s^{-1} x_1^{-1}$

 $$y=\pi(x, w)=\pi \big(x_1 b, \, b^{-1} s y s^{-1} x_1^{-1} \big).$$

 $$x=x_1 b \rightarrow b x_1,$$ 

$$u= t y_2 \, s^{-1} x_1^{-1} w_4 a \rightarrow  x_1^{-1} w_4 a t y_2 \, s^{-1} \equiv x_1^{-1} b^{-1} s y_1 w_3^{-1} a t y_2 \, s^{-1}$$

since $w_4 \equiv b^{-1} s y_1 w_3^{-1}$. The last word of the preceding line is not 

necessarily reduced.

$$\pi \big(b x_1, \, x_1^{-1} b^{-1} s y_1 w_3^{-1} a t y_2 \, s^{-1} \big) = y_1 w_3^{-1} a t y_2$$

because $y_2 y_1$ is reduced being a cyclic conjugate of $y$.

$$y_1 w_3^{-1} a t y_2 \rightarrow y_2 y_1 w_3^{-1} a t,$$

$$v=a^{-1} w_3 t^{-1} \rightarrow t^{-1} a^{-1} w_3$$

and

$$\pi \big(y_2 y_1 w_3^{-1} a t, \, t^{-1} a^{-1} w_3 \big) = y_2 y_1.$$

Finally 

$$y_2 y_1 \rightarrow y.$$

\bigskip

\textbf{2 a IV.} 

$x=x_1 b$, \,  $u=t s_1^{-1} x_1^{-1} w_4 a$, \, $v=a^{-1} w_3 t^{-1},$

$w_4 w_3= b^{-1} s_1 \, s_2 \, y  \, s_2^{-1}$, \, $w=w_4 w_3 w_2=b^{-1} s_1 s_2 y s_2^{-1}  s_1^{-1} x_1^{-1}$.

  $$y=\pi(x, w)=\pi \big(x_1 b, \, b^{-1} s_1 s_2 y s_2^{-1}  s_1^{-1} x_1^{-1} \big).$$

 $$x=x_1 b \rightarrow b x_1,$$ 

$$u= t s_1^{-1} x_1^{-1} w_4 a \rightarrow  x_1^{-1} w_4 a t s_1^{-1} \equiv x_1^{-1} b^{-1} s_1 \, s_2 \, y  \, s_2^{-1} w_3^{-1} a t s_1^{-1}$$

since $w_4 \equiv b^{-1} s_1 \, s_2 \, y  \, s_2^{-1} w_3^{-1}$. The last word of the preceding line is not necessarily reduced.

$$\pi \big(b x_1, \, x_1^{-1} b^{-1} s_1 \, s_2 \, y  \, s_2^{-1} w_3^{-1} a t s_1^{-1} \big) =   s_2 \, y  \, s_2^{-1} w_3^{-1} a t$$

since $t w_3^{-1}$ is reduced being a subword of $v^{-1}$, thus $t s_2$ is reduced because the first letter of $w_3^{-1}$ is equal to the first of $s_2$.

$$v=a^{-1} w_3 t^{-1} \rightarrow t^{-1} a^{-1} w_3$$

and

$$\pi \big(s_2 \, y  \, s_2^{-1} w_3^{-1} a t, \, t^{-1} a^{-1} w_3 \big) =  y.$$

\bigskip

\textbf{2 a V.} 

$x=w_2^{-1} x_3 b$, \, $u=t w_2 w_4 a$, \, $v=a^{-1} w_3 t^{-1}$,

$w_4 w_3= b^{-1} s \, y \, s^{-1} x_3^{-1}$, \, $w=w_4 w_3 w_2=b^{-1} s y s^{-1} x_3^{-1} \, x_2^{-1}$.

  $$y=\pi(x, w)=\pi \big(w_2^{-1} x_3 b, \, b^{-1} s y s^{-1} x_3^{-1} \, x_2^{-1} \big).$$

 $$x=w_2^{-1} x_3 b \rightarrow x_3 b w_2^{-1},$$ 

$$u= t w_2 w_4 a \rightarrow  w_2 w_4 a t \equiv w_2 b^{-1} s \, y \, s^{-1} x_3^{-1} w_3^{-1} a t$$

since $w_4 \equiv b^{-1} s \, y \, s^{-1} x_3^{-1} w_3^{-1}$. The last word of the preceding line is not necessarily reduced.

$$\pi \big(x_3 b w_2^{-1}, \, w_2 b^{-1} s \, y \, s^{-1} x_3^{-1} w_3^{-1} a t \big) = x_3 s \, y \, s^{-1} x_3^{-1} w_3^{-1} a t$$

because $x_3 s$ is reduced being a subword of $w^{-1}$; moreover $t w_3$ is reduced being a subword of $v^{-1}$, thus $t x_3$ is reduced because the first letter of $x_3$ is equal to the first of $v$.

$$v=a^{-1} w_3 t^{-1} \rightarrow t^{-1} a^{-1} w_3$$

and

$$\pi \big(x_3 s \, y \, s^{-1} x_3^{-1} w_3^{-1} a t, \, t^{-1} a^{-1} w_3 \big) =  y.$$

\bigskip

\textbf{2 b I.} 

$x=s y_1 b_1 b_2$, \, $u=t b_1^{-1} y_2 s^{-1} w_4 a$, \, $v=a^{-1} w_3 t^{-1}$.

$w_4 w_3= w_1=b_2^{-1}$, \, $w=b_2^{-1} b_1^{-1} y_2 s^{-1}$.

 $$y=\pi(x, w)=\pi \big(s y_1 b_1 b_2, \, b_2^{-1} b_1^{-1} y_2 s^{-1} \big).$$

 $$x=s y_1 b_1 b_2 \rightarrow b_2 s y_1 b_1,$$ 

$$u= t b_1^{-1} y_2 s^{-1} w_4 a \rightarrow  b_1^{-1} y_2 s^{-1} w_4 a t.$$

$$\pi \big(b_2 s y_1 b_1, \, b_1^{-1} y_2 s^{-1} w_4 a t \big) = b_2 s y_1 y_2 s^{-1} w_4 a t=b_2 s y s^{-1} w_4 a t$$

since $t b_2$ is reduced because $t w_3^{-1}$ is reduced (being a subword of $v^{-1}$) and the first letter of $w_3^{-1}$ is equal to the first of $b_2$.

$$v=a^{-1} w_3 t^{-1} \rightarrow t^{-1} a^{-1} w_3$$

and

$$\pi \big(b_2 s y s^{-1} w_4 a t, \, t^{-1} a^{-1} w_3 \big) =  b_2 s y s^{-1} w_4  w_3=y$$

since

$w_4 w_3=b_2^{-1}$.

\bigskip

\textbf{2 b II.} 

$x=s y_1 b$, \,  $u=t y_4 \, s^{-1} w_4 a$, \, $v=a^{-1} w_3 t^{-1}$, 

$w_4 w_3= b^{-1} y_3$, \,  $w=b^{-1} y_3 \, y_4 s^{-1}$

 $$y=\pi(x, w)=\pi \big(s y_1 b, \, b^{-1} y_3 \, y_4 s^{-1} \big).$$

 $$x=s y_1 b \rightarrow y_1 b s,$$ 

$$u= t y_4 \, s^{-1} w_4 a \rightarrow s^{-1} w_4 a t y_4  \equiv s^{-1} b^{-1} y_3 w_3^{-1} a t y_4$$

since $w_4 \equiv b^{-1} y_3 w_3^{-1}$. The last word of the preceding line is not necessarily reduced.

$$\pi \big(y_1 b s, \, s^{-1} b^{-1} y_3 w_3^{-1} a t y_4 \big)  \equiv  y_1 y_3 w_3^{-1} a t y_4.$$

$$y_1 y_3 w_3^{-1} a t y_4 \rightarrow y_4 y_1 y_3 w_3^{-1} a t$$

$$v=a^{-1} w_3 t^{-1} \rightarrow t^{-1} a^{-1} w_3$$

and

$$\pi \big(y_4 y_1 y_3 w_3^{-1} a t, \, t^{-1} a^{-1} w_3 \big) =  y_4 y_1 y_3.$$

Finally

$y_4 y_1 y_3 \rightarrow y$.

\bigskip

\textbf{2 b III.} 

$x=s_1 s_2 y_1 b$, \,  $u=t s_1^{-1} w_4 a$, \, $v=a^{-1} w_3 t^{-1}$,

$w_4 w_3= b^{-1} y_2  \, s_2^{-1}$, \, $w=b^{-1} y_2 s^{-1}$.

 $$y=\pi(x, w)=\pi \big(s_1 s_2 y_1 b, \, b^{-1} y_2 s^{-1} \big).$$

 $$x=s_1 s_2 y_1 b \rightarrow s_2 y_1 b s_1,$$ 

$$u= t s_1^{-1} w_4 a \rightarrow  s_1^{-1} w_4 a t \equiv s_1^{-1} b^{-1} y_2  \, s_2^{-1} w_3^{-1} a t$$

since $w_4 \equiv b^{-1} y_2  \, s_2^{-1} w_3^{-1}$. The last word of the preceding line is not necessarily reduced.

$$\pi \big(s_2 y_1 b s_1, \, s_1^{-1} b^{-1} y_2  \, s_2^{-1} w_3^{-1} a t \big) = s_2 y_1  y_2  \, s_2^{-1} w_3^{-1} a t=s_2 y  \, s_2^{-1} w_3^{-1} a t$$

because $t w_3^{-1}$ is reduced being a subword of $v^{-1}$, thus $t s_2$ is reduced since the first letter of $s_2$ is equal to the first of $w_3^{-1}$.

$$v=a^{-1} w_3 t^{-1} \rightarrow t^{-1} a^{-1} w_3$$

and

$$\pi \big(s_2 y  \, s_2^{-1} w_3^{-1} a t, \, t^{-1} a^{-1} w_3 \big) =  y.$$

\bigskip

\textbf{2 c I.} 

$x=w_5^{-1} s y  s^{-1} b_1 b_2$, \, $u=t b_1^{-1} w_5 w_4 a$, \, $v=a^{-1} w_3 t^{-1}$,

$w_4 w_3= b_2^{-1}$, \, $w=b_2^{-1} b_1^{-1} w_5$

 $$y=\pi(x, w)=\pi \big(w_5^{-1} s y  s^{-1} b_1 b_2, \, b_2^{-1} b_1^{-1} w_5 \big).$$

 $$x=w_5^{-1} s y  s^{-1} b_1 b_2 \rightarrow s y  s^{-1} b_1 b_2 w_5^{-1},$$ 

$$u= t b_1^{-1} w_5 w_4 a \rightarrow  w_5 w_4 a t b_1^{-1} \equiv w_5 b_2^{-1} w_3^{-1} a t b_1^{-1}$$

since $w_4 \equiv b_2^{-1} w_3^{-1}$. The last word of the preceding line is not necessarily reduced.

$$\pi \big(s y  s^{-1} b_1 b_2 w_5^{-1}, \, w_5 b_2^{-1} w_3^{-1} a t b_1^{-1} \big) = s y  s^{-1} b_1 w_3^{-1} a t b_1^{-1}$$

because $b_1 s$ is reduced being a subword of $w^{-1}$; moreover $w_3^{-1} w_4^{-1}=b_2$, thus $b_1 w_3^{-1}$ is reduced.

$$s y  s^{-1} b_1 w_3^{-1} a t b_1^{-1} \rightarrow b_1^{-1} s y  s^{-1} b_1 w_3^{-1} a t$$

$$v=a^{-1} w_3 t^{-1} \rightarrow t^{-1} a^{-1} w_3$$

and

$$\pi \big(b_1^{-1} s y  s^{-1} b_1 w_3^{-1} a t, \, t^{-1} a^{-1} w_3 \big) =  y.$$

\bigskip

\textbf{2 c II.} 

$x= w_2^{-1} w_6^{-1} s y  s^{-1} b$, \, $u=t w_2 w_4 a$, \, $v=a^{-1} w_3 t^{-1}$,

$w_4 w_3= b^{-1} w_6$, \, $w=b^{-1} w_6 w_2$.

There exists a word $w_6$ such that $w_5=w_6 w_2$ and $w_1=b^{-1} w_6$.

 $$y=\pi(x, w)=\pi \big(w_2^{-1} w_6^{-1} s y  s^{-1} b, \,b^{-1} w_6 w_2\big).$$

 $$x=w_2^{-1} w_6^{-1} s y  s^{-1} b \rightarrow w_6^{-1} s y  s^{-1} b w_2^{-1},$$ 

$$u= t w_2 w_4 a \rightarrow  w_2 w_4 a  t\equiv w_2 b^{-1} w_6 w_3^{-1} a t$$

since $w_4 \equiv b^{-1} w_6 w_3^{-1}$. The last word of the preceding line is not necessarily reduced.

$$\pi \big(w_6^{-1} s y  s^{-1} b w_2^{-1}, \, w_2 b^{-1} w_6 w_3^{-1} a t \big) \equiv w_6^{-1} s y  s^{-1}  w_6 w_3^{-1} a t.$$

$$v=a^{-1} w_3 t^{-1} \rightarrow t^{-1} a^{-1} w_3$$

and

$$\pi \big(w_6^{-1} s y  s^{-1}  w_6 w_3^{-1} a t, \, t^{-1} a^{-1} w_3 \big) =  y.$$

\bigskip

Finally for case 3 the calculations are analogous to case 2 and for case 4 to case 1. 
                                  \end{proof}

\begin{proposition} \label{amex} Let $\sigma$ be an $\overline{R}$-SLA; then there exists a $C$-SLA $\sigma'$ with the same set of base elements and result of $\sigma$ and such that $\eta(\sigma')=\eta(\sigma)$. In particular $\overline{R}_k=C_k$.
   \end{proposition}

  \begin{proof} Let $n$ be the number of steps of $\sigma$ which cannot be steps of a $C$-SLA, that is steps of the form $\pi(x, w)$ where $x$ and $w$ are not 1-corollas. We prove the claim by induction on $n$.
  
  Let $n=1$ and let $\pi(x, w)$ be the step such that $x$ and $w$ are not 1-corollas. Let $\sigma_1$ and $\sigma_2$ be the pSLsA's computing $x$ and $w$. Since $n=1$ then $\sigma_1$ and $\sigma_2$ are SLA's in $C$. This means in particular that $w$ is a cyclic conjugate of $\pi(u, v)$ where either $u$ or $v$ or both are 1-corollas. Let $\eta(\sigma_2)=k$; by Remark \ref{wonder}, $\pi(x, w)$ is equal to $\pi(p, u')$ or $\pi(p, v')$ where $p$ is a cyclic conjugate of $\pi(x', v')$ or of $\pi(x', u')$ respectively and $v'$ or respectively $u'$ is a 1-corolla. Thus if $\sigma_2$ is the SLsA of $\sigma_2$ computing $u$ or $v$ respectively, then $\eta(\sigma'_2)=k-1$. By iterating a finite number of times we have the claim.
  
  If $n>1$ then we repeat the proceeding finitely many times and we obtain the claim.    \end{proof}

Addresses: 

\textit{Institut de Mathématiques de Jussieu}

\textit{175, rue du Chevaleret}

\textit{75013 Paris - France}        

and

\textit{Dipartimento di Matematica dell'Università di Palermo}

\textit{Via Archirafi 34}

\textit{90123 Palermo - Italy}

\smallskip \smallskip

e-mail: \textsf{carmelovaccaro@yahoo.it}


\begin{thebibliography}{99}

\addcontentsline{toc}{section}{\textbf{Bibliography}}


\bibitem{reclo} http://en.wikipedia.org/wiki/Binary\_relation

\bibitem{bino} http://en.wikipedia.org/wiki/Binomial\_coefficient


\bibitem{compo} http://en.wikipedia.org/wiki/Composition\_(number\_theory)



\bibitem{pomi} http://en.wikipedia.org/wiki/Principle\_of\_mathematical\_induction


\bibitem{adjunc} http://planetmath.org/encyclopedia/AdjunctionSpace.html



\bibitem{zorn} http://planetmath.org/encyclopedia/ZornsLemma.html

\bibitem{Birg} J.-C. Birget, Infinite string rewrite systems and complexity, \textit{J. Symbolic Comput.} \textbf{25}, (1998),  no. 6, 759--793. 


\bibitem{Brids} M. R. Bridson, The geometry of the Word Problem, in M. R. Bridson and S. M. Salamon, eds., \emph{Invitations to Geometry and Topology}, Oxford Univ. Press, 2002, pp. 33--94. 

\bibitem{BriHae} M. R. Bridson, A. Haefliger, \textit{Metric Spaces of Non-Positive Curvature}, Springer-Verlag, 1999. 

\bibitem{BroHue} R. Brown, J. Huebschmann, Identities among relations in \emph{Low-dimensional topology (Bangor, 1979)}, London Math. Soc. Lecture Note Ser., vol. 48, Cambridge Univ. Press, 1982, pp. 153--202. 


\bibitem{CH} D. J. Collins, J. Huebschmann, Spherical Diagrams and Identities Among Relations, \textit{Math. Ann.} \textbf{261}, (1982), no. 2, 155--183.


\bibitem{Carb} A. Carbone, Group cancellation and resolution, \textit{Studia Logica} \textbf{82}, (2006), no. 1, 75--95.   

\bibitem{CDP} M. Coornaert, T. Delzant, A. Papadopoulos, \emph{Géométrie et théorie des groupes : les groupes hyperboliques de Gromov}, Springer-Verlag, 1990.


\bibitem{EC} R. L. Epstein, W. A. Carnielli, \emph{Computability. Computable Functions, Logic, and the Foundations of Mathematics}, Wadsworth Publishing, 1999.


\bibitem{Gerst} S. M. Gersten, Dehn functions and $l\sb 1$-norms of finite presentations, in \textit{Algorithms and classification in combinatorial group theory (Berkeley, CA, 1989)}, Math. Sci. Res. Inst. Publ., 23, Springer, 1992, pp. 195--224. 


\bibitem{GI} R. I. Grigorchuk, S. V. Ivanov, On Dehn functions of infinite presentations of groups, \textit{Geom. Funct. Anal.} \textbf{18}, (2009), no. 6, 1841--1874.

\bibitem{Grom} M. Gromov, Hyperbolic groups, in \emph{Essays in group theory}, Springer, 1987, pp. 75--263.

\bibitem{H} A. Hatcher, \emph{Algebraic Topology}, Cambridge University Press, 2002 or http://www.math.cornell.edu/$\sim$hatcher/AT/ATpage.html

\bibitem{Holt} D. F. Holt, B. Eick, E. O'Brien, \emph{A Handbook of computational group theory}, Chapman \& Hall/CRC, 2005.


\bibitem{LS} R. C. Lyndon, P. E. Schupp, \emph{Combinatorial Group Theory}, Springer-Verlag, 1977.

\bibitem{MO}  K. Madlener,  F. Otto, Pseudonatural algorithms for the word problem for finitely presented monoids and groups, \textit{J. Symbolic Comput.} \textbf{1}, (1985), no. 4, 383--418.



\bibitem{MKS} W. Magnus, A. Karrass, D. Solitar, \emph{Combinatorial Group Theory}, Dover Publications, 1976.





\bibitem{Miller} C. F. Miller III, Decision problems for groups - survey and reflections, in \textit{Algorithms and classification in combinatorial group theory (Berkeley, CA, 1989)}, 1--59, Math. Sci. Res. Inst. Publ., 23, Springer, 1992.


\bibitem{Short} H. Short, Diagrams and groups, in \emph{The Geometry of the Word Problem for Finitely Generated Groups}, Birkhäuser, 2007.

\bibitem{Shen} A. Shen, N. K. Vereshchagin, \textit{Computable functions}, American Mathematical Society, 2003.


\bibitem{Sims} C. C. Sims, \emph{Computation with finitely presented groups}, Cambridge University Press, 1994.






\bibitem{VaccDecisProb} C. Vaccaro, Decision problems for finite and infinite presentations of groups and monoids, https://arxiv.org/abs/1001.3981

\bibitem{Whitn} H. Whitney, Non-separable and planar graphs, \textit{Trans. Amer. Math. Soc.} \textbf{34}, (1932), no. 2, 339--362.

\end{thebibliography}
\end{document}